\documentclass[12pt,reqno]{amsart}

\usepackage{amscd,amssymb,epsfig,float,url}
\usepackage[all]{xy}

\newcommand{\A}{{\mathbb A}}

\newcommand{\Z}{{\mathbb Z}}
\newcommand{\Q}{{\mathbb Q}}
\newcommand{\C}{{\mathbb C}}
\newcommand{\R}{{\mathbb R}}
\renewcommand{\P}{{\mathbb P}}
\renewcommand{\H}{{\mathbb H}}

\newcommand{\BB}{{\mathcal B}}
\newcommand{\DD}{{\mathcal D}}

\newcommand{\KK}{{\mathcal K}}
\newcommand{\LL}{{\mathcal L}}

\newcommand{\OO}{{\mathcal O}}

\newcommand{\RR}{{\mathcal R}}

\newcommand{\TT}{{\mathcal T}}

\newcommand{\YY}{{\mathcal Y}}
\newcommand{\aaa}{{\bf a}}
\newcommand{\mmm}{{\bf m}}

\newcommand{\www}{\widetilde}

\newcommand{\oooo}{\overline}
\newcommand{\uuuu}{\underline}

\newcommand{\paa}{\partial}

\DeclareMathOperator{\Aut}{Aut}

\DeclareMathOperator{\id}{id}
\DeclareMathOperator{\Imm}{Im}
\DeclareMathOperator{\Lie}{Lie}
\DeclareMathOperator{\lcm}{lcm}
\DeclareMathOperator{\mult}{mult}

\DeclareMathOperator{\pr}{pr}
\DeclareMathOperator{\rank}{rank}
\DeclareMathOperator{\rk}{rk}
\DeclareMathOperator{\Rad}{Rad}
\DeclareMathOperator{\Ree}{Re}

\DeclareMathOperator{\Stab}{Stab}

\begin{document}

\theoremstyle{plain}
\newtheorem{lemma}{Lemma}[section]
\newtheorem{definition/lemma}[lemma]{Definition/Lemma}
\newtheorem{theorem}[lemma]{Theorem}
\newtheorem{proposition}[lemma]{Proposition}
\newtheorem{corollary}[lemma]{Corollary}
\newtheorem{conjecture}[lemma]{Conjecture}
\newtheorem{conjectures}[lemma]{Conjectures}

\theoremstyle{definition}
\newtheorem{definition}[lemma]{Definition}
\newtheorem{withouttitle}[lemma]{}
\newtheorem{remark}[lemma]{Remark}
\newtheorem{remarks}[lemma]{Remarks}
\newtheorem{example}[lemma]{Example}
\newtheorem{examples}[lemma]{Examples}

\title[Distinguished bases and Stokes regions]
{Distinguished bases and Stokes regions
for the simple and the simple elliptic singularities} 

\author{Claus Hertling \and C\'eline Roucairol}

\address{Claus Hertling\\
Lehrstuhl f\"ur Mathematik VI, Universit\"at Mannheim, Seminargeb\"aude
A 5, 6, 68131 Mannheim, Germany}

\email{hertling@math.uni-mannheim.de}

\address{C\'eline Roucairol\\ Avignon, France}

\email{roucairol@googlemail.com}


\thanks{This work was supported by the DFG grants 
He2287/2-2 and He2287/4-1 
}

\date{June 04, 2018}

%\maketitle

\begin{abstract}
{\small Isolated hypersurface singularities come equipped with a
Milnor lattice, a $\Z$-lattice of finite rank, and a set of
{\it distinguished} $\Z$-bases of this lattice. 
Usually these bases are constructed from {\it one} morsification 
and {\it all possible} choices of distinguished systems of paths. 
But what does one obtain if one considers {\it all possible} morsifications
and {\it one} fixed distinguished system of paths? Looijenga asked
this question 1974 for the simple singularities. He and Deligne
found that one obtains a bijection between Stokes regions in a 
universal unfolding and the set of distinguished bases modulo signs.
This allows to see the base space of the universal unfolding
as an atlas of Stokes data.
Here we reprove their result and extend it to the simple
elliptic singularities. We use more conceptual arguments,
moduli spaces of marked singularities (i.e. Teichm\"uller spaces
for singularities), extensions of them to F-manifolds,
and the actions of symmetries of singularities on
the Milnor lattices and these moduli spaces.
We use and extend results of Jaworski on the Lyashko-Looijenga
maps for the simple elliptic singularities.
The sections 2 and 3 give a survey on singularities and the
associated objects which allows to read the paper independently
of other sources.
}
\end{abstract}

\maketitle

\tableofcontents

\setcounter{section}{0}

\section{Introduction}\label{s1}
\setcounter{equation}{0}

\noindent
This paper has three parts.
The first part is an introduction to isolated hypersurface
singularities. It makes the paper readable independently
of other sources, and it gives a basis for the other
two parts. 

The second part is an extension to the simple elliptic 
singularities of work which Looijenga and Deligne did 1974
for the simple singularities. This is the central part of 
the paper.

The third part is an extension and refinement of work of
Jaworski 1986--1988 on the Lyashko-Looijenga maps for
the simple elliptic singularities. The arguments are 
less conceptual and more computational and more laborious
than the arguments in the second part.
We need it in order to determine the sizes of certain
finite sets which are in bijection by the second part.

An {\it isolated hypersurface singularity} is a holomorphic
function germ $f:(\C^{n+1},0)\to (\C,0)$ with an 
isolated singularity at 0.
In order to see its topology, one chooses a 
{\it good representative} $f:Y\to T$ with 
$Y\subset\C^{n+1}$ a suitable neighborhood of 0 and
$T\subset \C$ a small disk around 0.
The {\it Milnor lattice} $Ml(f)$ is the 
(reduced for $n=0$) middle homology 
$H_n^{(red)}(f^{-1}(\tau),\Z)\cong\Z^\mu$ 
of a regular fiber $f^{-1}(\tau)$ of $f:Y\to T$ for some
$\tau\in T\cap\R_{>0}$. 
Here $\mu\in\Z_{>0}$ is the {\it Milnor number} of $f$.
The Milnor lattice comes equipped
with a monodromy $M_h$, an intersection form
$I$, a {\it Seifert form} $L$, and a set $\BB(f)$ of
certain $\Z$-bases of $Ml(f)$, the {\it distinguished bases}.
$M_h$ is a quasiunipotent automorphism,
$I$ is a $(-1)^n$-symmetric bilinear form,
$L$ is a unipotent bilinear form, and $L$ determines
$M_h$ and $I$. The group 
$G_\Z(f):=\Aut(Ml(f),M_h,I,L)=\Aut(Ml(f),L)$ will be important.

The distinguished bases are constructed as follows. First,
one chooses a morsification $F^{(mor)}:Y^{(mor)}\to T$
of $f:Y\to T$, that is a deformation such that the one
singularity of $f:Y\to T$ with Milnor number $\mu$ 
splits into $\mu$ $A_1$-singularities $x^{(j)}$, 
$j=1,...,\mu,$
of $F^{(mor)}:Y^{(mor)}\to T$ with pairwise different
critical values $u_j=F^{(mor)}(x^{(j)})$, $j=1,...,\mu$,
with $|u_j|<\tau$. 
Second, one chooses a {\it distinguished system of paths}.
That is a system of $\mu$ paths $\gamma_j$, $j=1,...,\mu$,
from $u_{\sigma(j)}$ to $\tau$ for some permutation
$\sigma\in S_\mu$, which do not intersect except at $\tau$
and which arrive at $\tau$ in clockwise order.
Third, one shifts from the $A_1$-singularity above each
value $u_{\sigma(j)}$ the (up to the sign unique)
vanishing cycle along $\gamma_j$ to
$H_n^{(red)}((F^{(mor)})^{-1}(\tau),\Z)$ and then by a 
canonical isomorphism to $Ml(f)$ and calls the image
$\delta_j$. The tuple $\uuuu\delta =(\delta_1,...,\delta_\mu)$
turns out to be a $\Z$-basis of $Ml(f)$ and is called
a {\it distinguished basis}.
One morsification, all possible choicees of distinguished
systems of paths and both possible choices $\pm\delta_j$
of signs for each cycle give all distinguished bases.

The Stokes matrix of one distinguished basis 
$\uuuu\delta$ is the matrix
$S=(-1)^{(n+1)(n+2)/2}L(\uuuu\delta^t,\uuuu\delta)^t$.
It is an upper triangular integer matrix with 1's on
the diagonal. The following table gives some information
on the sizes of the sets $\BB(f)$ and 
$|\{\textup{Stokes matrices}\}|$.
\begin{eqnarray}\label{1.1}
\begin{array}{lll}
f & |\BB(f)| & |\{\textup{Stokes matrices}\}| \\ \hline 
\textup{simple singularity} & \textup{finite} & \textup{finite} \\
\textup{simple elliptic singularity} & \textup{infinite} & \textup{finite} \\
\textup{any other singularity} & \textup{infinite} & \textup{infinite}
\end{array}
\end{eqnarray}
The last line of it was proved only recently by Ebeling
\cite{Eb18}. The other two lines are explained for example
in \cite{Eb18} or in remark \ref{t7.2} (i) below.

The simple singularities $A_\mu\ (\mu\geq 1)$, 
$D_\mu\ (\mu\geq 4)$, $E_6$, $E_7$ and $E_8$ and the
simple elliptic singularities 
$\www E_6$, $\www E_7$ and $\www E_8$ are the first ones
in Arnold's lists \cite[ch. 15.1]{AGV85} of isolated
hypersurface singularities.
The simple singularities have no $\mu$-constant parameter.
The simple elliptic singularities are 1-parameter families.
See subsection \ref{s4.1} for normal forms for all of them.

Deligne \cite{De74} characterized $\BB(f)$ and calculated
the number $|\BB(f)|$ for the simple singularities.
Yu \cite{Yu90}\cite{Yu96}\cite{Yu99} derived from that
the number $|\{\textup{Stokes matrices}\}|$
for the simple singularities.
Kluitmann characterized $\BB(f)$ for the simple elliptic
singularities. He calculated the number 
$|\{\textup{Stokes matrices}\}|$
for $\www E_6$ in \cite{Kl83} and for $\www E_7$
in \cite{Kl87}, by huge combinatorial efforts.
The number $|\{\textup{Stokes matrices}\}|$ for $\www E_8$
was not calculated before this paper.
In corollary \ref{t7.3} we recover Kluitmann's numbers
for $\www E_6$ and $\www E_7$, and we give the number
for $\www E_8$, by a completely different method.
Our method combines a natural bijection in the second
part with the calculation of three numbers in the third
part, the degrees of certain Lyashko-Looijenga maps
for $\www E_6$, $\www E_7$ and $\www E_8$.

The simple singularities $f=f(x_0,...,x_n)$ have
{\it universal unfoldings} 
\begin{eqnarray}\label{1.2}
F^{alg}(x_0,...,x_n,t_1,...,t_\mu)=F(x,t)=F_t(x)
=f(x)+\sum_{j=1}^\mu t_jm_j,
\end{eqnarray}
with $m_1,...,m_\mu\in\C[x]$ the monomials in 
table \eqref{4.4} and with parameters
$t\in M^{alg}:=\C^\mu$.

1974 Looijenga \cite{Lo74} and Lyashko (but his work
was published only later in \cite{Ly79}\cite{Ly84})
considered the {\it Lyashko-Looijenga map}
\begin{eqnarray}\label{1.3}
LL^{alg} : M^{alg}\to  M_{LL}^{(\mu)} :=
\{y^\mu+\sum_{j=0}^{\mu-1}s_jy^j\, |
\, (s_1,...,s_\mu)\in\C^\mu\} \hspace*{1cm} \\
 t\mapsto  \prod_{j=1}^\mu (y-u_j)\ \textup{with }
(u_1,...,u_\mu)\textup{ the critical values of }F^{alg}_t.
\nonumber
\end{eqnarray}
for the simple singularities.
It is a branched covering of a finite degree $\deg LL^{alg}$,
see theorem \ref{t6.1} for details. 

Looijenga posed the following problem: 
Consider a generic polynomial $p(y)=\prod_{j=1}^\mu(y-u_j)
\in M_{LL}^{(\mu)}$. Then $F_t$ for any 
$t\in (LL^{alg})^{-1}(p(y))$ is a morsification
of $f$ with the same critical values $u_1,...,u_\mu$.
Now fix {\it one distinguished system of paths} 
from $u_1,...,u_\mu$ to $\tau$. Each morsification $F^{alg}_t$
with $t\in (LL^{alg})^{-1}(p(y))$ gives one distinguished basis
$\uuuu\delta=(\delta_1,...,\delta_\mu)$ up to signs.
One obtains a map
\begin{eqnarray}\label{1.4}
LD:(LL^{alg})^{-1}(p(y))\to \BB(f)/G_{sign,\mu},
\end{eqnarray}
where the group $G_{sign,\mu}:=\{\pm 1\}^\mu$ acts on
$\uuuu\delta=(\delta_1,...,\delta_\mu)$ by sign changes.
An easy argument in \cite{Lo74} which uses that $LL^{alg}$
is a branched covering, shows that the map $LD$ is surjective.
Looijenga asked whether $LD$ is injective. He proved
this for the $A_\mu$-singularities.
Then Deligne \cite{De74} calculated $|\BB(f)|$ for all simple
singularities and showed $|\BB(f)/G_{sign,\mu}|=\deg LL^{alg}
=|(LL^{alg})^{-1}(p(y))|$. This proved that $LD$ is a bijection
for all simple singularities.
Deligne's letter \cite{De74} to Looijenga is not published.
The result that $LD$ is a bijection 
is stated in \cite{Mi89}, \cite{Yu90} and below in 
theorem \ref{t7.1}.

A central part of this paper is an extension of this result
to the simple elliptic singularities. Here a universal 
unfolding of a single simple elliptic singularity is not
good enough. In subsection \ref{s4.2} we present a global
family of functions
\begin{eqnarray}\label{1.5}
F^{alg}(x_0,...,x_n,t_1,...,t_{\mu-1},\lambda)&=&
F^{alg}(x,t',\lambda)=F^{alg}_{t',\lambda}(x)\\
&=& f_\lambda(x) + \sum_{j=1}^{\mu-1} t_jm_j,
\nonumber
\end{eqnarray}
with $m_1,...,m_{\mu-1}\in \C[x]$ the 
monomials in table \eqref{4.7},
$f_\lambda(x)$ the 1-parameter families in Legendre normal form
in \eqref{4.2} of the simple elliptic singularities and
with parameters $(t',\lambda)\in M^{alg}:=\C^{\mu-1}\times
(\C-\{0,1\})$. For each $\lambda\in\C-\{0,1\}$, the 
family $F^{alg}$ is (locally) a universal unfolding 
of $f_\lambda$.
The family $F^{alg}$ is not completely canonical.
But the family $F^{mar}$ with parameter space 
$M^{mar}:=\C^{\mu-1}\times \H$, where 
$\H\to\C-\{0,1\}$ is the universal covering, is canonical.
This is made precise in theorem \ref{t4.3} in a way
which uses marked singularities, the fact that the
parameter space of each of the three families 
of marked simple elliptic 
singularities is $M^{mar}_\mu\cong\H$ \cite{GH17-1},
and a thickening of this space to $M^{mar}$. 
We obtain Lyashko-Looijenga maps
$LL^{alg}:M^{alg}\to M_{LL}^{(\mu)}$ and
$LL^{mar}:M^{mar}\to M_{LL}^{(\mu)}$. 

Analogously to $LD$ for the simple singularities, we obtain
a {\it Looijenga-Deligne map}
\begin{eqnarray}\label{1.6}
LD: (LL^{mar})^{-1}(p(y))\to \BB(f)/G_{sign,\mu}
\end{eqnarray}
for generic $p(y)\in M_{LL}^{(\mu)}$. 
A main result of this paper is that this map is a bijection,
see theorem \ref{t7.1}.
But our arguments are more involved than the arguments
for the simple singularities.
The surjectivity follows by the same easy argument in 
\cite{Lo74}, as soon as one has that the map
\begin{eqnarray}\label{1.7}
LL^{alg}:M^{alg}-(\textup{caustic}\cup\textup{Maxwell stratum})\\
\to M_{LL}^{(\mu)}-(\textup{discriminant }\DD_{LL}^{(\mu)})
\nonumber
\end{eqnarray}
is a finite covering. This is a hard theorem of Jaworski
\cite[Theorem 2]{Ja86} \cite[Proposition 1]{Ja88},
see theorem \ref{t6.2}.

As both sides of \eqref{1.6} are infinite, the injectivity
of $LD$ in \eqref{1.6} does not follow by a comparison
of numbers. We need an action of $G_\Z$ on $M^{mar}$
and on the middle homology bundle above $M^{mar}-\DD^{mar}$.
We need that $M^{mar}$ is an F-manifold with Euler field.
And we need that and how a Stokes matrix $S$ of
$LD(t)\in \BB(f)/G_{sign,\mu}$ for $t\in (LL^{mar})^{-1}(p(y))$
encodes the covering in \eqref{1.7}.
See the proof of theorem \ref{t7.1} for the details.

The  bijection in \eqref{1.6} induces also a bijection
\begin{eqnarray}\label{1.8}
(LL^{mar})^{-1}(p(y))/G_\Z 
\to \{\textup{Stokes matrices}\}/G_{sign,\mu}.
\end{eqnarray}
Both sides are finite, and the number
$|(LL^{mar})^{-1}(p(y))/G_\Z|$ is in a simple way related
to $\deg LL^{alg}$. Though Jaworski's proofs in 
\cite{Ja86} and \cite{Ja88} that $LL^{alg}$ in \eqref{1.7}
is a covering, do not allow to calculate 
its degree $\deg LL^{alg}$.

The main task in the third part of the paper is to extend
Jaworski's work and calculate $\deg LL^{alg}$.
In theorem \ref{t6.3} we obtain an extension of 
$M^{alg}$ above $\C-\{0,1\}$ to an orbibundle
$M^{orb}\stackrel{\pi_{orb}}{\to}\P^1$ 
such that $LL^{alg}$ extends to a  holomorphic map 
$LL^{orb}:M^{orb}\to M_{LL}^{(\mu)}$
which is outside of the $\mu$-constant stratum 
(and its translates by the unit field)
a branched covering and which maps 
$(\textup{caustic})\cup(\textup{Maxwell stratum})\cup
\pi_{orb}^{-1}(\{0,1,\infty\})$ to the discriminant
$\DD_{LL}^{(\mu)}\subset M_{LL}^{(\mu)}$. 
Detailed information about $M^{orb}$ and $LL^{orb}$
allows us to calculate the degree 
$\deg LL^{orb} \, (=\deg LL^{alg})$.

The first part of the paper consists of section \ref{s2},
the subsections \ref{s3.1} and \ref{s3.2} and the
first three pages of section \ref{s4}.
Section \ref{s2} recalls classical data and facts around
isolated hypersurface singularities, namely
Milnor fibrations, Milnor lattices, 
universal unfoldings, the base spaces as F-manifolds
with Euler fields, Lyashko-Looijenga maps, distinguished bases,
Stokes matrices, and Thom-Sebastiani type results.
Subsection \ref{s3.1} reviews results in 
\cite[Theorem 13.11 and Theorem 13.13]{He02} on symmetries
of singularities. Subsection \ref{s3.2} reviews
results in \cite[Theorem 4.3]{He11} on the moduli spaces
$M^{mar}_\mu$ of marked singularities.
The first three pages of section \ref{s4} give normal forms
for the simple and the simple elliptic singularities
and for the unfoldings. 

The second part of the paper consists of the subsections
\ref{s3.3} and \ref{s3.4}, the latter part of section \ref{s4}
and the sections \ref{s6} and \ref{s7}. 
Subsection \ref{s3.3} describes a thickening 
$M^{mar}\supset M^{mar}_\mu$ of the moduli spaces of marked
singularities. Theorem \ref{t4.3} in the latter part of
section \ref{s4} proves the claims about this thickening in 
the cases of the simple and the simple elliptic singularities.
Subsection \ref{s3.4} defines and discusses 
Looijenga-Deligne maps $LD$ in a general setting.
Section \ref{s7} states and proves the main result
theorem \ref{t7.1} that $LD$ is a bijection for each
simple singularity and each family of simple elliptic 
singularities. Corollary \ref{t7.3} provides the finite
numbers $|\{\textup{Stokes matrices}\}|$.
Section \ref{s6} states the old and new results on the 
Lyashko-Looijenga maps for the simple singularities
(theorem \ref{t6.1}, Lyashko and Looijenga)
and the simple elliptic singularities
(theorem \ref{t6.2}, Jaworski, and theorem \ref{t6.3}, new).
The most beautiful formula in section \ref{s6} is formula
\eqref{6.7} for $\deg LL^{alg}$ for the simple elliptic
singularities.

The third part of the paper consists of the sections 
\ref{s5}, \ref{s8}, \ref{s9} and \ref{s10}.
Section \ref{s5} makes the general discussion of the
symmetries of singularities in subsection \ref{s3.1} 
explicit in the cases of the simple and the simple
elliptic singularities.
Section \ref{s8} follows Fulton's book \cite{Fu84}
and extends some results there to the case of 
{\it smooth cone bundles} (definition \ref{t8.1}).
We need this for the proof of corollary \ref{t8.6}
which is used in the proof of formula \eqref{6.7}
in section \ref{s10}.
The long section \ref{s9} provides for the simple
elliptic singularities the extension of $M^{alg}
=\C^{\mu-1}\times (\C-\{0,1\})$ to $\lambda=0$
such that $LL^{alg}$ extends well. Finding the right way
to glue into $M^{alg}$ a fiber above $\lambda=0$
was the most laborious part of this paper.
Section \ref{s10} combines this with the symmetries in
section \ref{s5} and provides the right extensions of 
$M^{alg}$ to $\lambda=1$ and $\lambda=\infty$,
and it uses corollary \ref{t8.6} to prove the
formula \eqref{6.7} for $\deg LL^{alg}$ for the simple
elliptic singularities.

The first author thanks the organizers of the conference
"Moduli spaces and applications in geometry, topology, 
analysis and mathematical physics" in Beijing February 
27 -- March 3, 2017, for the invitation to the conference, 
and both authors thank especially Lizhen Ji for a lot
of patience during the preparation of this paper.

\section{Topology and unfoldings of isolated hypersurface singularities}\label{s2}
\setcounter{equation}{0}

\noindent
An {\it isolated hypersurface singularity} (short: {\it singularity})
is a holomorphic function germ $f:(\C^{n+1},0)\to (\C,0)$ with an isolated 
singularity at $0$. 
Such objects were studied intensively since the end of the 1960'ies.
In this section, we review classical facts on their topology and their
unfoldings and fix some notations.
For the topology compare \cite{AGV88}
and \cite{Eb07}. For the unfoldings compare \cite{AGV85}
and (especially for F-manifolds) \cite{He02}.

The {\it Jacobi ideal} of $f$ is the ideal 
$J_f:=(\frac{\paa f}{\paa x_i}) \subset\OO_{\C^{n+1},0}$, 
its {\it Jacobi algebra} is the quotient
$\OO_{\C^{n+1},0}/J_f$, its {\it Milnor number} is the finite number
$\mu:=\dim\OO_{\C^{n+1},0}/J_f.$ 

\subsection{Topology of singularities}\label{s2.1}
A {\it good representative} of $f$ has to be defined with some 
care \cite{Mi68}\cite{AGV88}\cite{Eb07}. It is $f:Y\to T$
with $Y\subset\C^{n+1}$ a suitable small neighborhood of 0 and 
$T\subset\C$ a small disk around 0.
Then $f:Y'\to T'$ with $Y'=Y-f^{-1}(0)$ and 
$T'=T-\{0\}$ is a locally trivial $C^\infty$-fibration,
the  {\it Milnor fibration}. Each fiber has the
homotopy type of a bouquet of $\mu$ $n$-spheres \cite{Mi68}.

Therefore the (reduced for $n=0$) middle homology groups are {}\\{}
$H_n^{(red)}(f^{-1}(\tau),\Z) \cong \Z^\mu$ for $\tau\in T'$.
Each comes equipped with an {\it intersection form} $I$, 
which is a datum of one fiber,
a {\it monodromy} $M_h$ and a {\it Seifert form} $L$, which come from the 
Milnor fibration,
see \cite[I.2.3]{AGV88} for their definitions 
(for the Seifert form, there are several
conventions in the literature, we follow \cite{AGV88}). 
$M_h$ is a quasiunipotent automorphism, $I$ and $L$ are 
bilinear forms with values in $\Z$,
$I$ is $(-1)^n$-symmetric, and $L$ is unimodular. $
L$ determines $M_h$ and $I$ because of the formulas
\cite[I.2.3]{AGV88}
\begin{eqnarray}\label{2.1}
L(M_ha,b)&=&(-1)^{n+1}L(b,a),\\ \label{2.2}
I(a,b)&=&-L(a,b)+(-1)^{n+1}L(b,a).
\end{eqnarray}
The lattices $H_n(f^{-1}(\tau),\Z)$ for all Milnor fibrations
$f:Y'\to T'$ and then all $\tau\in\R_{>0}\cap T'$ are canonically isomorphic,
and the isomorphisms respect $M_h$, $I$ and $L$. 
This follows from Lemma 2.2 in \cite{LR73}. 
These lattices are identified and called {\it Milnor lattice} $Ml(f)$.
The group $G_\Z$ is 
\begin{eqnarray}\label{2.3}
G_\Z=G_\Z(f):= \Aut(Ml(f),L)=\Aut(Ml(f),M_h,I,L),
\end{eqnarray}
the second equality is true because $L$ determines $M_h$ and $I$.
We will use the notation $Ml(f)_\C:=Ml(f)\otimes_\Z \C$,
and analogously for other rings $R$ with $\Z\subset R\subset \C$,
and the notations
\begin{eqnarray*}
Ml(f)_\lambda&:=&\ker((M_h-\lambda\id)^\mu:Ml(f)_\C\to Ml(f)_\C)
\subset Ml(f)_\C,\\
Ml(f)_{1,\Z}&:=& Ml(f)_1\cap Ml(f)\subset Ml(f),\\
Ml(f)_{\neq 1}&:=&\bigoplus_{\lambda\neq 1}Ml(f)_\lambda \subset 
Ml(f)_\C,\\
Ml(f)_{\neq 1,\Z}&:=& Ml(f)_{\neq 1}\cap Ml(f)\subset Ml(f).
\end{eqnarray*}

The formulas \eqref{2.1} and \eqref{2.2} show
$I(a,b)= L((M_h-\id)a,b)$. Therefore the eigenspace with eigenvalue
1 of $M_h$ is the radical $\Rad(I)\subset Ml(f)$ of $I$. 
By \eqref{2.2} $L$ is 
$(-1)^{n+1}$-symmetric on the radical of $I$.

\subsection{Unfoldings of singularities}\label{s2.2}

The notion of an unfolding of an isolated hypersurface singularity $f$ goes
back to Thom and Mather. An {\it unfolding} of $f$ is a holomorphic
function germ $F:(\C^{n+1}\times M,0)\to (\C,0)$ such that
$F|_{(\C^{n+1},0)}=f$ and such that $(M,0)$ 
is the germ of a complex manifold.
Its Jacobi ideal is 
$J_F:=(\frac{\paa F}{\paa x_i})\subset \OO_{\C^{n+1}\times M,0}$,
its critical space is the germ $(C,0)\subset (\C^{n+1}\times M,0)$
of the zero set of $J_F$ with the canonical complex structure.
The projection $(C,0)\to (M,0)$ is finite and flat of degree $\mu$.
A kind of Kodaira-Spencer map is the $\OO_{M,0}$-linear map
\begin{eqnarray}\label{2.4}
\aaa_C:\TT_{M,0}\to \OO_{C,0},\quad X\mapsto\www X(F)|_{(C,0)}
\end{eqnarray}
where $\www X$ is an arbitrary lift of $X\in\TT_{M,0}$ to $(\C^{n+1}\times M,0)$.

We will use the following notion of morphism between unfoldings. 
Let $F_i:(\C^{n+1}\times M_i,0)\to(\C,0)$ for $i\in\{1,2\}$ be
two unfoldings of $f$ with projections $\pr_i:(\C^{n+1}\times M_i,0)\to(M_i,0)$.
A {\it morphism} from $F_1$ to $F_2$ is a pair 
$(\Phi,\varphi)$ of map germs such that the following diagram
commutes,
\begin{eqnarray}\label{2.5}
\begin{xy}
\xymatrix{ (\C^{n+1}\times M_1,0) \ar[r]^\Phi \ar[d]^{\pr_1} 
& (\C^{n+1}\times M_2,0) \ar[d]^{\pr_2}\\
(M_1,0) \ar[r]^{\varphi} & (M_2,0)  }
\end{xy}
\end{eqnarray}
and
\begin{eqnarray}\label{2.6}
\Phi|_{(\C^{n+1}\times\{0\},0)}&=&\id,\\
F_1 &=& F_2\circ\Phi \label{2.7}
\end{eqnarray}
hold. Then one says that $F_1$ {\it is induced} by $(\Phi,\varphi)$ from $F_2$.
An unfolding is {\it versal} if any unfolding is induced from it by a 
suitable morphism. A versal unfolding $F:(\C^{n+1}\times M,0)\to(\C,0)$ is
{\it universal} if the dimension of the parameter space $(M,0)$ is
minimal. Universal unfoldings exist by work of Thom and Mather.
More precisely, an unfolding is versal if and only if the map 
$\aaa_C$ is surjective, and it is universal if and only if the map $\aaa_C$
is an isomorphism (see e.g. \cite[ch. 8]{AGV85} for a proof). 
Observe that $\aaa_C$ is surjective/an isomorphism
if and only if its restriction to 0, the map
\begin{eqnarray}\label{2.8}
\aaa_{C,0}:T_0M\to \OO_{\C^{n+1},0}/J_f
\end{eqnarray}
is surjective/an isomorphism. Therefore an unfolding
\begin{eqnarray}\label{2.9}
F(x_0,...,x_n,t_1,...,t_\mu)=F(x,t)=F_t(x)=f(x)+\sum_{j=1}^\mu m_it_i,
\end{eqnarray}
with $(M,0)=(\C^\mu,0)$ with coordinates $t=(t_1,...,t_\mu)$ 
where $m_1,...,m_\mu\in\OO_{\C^{n+1},0}$ represent a basis of the 
Jacobi algebra $\OO_{\C^{n+1},0}/J_f$, is universal.

\subsection{F-manifolds}\label{s2.3}

The base space of a universal unfolding 
is an {\it F-manifold} with {\it Euler field}. 

\begin{definition}\label{t2.1}\cite{HM99}\cite{He02}
(a) An {\it F-manifold} is a complex manifold $M$ together
with a holomorphic commutative and associative multiplication
$\circ$ on its holomorphic tangent bundle $TM$ and with a unit field
$e\in \TT_M$ such that the integrability condition 
\begin{eqnarray}\label{2.10}
\Lie_{X\circ Y}(\circ)= X\circ \Lie_Y(\circ)+Y\circ \Lie_X(\circ)
\end{eqnarray}
holds.

(b) Let $(M,\circ,e)$ be an F-manifold. An Euler field (of weight 1)
is a global holomorphic vector field $E\in\TT_M$ with $\Lie_E(\circ)=\circ$.
\end{definition}

F-manifolds are studied in \cite[ch. 2--5]{He02}. In the case of a 
universal unfolding $F:(\C^{n+1}\times M,0)\to(\C,0)$,
its base space inherits from the isomorphism $\aaa_C^{-1}:\OO_{C,0}\to\TT_M$
a multiplication. It satisfies \eqref{2.10}, and $e:=\aaa_C^{-1}([1])$
and $E:=\aaa_C^{-1}([F])$ are the unit field and an Euler field
\cite[Theorem 5.3]{He02}, so it is an F-manifold with Euler field.

We call a universal unfolding {\it universal} and not just {\it semiuniversal}, because
the morphism $\varphi$ in \eqref{2.5} between the base spaces of 
any two universal unfoldings is unique \cite[Theorem 5.4]{He02}
(but $\Phi$ on the total spaces is not unique). 
Therefore the base space of a universal unfolding is (with its structure 
as F-manifold with Euler field) unique up to unique isomorphism.

The following result on decompositions of germs of F-manifolds 
is a very instructive application of the integrability condition \eqref{2.10}, 
and it is especially telling in the case of isolated hypersurface singularities.

\begin{theorem}\label{t2.2}\cite[Theorem 2.11]{He02}
Let $(M,p)$ be the germ in $p\in M$ of an F-manifold.
It is an elementary fact from commutative algebra that the algebra
$T_pM$ decomposes into a direct sum $\bigoplus_{k=1}^l (T_pM)_k$
of irreducible local algebras (it is just the decomposition into
simultaneous generalized eigenspaces with respect to all (commuting!)
multiplication endomorphisms). 

This decomposition extends uniquely to a decomposition $(M,p)=\prod_{k=1}^l (M_k,p_k)$
of germs of F-manifolds. These germs are irreducible germs of F-manifolds.
If $(M,p)$ has an Euler field, the Euler field decomposes into a sum of Euler fields.
\end{theorem}

In the case of a good representative $F:\YY\to T$ of a universal unfolding $F$,
for any $t\in M$, the canonical decomposition from theorem \ref{t2.2}
of $(M,t)$ into a product of germs of F-manifolds is a canonical decomposition
into a product of germs of base spaces of universal unfoldings of the germs
of $F_t$ at all its critical points.

At generic $t$, this is a decomposition into 1-dimensional F-manifolds, and
the eigenvalues $u_1,...,u_\mu$ of the Euler field form there local coordinates,
Dubrovin's {\it canonical coordinates}. The Euler field has there the shape
$E=\sum_{j=1}^\mu u_j e_j$ where $e_j=\frac{\paa}{\paa u_j}$, 
the multiplication is given by $e_i\circ e_j=\delta_{ij}e_i$, the
unit field is $e=\sum_{j=1}^\mu e_j$, and the values $u_1,...,u_\mu$
are the critical values of $F_t$, i.e. the values of $F_t$ at its critical
points. Up to isomorphism there is only one germ of a 1-dimensional
F-manifold, which is called $A_1$. Then the germ $(M,t)$ at generic $t$
is as a germ of an F-manifold of the type $A_1^\mu$.

\subsection{Lyashko-Looijenga map}\label{s2.4}

Looijenga \cite{Lo74} was close to the notion of an F-manifold. 
He had already the canonical coordinates at generic points. 
And he and Lyashko \cite{Ly79}\cite{Ly84} studied the Lyashko-Looijenga map 
and its behaviour near the
caustic and the Maxwell stratum. For $\mu\in\Z_{\geq 1}$ define 
\begin{eqnarray}\label{2.11}
M_{LL}^{(\mu)}&=&\{y^\mu+\sum_{j=1}^\mu s_jy^{j-1}\, |\, 
(s_1,...,s_\mu)\in\C^\mu\} \cong\C^\mu,\\
\DD_{LL}^{(\mu)}&:=& \{p(y)\in 
M_{LL}^{(\mu)}\, | \, p(y)\textup{ has multiple roots}\}.
\label{2.12}
\end{eqnarray}
$\DD_{LL}^{(\mu)}$ is a hypersurface in $M_{LL}^{(\mu)}$.
Let $F:(\C^{n+1}\times M,0)\to(\C,0)$ be a universal unfolding of
a singularity $f$. 
Let $F:\YY\to T$ be a good representative of it with base space $M$. 
Then its {\it Lyashko-Looijenga map} is the holomorphic map
\begin{eqnarray}\label{2.13}
LL:M&\to& M_{LL}^{(\mu)},\quad t\mapsto \prod_{j=1}^\mu (y-u_j),\quad\textup{where }
u_1,...,u_\mu\\
&&\textup{are the critical values of }
F_t\textup{ (with multiplicities).}\nonumber
\end{eqnarray}
Define the {\it caustic} $\KK_3\subset M$ 
and the {\it Maxwell stratum} $\KK_2\subset M$ by 
\begin{eqnarray}\label{2.14}
\KK_3&:=& \{t\in M\, |\, F_t
\textup{ has less than }\mu \textup{ singularities}\},\\
\KK_2&:=& \textup{the closure in }M\textup{ of the set }\{t\in M\, |\, 
F_t\textup{ has }\mu \label{2.15}\\
&&\textup{ singularities, but less than }\mu \textup{ critical values}\} . \nonumber
\end{eqnarray}
They are hypersurfaces in $M$. 
%In fact, they are irreducible hypersurfaces
%(see \cite[Lemma 2.3]{Gr78} for $\KK_3$ and \cite{AH18} for $\KK_2$).

The Lyashko-Looijenga map $LL$ restricts to a locally biholomorphic map
$LL:M-(\KK_3\cup \KK_2)\to M_{LL}^{(\mu)}-\DD_{LL}^{(\mu)}$
(because the $u_1,...,u_\mu$ are
local coordinates on $M-\KK_3$), 
it maps $\KK_3\cup \KK_2$ to $\DD_{LL}^{(\mu)}$,
and it is a branched covering of order 3 respectively 2 at generic points of
$\KK_3$ respectively $\KK_2$. All of this was proved by Lyashko \cite{Ly79}\cite{Ly84} 
and Looijenga \cite{Lo74}. Now it is an easy consequence of the F-manifold structure.
At a generic point $t$ of $\KK_3$, the germ of the F-manifold is of the type
$A_2A_1^{\mu-2}$. Here $A_2$ is the first in the countable series of irreducible
germs of massive F-manifolds \cite[ch. 4.2]{He02}.

\subsection{Distinguished bases and Stokes matrices of singularities}\label{s2.5}

Good references for distinguished bases are \cite{AGV88} and \cite{Eb07}.
We sketch their construction and properties.

Choose a universal unfolding of $f$,
a {\it good representative} $F:\YY\to T$ of it with base space $M$,
and a generic parameter $t\in M$. Then $F_t:Y_t\to T$
with $T\subset \C$ the same disk as that for a Milnor fibration $f:Y\to T$ 
above and $Y_t\subset\C^{n+1}$
is a {\it morsification} of $f$. 
It has $\mu$ $A_1$-singularities,
and their critical values $u_1,...,u_\mu\in T$ 
are pairwise different. Their numbering is also a choice. 
Now choose a value $\tau\in T\cap\R_{>0}-\{u_1,...,u_\mu\}$ and
a {\it distinguished system of paths}. That is
a system of $\mu$ paths $\gamma_j$, $j=1,...,\mu$, from
$u_j$ to $\tau$ which do not intersect except at $\tau$
and which arrive at $\tau$ in clockwise order.
Finally, shift from the $A_1$ singularity above each value $u_j$
the (up to sign unique) vanishing cycle along $\gamma_j$
to the Milnor lattice $Ml(f)=H_n(f^{-1}(\tau),\Z)$,
and call the image $\delta_j$. 

The tuple $\uuuu{\delta}=(\delta_1,...,\delta_\mu)$ 
is a $\Z$-basis of 
$Ml(f)$. All such bases are called {\it distinguished bases}.
They form one orbit $\BB(f)$ of an action of a semidirect product 
$G_{sign,\mu}\rtimes\textup{Br}_\mu$. 
Here $\textup{Br}_\mu$ is the braid group with $\mu$ strings, 
see \cite{AGV88} or \cite{Eb07} for its action. 
The {\it sign change group} $G_{sign,\mu}:=\{\pm 1\}^\mu$ acts simply
by changing the signs of the entries of the tuples
$(\delta_1,...,\delta_\mu)$. 
The members of the distinguished bases are called 
{\it vanishing cycles}.

The matrix $L(\uuuu{\delta}^t,\uuuu{\delta})
=(L(\delta_i,\delta_j))_{i,j=1,...,\mu}$ 
of the Seifert form with respect to a distinguished basis 
$\uuuu{\delta}=(\delta_1,...,\delta_\mu)$ 
is a lower triangular matrix with $(-1)^{(n+1)(n+2)/2}$ on the diagonal.
The {\it Stokes matrix} of the distinguished basis $\uuuu{\delta}$
is by definition the upper triangular matrix in $M(\mu\times\mu,\Z)$
\begin{eqnarray}\label{2.16}
S:=(-1)^{(n+1)(n+2)/2}\cdot L(\uuuu{\delta}^t,\uuuu{\delta})^t
\end{eqnarray}
with 1's on the diagonal. Then \eqref{2.1} and \eqref{2.2} give 
\begin{eqnarray}\label{2.17}
M_h(\uuuu{\delta})&=&\uuuu{\delta}\cdot (-1)^{n+1}\cdot S^{-1}S^t,\\
I(\uuuu{\delta}^t,\uuuu{\delta})&=&(-1)^{n(n+1)/2}\cdot 
(S+(-1)^nS^t).\label{2.18}
\end{eqnarray}

The {\it Coxeter-Dynkin diagram} of the distinguished basis $\uuuu{\delta}$ 
encodes $S$ in a geometric way. It has $\mu$ vertices which are numbered
from 1 to $\mu$. Between two vertices $i$ and $j$ with $i<j$
one draws

\begin{tabular}{ll}
no edge & if $S_{ij}=0$, \\
$|S_{ij}|$ edges & if $S_{ij}<0$, \\
$S_{ij}$ dotted edges & if $S_{ij}>0$. 
\end{tabular}

Coxeter-Dynkin diagrams of many singularities were calculcated
by A'Campo, Ebeling, Gabrielov and Gusein-Zade. Some of them
can be found in \cite{Ga74}, \cite{Eb83} and \cite{Eb07}.
Each Coxeter-Dynkin diagram of any singularity is connected.
We will use this important result in lemma \ref{t2.3}.
There are three proofs of it, by Gabrielov \cite{Ga74}, 
Lazzeri \cite{La73} and L\^e \cite{Le73}.

The Picard-Lefschetz transformation on $Ml(f)$ 
of a vanishing cycle $\delta$ is
\begin{eqnarray}\label{2.19}
s_{\delta}(b)&:=&b-(-1)^{n(n+1)/2}\cdot I(\delta,b)\cdot \delta.
\end{eqnarray}
For $n$ even $I(\delta,\delta)=(-1)^{n(n+1)/2}\cdot 2$ and $s_\delta$
is the identity on the subspace in $Ml(f)$ orthogonal to $\delta$
and $-\id$ on $\Z\cdot\delta$. For $n$ odd $s_\delta$ is unipotent
with kernel of $s_\delta-\id$ of rank $\mu-1$.
In both cases $s_\delta$ determines $\delta$ up to the sign.

The monodromy $M_h$ is
\begin{eqnarray}\label{2.20}
M_h &=& s_{\delta_1}\circ ...\circ s_{\delta_\mu}
\end{eqnarray}
for any distinguished basis $\uuuu{\delta}=(\delta_1,...,\delta_\mu)$.

Let us formulate a correspondence for later use.

\begin{lemma}\label{t2.3}
The orbit under $G_{sign,\mu}\rtimes \textup{Br}_\mu$ of a tuple
\begin{eqnarray}\label{2.21}
((u_1,...,u_\mu),\textup{a distinguished system of paths}, 
\textup{a Stokes matrix }S)
\end{eqnarray}
where $u_1,...,u_\mu\in\C$ are pairwise different
is equivalent to the isomorphism class of a $\Z$-lattice bundle
of rank $\mu$ above $\C-\{u_1,...,u_\mu\}$.
The only automorphisms (which fix the basis $\C-\{u_1,...,u_\mu\}$) 
of this $\Z$-lattice bundle are $\pm\id$.
\end{lemma}

{\bf Proof:}
If a morsification $F_t$ with critical values $u_1,...,u_\mu$ 
and distinguished basis above the distinguished system of paths 
with Stokes matrix $S$ exists,
then the $\Z$-lattice bundle is up to isomorphism
the flat extension to $\C-\{u_1,...,u_\mu\}$ of the middle homology bundle
\begin{eqnarray}\label{2.22}
\bigcup_{\tau\in T-\{u_1,...,u_\mu\}}H_n(F_t^{-1}(\tau),\Z).
\end{eqnarray}
If not, the $\Z$-lattice bundle is obtained from a case in \eqref{2.22}
by a suitable deformation.

The vanishing cycle near $u_j$ is up to its
sign uniquely determined by its Picard-Lefschetz transformation.

Any automorphism of the $\Z$-lattice bundle maps each of these 
vanishing cycles to $\pm$
itself. As the Coxeter-Dynkin diagram is connected, the only
automorphisms of the $\Z$-lattice bundle are $\pm\id$.
\hfill $\Box$

\subsection{Thom-Sebastiani type results}\label{s2.6}

A result of Thom and Sebastiani
compares the Milnor lattices and monodromies of 
the singularities $f=f(x_0,...,x_n),g=g(y_0,...,y_m)$ and
$f+g=f(x_0,...,x_n)+g(x_{n+1},...,x_{n+m+1})$.
There are extensions by Deligne for the Seifert form and 
by Gabrielov for distinguished bases. 
All results can be found in \cite[I.2.7]{AGV88}.
They are restated here.
There is a canonical isomorphism
\begin{eqnarray}\label{2.23}
\Phi:Ml(f+g)&\stackrel{\cong}{\longrightarrow} &Ml(f)\otimes Ml(g),\\
\textup{with } M_h(f+g)&\cong & M_h(f)\otimes M_h(g) \label{2.24}\\
\textup{and } 
L(f+g)&\cong& (-1)^{(n+1)(m+1)}\cdot L(f)\otimes L(g).\label{2.25}
\end{eqnarray}
If $\uuuu{\delta}=(\delta_1,...,\delta_{\mu(f)})$
and $\uuuu{\gamma}=(\gamma_1,...,\gamma_{\mu(g)})$ are
distinguished bases of $f$ and $g$  with Stokes matrices
$S(f)$ and $S(g)$, then 
$$\Phi^{-1}(\delta_1\otimes \gamma_1,...,
\delta_1\otimes \gamma_{\mu(g)},
\delta_2\otimes \gamma_1,...,
\delta_2\otimes \gamma_{\mu(g)},
...,
\delta_{\mu(f)}\otimes \gamma_1,...,
\delta_{\mu(f)}\otimes \gamma_{\mu(g)})$$
is a distinguished basis of $Ml(f+g)$,
that means, one takes the vanishing cycles 
$\Phi^{-1}(\delta_i\otimes \gamma_j)$ in the lexicographic order.
Then by \eqref{2.16} and \eqref{2.25}, the matrix 
\begin{eqnarray}\label{2.26}
S(f+g)=S(f)\otimes S(g)
\end{eqnarray}
(where the tensor product is defined
so that it fits to the lexicographic order) 
is the Stokes matrix of this distinguished basis.

In the special case $g=x_{n+1}^2$,
the function germ $f+g=f(x_0,...,x_n)+x_{n+1}^2\in \OO_{\C^{n+2},0}$
is called {\it stabilization} or {\it suspension} of $f$. 
As there are only two isomorphisms $Ml(x_{n+1}^2)\to\Z$, 
and they differ by a sign, there are two equally canonical
isomorphisms $Ml(f)\to Ml(f+x_{n+1}^2)$,
and they differ just by a sign. 
Therefore automorphisms and bilinear forms on $Ml(f)$ 
can be identified with automorphisms and bilinear forms on 
$Ml(f+x_{n+1}^2)$. In this sense
\begin{eqnarray}\label{2.27}
L(f+x_{n+1}^2) = (-1)^n\cdot L(f)\quad\textup{and}\quad 
M_h(f+x_{n+1}^2)= - M_h(f)
\end{eqnarray}
\cite[I.2.7]{AGV88}, and $G_\Z(f+x_{n+1}^2)=G_\Z(f)$.
The Stokes matrix $S$ does not change under stabilization.

\section{Marked singularities and their symmetries}\label{s3}
\setcounter{equation}{0}

\subsection{Symmetries of singularities}\label{s3.1}
Here we will review results from \cite[13.1 and 13.2]{He02}
on symmetries of singularities. A review with slightly
simplified proofs was already given in \cite[ch. 6]{He11}.
Let $f:(\C^{n+1},0)\to(\C,0)$ be a singularity, and 
let $F:\YY\to T$ be a good representative with base space
$M\subset\C^\mu$ (with coordinates $t=(t_1,...,t_\mu)$)
of a universal unfolding $(\C^{n+1}\times M,0)\to(\C,0)$.
Let 
\begin{eqnarray*}
\RR:=\{\varphi:(\C^{n+1},0)\to(\C^{n+1},0)\, |\, \varphi
\textup{ biholomorphic}\}
\end{eqnarray*}
be the group of all germs of coordinate changes, and let
\begin{eqnarray}\label{3.1}
\RR^f:=\{\varphi\in\RR\, |\, f\circ\varphi=f\}
\end{eqnarray}
be the group of symmetries of $f$. It is possibly
$\infty$-dimensional, but the group $j_k\RR^f$
of $k$-jets in $\RR^f$ is an algebraic group for any $k\in\Z_{\geq 1}$.
Let 
\begin{eqnarray}\label{3.2}
R_f:= j_1\RR^f/(j_1\RR^f)^0
\end{eqnarray}
be the finite group of components of $j_1\RR^f$. It is easy to see
that $R_f=j_k\RR^f/(j_k\RR^f)^0$ for any $k\in\Z_{\geq 1}$
\cite[Lemma 13.10]{He02}. Recall the definition of $G_\Z(f)$ in \eqref{2.3}.
There is a natural homomorphism
\begin{eqnarray}\label{3.3}
()_{hom}:\RR^f\to G_\Z(f),\quad \varphi\mapsto (\varphi)_{hom}.
\end{eqnarray}
Let $\Aut_M:=\Aut((M,0),\circ,e,E)$ be the group of automorphisms
of the germ $(M,0)$ as a germ of an F-manifold with Euler field.
It is a finite group because $M$ is a massive F-manifold
\cite[Theorem 4.14]{He02}. We claim that there is also a natural
homomorphism
\begin{eqnarray}\label{3.4}
()_M:\RR^f\to \Aut_M,\quad \varphi\mapsto (\varphi)_M.
\end{eqnarray}
It arises as follows. $F\circ\varphi^{-1}$ is a universal
unfolding of $f$ with the same base space $(M,0)$ as the universal
unfolding $F$. A morphism which induces $F\circ \varphi^{-1}$ by $F$
is given by a pair $(\Phi,(\varphi)_M)$ where $(\varphi)_M\in\Aut_M$ and where
$\Phi:(\C^{n+1}\times M,0)\to(\C^{n+1}\times M,0)$ is a 
biholomorphic map germ with 
\begin{eqnarray}\label{3.5}
F\circ \varphi^{-1}=F\circ\Phi,\quad \Phi|_{(\C^{n+1}\times\{0\},0)}=\id,
\quad \pr_M\circ \Phi = (\varphi)_M\circ \pr_M.
\end{eqnarray}
Here $\Phi$ is not unique, but $(\varphi)_M$ is unique because
$\Aut_M$ is finite and the differential of $(\varphi)_M$ at $T_0M$
is determined by the action of $\varphi$ on the Jacobi algebra 
$\OO_{\C^{n+1},0}/J_f\cong T_0M$. The map $\Phi\circ\varphi$
satisfies 
\begin{eqnarray}\label{3.6}
F\circ (\Phi\circ\varphi)=F\quad\textup{and}\quad
\pr_M\circ (\Phi\circ\varphi)=(\varphi)_M\circ \pr_M
\end{eqnarray} 
and is an extension of the symmetry $\varphi$ of $f$ to a symmetry of $F$.

The following theorem is contained in \cite[Theorem 13.11]{He02}
and is rewritten in \cite[Theorem 6.1]{He11}.

\begin{theorem}\label{t3.1}
As above, let $f:(\C^{n+1},0)\to(\C,0)$ be an isolated hypersurface singularity,
and let $F:(\C^{n+1}\times M,0)\to(\C,0)$ be a universal unfolding 
with base space $(M,0)$.

\medskip
(a) The homomorphism $()_M:\RR^f\to\Aut_M$ factors through $R_f$ 
to a homomorphism $()_M:R_f\to\Aut_M$. 
If $\mult f\geq 3$ then $()_M:R_f\to\Aut_M$ is 
an isomorphism and then $j_1\RR^f=R_f$.
If $\mult f=2$ then $()_M:R_f\to\Aut_M$ is surjective with kernel of
order 2. If $f=g(x_0,...,x_{n-1})+x_n^2$ then the kernel is generated by
the class of the symmetry $(x\mapsto (x_0,...,x_{n-1},-x_n))$.

\medskip
(b) The homomorphism $()_{hom}:\RR^f\to G_\Z(f)$ factors through $R_f$
to an injective homomorphism $()_{hom}:R_f\to G_\Z(f)$.
Let $G^{smar}_\RR(f)\subset G_\Z(f)$ denote its image.
It contains $-\id$ if and only if $\mult f=2$.  
If $f=g(x_0,...,x_{n-1})+x_n^2$ then
$-\id=(x\mapsto (x_0,...,x_{n-1},-x_n))_{hom}$.

\medskip
(c) The homomorphism 
\begin{eqnarray}\label{3.7}
()_M\circ ()_{hom}^{-1}:G^{smar}_\RR(f)\to\Aut_M
\end{eqnarray}
is an isomorphism if $\mult f\geq 3$. It is surjective with kernel
$\{\pm\id\}$ if $\mult f=2$.
\end{theorem}

Consider a singularity $f:(\C^{n+1},0)\to(\C,0)$ and a good representative
$F:\YY\to T$ of a universal unfolding with base space $M$.
One can choose it such that any element of $R_f$ 
lifts to an automorphism of $F$. 
Consider the $\Z$-lattice bundle 
\begin{eqnarray}\label{3.8}
\bigcup_{(\tau,t)\in T\times M-\DD} H_n(F_t^{-1}(\tau),\Z).
\end{eqnarray}

\begin{definition/lemma}\label{t3.2}
(a) (Definition)
We call the flat extension of the $\Z$-lattice bundle in \eqref{3.8} 
to $\C\times M-\DD$ the canonical $\Z$-lattice bundle of $M$.

\medskip
(b) (Lemma) Any element of $\Aut_M$ lifts to an automorphism
of the canonical $\Z$-lattice bundle of $M$.
The lift is unique up to $\pm 1$.
\end{definition/lemma}

{\bf Proof:}
The surjectivity of the homomorphism $()_M:R_f\to \Aut_M$ 
implies that any automorphism of $((M,0),\circ,e,E)$ lifts
to an automorphism of the unfolding and thus to an automorphism
of the $\Z$-lattice bundle in \eqref{3.10}.
Because of lemma \ref{t2.3}, the only automorphisms of the $\Z$-lattice
bundle which fix $T\times M$, are $\pm\id$. 
Therefore any element of $\Aut_M$ has only two lifts, and they 
differ by $\pm 1$.
\hfill$\Box$ 

\bigskip
Part (b) justifies part (a): The bundle depends only on the isomorphism
class of the germ $((M,0),\circ,e,E)$.  
Instead of lemma \ref{t2.3}, we could have used theorem \ref{t3.1} (c).
But that would give only that any element of $\Aut_M$ has two canonical
lifts, which differ by $\pm 1$, not that they are the only lifts.

In the case of a quasihomogeneous singularity, 
the finite group of quasihomogeneous symmetries
is a natural lift of $R_f$. This will be useful for the calculations
in section \ref{s5}.

\begin{theorem}\label{t3.3}\cite[Theorem 13.13]{He02}
\cite[Theorem 6.2]{He11}
Let $f\in \C[x_0,...,x_n]$ be a quasihomogeneous polynomial with an isolated singularity
at $0$ and weights $w_0,...,w_n\in\Q\cap(0,\frac{1}{2}]$ and weighted degree $1$.
Suppose that $w_0\leq ...\leq w_{n-1}<\frac{1}{2}$ (then $f\in \mmm^3$
if and only if $w_n<\frac{1}{2})$. Let $G_w$ be the algebraic group of 
quasihomogeneous coordinate changes, that means, those which respect $\C[x_0,...,x_n]$
and the grading by the weights $w_0,...,w_n$ on it. Then
$$R_f \cong \Stab_{G_w}(f).$$
\end{theorem}

\begin{remark}\label{t3.4}
Let $f\in\C[x_0,...,x_n]$ be a quasihomogeneous polynomial
with an isolated singularity at 0 and weights
$w_0,...,w_n\in(0,\frac{1}{2}]$ and weighted degree 1.
Then
\begin{eqnarray*}
\varphi_1:=(x\mapsto (e^{2\pi iw_0}x_0,...,e^{2\pi iw_n}x_n))
\in\textup{Stab}_{G_{\bf w}}(f)
\end{eqnarray*}
satisfies $(\varphi_1)_{hom}=M_h$. 
Now let $F(x,t)=f(x)+\sum_{i=1}^\mu t_im_i$
be a universal unfolding as in \eqref{2.9} with 
$m_i$ a weighted homogeneous polynomial of weighted degree
$\deg_{\bf w}m_i$. Then $\deg_{\bf w}t_i:=1-\deg_{\bf w}m_i$,
\begin{eqnarray*}
(\varphi_1)_M=(t\mapsto (e^{2\pi i \deg_{\bf w}t_1}t_1,...,
e^{2\pi i \deg_{\bf w}t_\mu}t_\mu)),
\end{eqnarray*}
and the pair $(\Phi_1,(\varphi_1)_M)$ with 
\begin{eqnarray*}
\Phi_1 = ((x,t)\mapsto (x,(\varphi_1)_M)
\end{eqnarray*}
induces $F\circ \varphi_1^{-1}$ by $F$, i.e.
\eqref{3.6} holds:
\begin{eqnarray*}
F\circ (\Phi_1\circ\varphi_1)=F,\quad 
\pr_M\circ (\Phi_1\circ\varphi_1) = (\varphi_1)_M\circ\pr_M.
\end{eqnarray*}
Especially, $M_h\in G^{smar}_\RR(f)$ and
$()_M\circ ()_{hom}^{-1}(M_h) =(\varphi_1)_M$.
\end{remark}

\subsection{Marked singularities and their moduli spaces}\label{s3.2}
In \cite{He11} the notion of a {\it marked} singularity was
defined and results from \cite{He02} on moduli spaces of right equivalence
classes of singularities were lifted to marked singularities. 
Here we recall the central notions and results.

\begin{definition}\label{t3.5}
Fix one reference singularity $f^{(0)}:(\C^{n+1},0)\to(\C,0)$.

\medskip
(a) Its $\mu$-homotopy class is the set of all singularities
$f:(\C^{n+1},0)\to(\C,0)$ for which a $C^\infty$-family $f_s$, $s\in[0,1]$,
of singularities with $\mu(f_s)=\mu(f^{(0)})$ and $f_0=f^{(0)}$ and $f_1=f$ exists.

\medskip
(b) A {\it marked} singularity is a pair $(f,\pm\rho)$ with $f$ in the
$\mu$-homotopy class of $f^{(0)}$ and $\rho:Ml(f,L)\to Ml(f^{(0)},L)$ an isomorphism
between Milnor lattices with Seifert forms. Here $\pm\rho$ means
the set $\{\rho,-\rho\}$, so neither $\rho$ nor $-\rho$ is preferred.

\medskip
(c) Two singularities $f_1$ and $f_2$ are right equivalent if a 
coordinate change $\varphi\in\RR$ with $f_1=f_2\circ \rho$ exists.
Notation: $f_1\sim_\RR f_2$.

Two marked singularities $(f_1,\pm\rho_1)$ and $(f_2,\pm\rho_2)$ 
are right equivalent if a coordinate change $\varphi\in\RR$ with
$f_1=f_2\circ \varphi$ and $\rho_1=\varepsilon\cdot \rho_2\circ(\varphi)_{hom}$
for some $\varepsilon\in\{\pm 1\}$ exists.
Notation: $(f_1,\pm\rho_1)\sim_\RR(f_2,\pm\rho_2)$.

\medskip
(d) The moduli spaces $M_\mu(f^{(0)})$ and $M_\mu^{mar}(f^{(0)})$ are defined
as the sets
\begin{eqnarray}\label{3.9}
M_\mu(f^{(0)})&:=& \{\textup{the }\mu\textup{-homotopy class of }f^{(0)}\}/\sim_\RR,\\
M_\mu^{mar}(f^{(0)})&:=& \{\textup{the marked singularities }(f,\pm\rho)\}/\sim_\RR.
\label{3.10}
\end{eqnarray}
\end{definition}

A central result in \cite{He02} is that the moduli space $M_\mu(f^{(0)})$ 
has the structure of an analytic geometric quotient.
In \cite{He11} this result is extended to the space $M_{\mu}^{mar}(f^{(0)})$,
and it is shown that $M_{\mu}^{mar}(f^{(0)})$ is locally isomorphic
to a {\it $\mu$-constant stratum}. This is recalled in theorem \ref{t3.6} below.

Let $f:(\C^{n+1},0)\to(\C,0)$ be a singularity and let
$F:(\C^{n+1}\times M(f),0)\to(\C,0)$ be a universal unfolding of $f$
with base space $(M(f),0)$. Then the $\mu$-constant stratum is the
germ $(S_\mu(f),0)\subset (M(f),0)$ of the subset 
\begin{eqnarray}\label{3.11}
S_\mu(f):=\{t\in M\, |\, F_t \textup{ has only one singularity }x_0
\textup{ and }F_t(x_0)=0\}.
\end{eqnarray}
Here $F:\YY\to T$ is a good representative of the germ $F$
with base space $M$. Obviously $(S_\mu(f),0)$ carries a natural structure 
as a reduced complex space germ. 
In \cite[Theorem]{He02} it is equipped furthermore with a 
natural complex structure, which is not necessarily reduced.

\begin{theorem}\label{t3.6}
Fix one reference singularity $f^{(0)}$. 

\medskip
(a) (\cite[Theorem 13.11]{He02} and \cite[Theorem 4.3]{He11})
$M_\mu(f^{(0)})$ and $M_\mu^{mar}(f^{(0)})$ are in a 
natural way complex spaces. 
They can be constructed with the underlying reduced complex structures
as analytic geometric quotients.

\medskip
(b) For any $f$ in the $\mu$-homotopy class of $f^{(0)}$, the germ
$(M_\mu(f^{(0)}),[f])$ is isomorphic with its canonical complex structure
to the quotient $(S_\mu(f),0)/\Aut_{M(f)}$ (the action of $\Aut_{M(f)}$ on $(M(f),0)$
restricts to an action on $(S_\mu(f),0)$).

For any marked singularity $(f,\pm \rho)$, the germ 
$(M_\mu^{mar}(f^{(0)}),[(f,\pm\rho)])$ is isomorphic with its canonical complex
structure to the $\mu$-constant stratum $(S_\mu(f),0)$.

\medskip
(c) For any $\chi\in G_\Z(f^{(0)})$, the map
\begin{eqnarray}\label{3.12}
\chi_{mar}:M_\mu^{mar}(f^{(0)})&\to& M_\mu^{mar}(f^{(0)}),\\
{}[(f,\pm\rho)]&\mapsto& [(f,\pm\chi\circ\rho)],\nonumber 
\end{eqnarray}
is an automorphism of $M_\mu^{mar}(f^{(0)})$. The action
\begin{eqnarray}\label{3.13}
G_\Z(f^{(0)})\times M_\mu^{mar}(f^{(0)})&\to& M_\mu^{mar}(f^{(0)}),\\
(\chi,[(f,\pm\rho)])&\mapsto& \chi_{mar}([(f,\pm\rho)])=[(f,\pm\chi\circ\rho)],\nonumber
\end{eqnarray}
is a group action from the left.
It is properly discontinuous. The quotient $M_\mu^{mar}(f^{(0)})/G_\Z(f^{(0)})$
is the moduli space $M_\mu(f^{(0)})$ of unmarked singularities,
with its canonical complex structure.

\medskip
(d) (Definition) Recall the definition of $G_\RR^{smar}(f)$ in theorem \ref{t3.1} (b).
Define 
\begin{eqnarray}\label{3.14}
G^{mar}_\RR(f):=\{\pm\psi\, |\, \psi\in G^{smar}_\RR(f)\}\subset G_\Z(f),
\end{eqnarray}
Remark: By theorem \ref{t3.1} this is equal to $G^{smar}_\RR(f)$ if $\mult(f)=2$ and has $G^{smar}_\RR(f)$
as index 2 subgroup if $\mult(f)\geq 3$.

\medskip
(e) For any point $[(f,\pm\rho)]\in M_\mu^{mar}(f^{(0)})$, 
the stabilizer of it in $G_\Z(f^{(0)})$ is the finite group
\begin{eqnarray}\label{3.15}
\rho\circ G^{mar}_\RR(f)\circ \rho^{-1}\subset G_\Z(f^{(0)}).
\end{eqnarray}

\end{theorem}

\begin{remarks}\label{t3.7}
(i) In \cite{He11} also the notion of a {\it strongly marked} singularity
is defined. It is a pair $(f,\rho)$ with $f$, $\rho$ and 
$f^{(0)}$ as in definition 
\ref{t3.5} (b). The moduli space $M_\mu^{smar}(f^{(0)})$ 
of strongly marked singularities behaves
as well as $M_\mu^{mar}(f^{(0)})$ if the following holds:
Either any singularity in the $\mu$-homotopy class of $f^{(0)}$ has multiplicity $\geq 3$,
or any singularity in the $\mu$-homotopy class of $f^{(0)}$ has multiplicity 2.
We expect that to hold, but we don't know it. If it does not hold
then $M_\mu^{smar}(f^{(0)})$ is not Hausdorff, see \cite[Theorem 4.3 (e)]{He11}.
We do not need strongly marked singularities here.

\medskip
(ii) In \cite{He11}\cite{GH17-1} and \cite{GH17-2} the moduli space
$M_\mu^{mar}(f^{(0)})$ for any singularity $f^{(0)}$ with modality $\leq 2$
was studied. It turned out that almost all of them are connected, but not all,
namely not those for $f^{(0)}$ in the subseries 
$W_{1,12r}^\sharp,S_{1,10r}^\sharp,U_{1,9r},E_{3,18r},
Z_{1,14,r},Q_{2,12r},W_{1,12r},S_{1,10r}$
of the eight bimodal series. 
This disproved the conjecture 3.2 (a) in \cite{He11}
that $M_\mu^{mar}(f^{(0)})$ is connected for any singularity.
But for all other singularities $f^{(0)}$ with modality $\leq 2$
$M_\mu^{mar}(f^{(0)})$ is connected. 
An equivalent statement to $M_\mu^{mar}(f^{(0)})$ connected is 
because of \cite[Theorem 4.4 (a)]{He11} that any element of $G_\Z(f^{(0)})$
arises from geometry, namely it is 
$(\pm 1)\cdot$the transversal monodromy of 
a suitable $\mu$-constant family $f_s$, $s\in[0,1]$, with 
$f_0=f_1=f^{(0)}$ (a $C^\infty$ family of singularities $f_s$ with
$\mu(f_s)=\mu(f^{(0)})$).

\medskip
(iii) In the case of the ADE-singularities, which have modality 0,
$M_\mu^{mar}(f^{(0)})$ is simply a point \cite{He11}. In the case of the simple
elliptic singularities, which have modality one and which are 
parametrized by elliptic curves, $M_\mu^{mar}(f^{(0)})$ is isomorphic
to $\H$ \cite{GH17-1}. In both cases, the connectedness of
$M_\mu^{mar}(f^{(0)})$ will be important in the proof of the main
theorem in section \ref{s7}.
\end{remarks}

\subsection{A thickening of the moduli space of marked singularities}\label{s3.3}
Fix one reference singularity $f^{(0)}$. By theorem \ref{t3.6} (b), locally
at $[(f,\pm\rho)]$, the moduli space $M_\mu^{mar}(f^{(0)})$ 
is isomorphic to the $\mu$-constant stratum $(S_\mu(f),0)\subset (M(f),0)$
of the singularity $f$. In \cite{GH18} we will show the following.

\begin{theorem}\label{t3.8}
Fix one reference singularity $f^{(0)}$. 

\medskip
(a) The moduli space $M_\mu^{mar}(f^{(0)})$ of marked singularities can 
be extended globally to a $\mu$-dimensional F-manifold 
$M^{mar}(f^{(0)})\supset M_\mu^{mar}(f^{(0)})$ with the following properties.

\begin{list}{}{}
\item[(i)]
For any point $[(f,\pm\rho)]\in M_\mu^{mar}(f^{(0)})$, a certain neighborhood 
$U_{[(f,\pm\rho)]}$ of it in $M^{mar}(f^{(0)})$ 
and an isomorphism $\psi_{[(f,\pm\rho)]}: M(f)\to U_{[(f,\pm\rho)]}$ 
of $F$-manifolds exist, 
where $M(f)$ is the base space of a good representative 
of a universal unfolding of $f$.
\item[(ii)]
$M^{mar}(f^{(0)})$ is covered by these neighborhoods $U_{[(f,\pm\rho)]}$.
\item[(iii)]
The action of $G_\Z(f^{(0)})$ on $M_\mu^{mar}(f^{(0)})$ extends to an action
of $G_\Z(f^{(0)})$ on this F-manifold with Euler field, and the map
\begin{eqnarray}\label{3.16}
G_\Z(f^{(0)})\to \Aut(M^{mar}(f^{(0)}),\circ,e,E)
\end{eqnarray}
is surjective with kernel $\{\pm \id\}$.
\end{list}

\medskip
(b) Let $\DD^{mar}\subset \C\times M^{mar}(f^{(0)})$ be the discriminant
\begin{eqnarray*}
\DD^{mar}:=\{(\tau,t)\in \C\times M^{mar}(f^{(0)})\, |\, 
E\circ:T_tM\to T_tM\textup{ has eigenvalue }\tau\}.
\end{eqnarray*}
It is a hypersurface. 
$M^{mar}(f^{(0)})$ comes equipped with a $\Z$-lattice bundle
$H_\Z\to (\C\times M^{mar}(f^{(0)})-\DD^{mar})$ of rank $\mu$
with the following properties.

\begin{list}{}{}
\item[(i)]
For any point  $[(f,\pm\rho)]\in M_\mu^{mar}(f^{(0)})$
and a good representative $F:\YY\to T$ of a universal unfolding of $f$ 
with base space $M(f)$ and the isomorphism 
$\psi_{[(f,\pm\rho)]}:M(f)\to U_{[(f,\pm\rho)]}$ as in (a)(i),
this isomorphism lifts to an isomorphism from the canonical
$\Z$-lattice bundle above $M(f)$ in definition \ref{t3.2} (a) 
to the restriction of $H_\Z$ to  $\C\times U_{[(f,\pm\rho)]}-\DD^{mar}$.
(Because of lemma \ref{t3.2} (b), this lift is unique up to $\pm 1$).
\item[(ii)]
Let $r:M^{mar}(f^{(0)})\to\R_{>0}$ be a $C^\infty$ function with 
$\DD^{mar}\subset \bigcup_{t\in M^{mar}(f^{(0)})}\Delta_{r(t)}\times \{t\}$.
Then the restriction of $H_\Z$ to 
$\bigcup_{t\in M^{mar}(f^{(0)})}\R_{>r(t)}\times \{t\}$ has trivial monodromy,
i.e. it is a trivial flat $\Z$-lattice bundle.
\item[(iii)]
Write $t^{(0)}:= [(f^{(0)},\pm\id)]\in M^{mar}_\mu(f^{(0)})$.
For any $t=[(f,\pm\rho)]\in M^{mar}_\mu(f^{(0)})$  and any small $\tau>0$, 
the following diagram of isomorphisms commutes,
\begin{eqnarray}\label{3.17}
\begin{xy}
\xymatrix{ H_{\Z,(\tau,t)} \ar[r]^{\cong}_{(i)}  \ar[d]_{\cong}^{(ii)}
& Ml(f)  \ar[d]_{\cong}^{\pm\rho}\\
H_{\Z,(\tau,t^{(0)})} \ar[r]^{\cong}_{(i)} & Ml(f^{(0)})  }
\end{xy}
\end{eqnarray}

\item[(iv)]
The action of $G_\Z(f^{(0)})$ on $M^{mar}(f^{(0)})$ extends to an action on the
$\Z$-lattice bundle $H_\Z$ (the action of $G_\Z(f^{(0)})$ on the first factor $\C$
of $\C\times M^{mar}(f^{(0)})$ is trivial by definition).
\end{list}
%\medskip
%(c) For any $t\in M^{mar}(f^{(0)})$ and any $\tau\in\R_{>0}$ with
%$\DD^{mar}\cap \C\times\{t\}\subset \Delta_\tau\times\{t\}$ 
%an isomorphism $\rho_{(\tau,t)}:H_{\Z,(\tau,t)}\to Ml(f^{(0)})$ exists
%such that the family of all isomorphisms glues modulo $\pm\id$
%to a flat family of isomorphisms. 
%This gives a family of markings modulo $\pm\id$ on $M^{mar}(f^{(0)})$.
%
%And for $(\tau,t)=(\tau,[(f,\pm\rho)])\in \R_{>0}\times M^{mar}_\mu(f^{(0)}) 
%\subset\R_{>0}\times M^{mar}(f^{(0)})$, the isomorphism $\rho_{(\tau,t)}$
%is equal to $\pm\rho$. That means, the family of markings modulo $\pm\id$ 
%on $M^{mar}_\mu(f^{(0)})$ extends to the family of markings modulo $\pm\id$ on 
%$M^{mar}(f^{(0)})$.
\end{theorem}

As the paper \cite{GH18} is not yet available, we will prove this theorem for
the cases of interest here, the simple and the simple elliptic singularities, 
directly in section \ref{s4}. See theorem \ref{t4.3}.

$M^{mar}(f^{(0)})$ contains besides $\DD^{mar}$ also two other hypersurfaces,
the caustic $\KK_3^{mar}$ and the Maxwell stratum $\KK_2^{mar}$,
\begin{eqnarray}
\KK_3^{mar}&:=& \{t\in M^{mar}(f^{(0)})\, |\, T_tM^{mar}(f^{(0)})
\textup{ decomposes into }\nonumber \\
&&\textup{less than }\mu\textup{ irreducible local algebras}\},\label{3.18}\\
\KK_2^{mar}&:=& \textup{the closure of the set }
\{t\in M^{mar}(f^{(0)})-\KK_3^{mar}\, |\, \textup{some} \nonumber \\
&&\textup{eigenvalues of }E\circ:T_tM\to T_tM\textup{ coincide}\}.\label{3.19}
\end{eqnarray}

\subsection{A Looijenga-Deligne map for distinguished bases}\label{s3.4}
Looijenga \cite{Lo74} studied in the case of the simple singularities
a relationship between distinguished bases and the base space of a universal
unfolding. His idea carries over to the F-manifold $M^{mar}(f^{(0)})$ 
in theorem \ref{t3.8} of an arbitrary $\mu$-homotopy class of singularities. 
We describe the idea here in this generality. 
In section \ref{s7} we will study it in the cases
of the simple and the simple elliptic singularities. 
In \cite{GH18} we will study it in the general case. 

The set $\BB(f)$ of distinguished bases of the 
Milnor lattice $Ml(f)$ of a singularity $f$ was constructed 
in subsection \ref{s2.5} 
by choosing {\it one} morsification $F_t$ of $f$
and considering {\it all possible} distinguished systems of paths.
Following Looijenga, now we want to fix {\it one} distinguished system of paths
and consider {\it all possible} morsifications. 
The following two definitions make this precise.

\begin{definition}\label{t3.9}
(a) Fix a tuple  $(u_1,...,u_\mu)\subset\C^\mu$ with $u_i\neq u_j$ for $i\neq j$.
The {\it good ordering} of it is the lexicographic ordering
$(u_{\sigma(1)},...,u_{\sigma(\mu)})$ by 
(imaginary part,$-$real part). That means, the corresponding {\it good} permutation
$\sigma\in S_\mu$ is uniquely determined by
\begin{eqnarray}\label{3.20}
i<j\iff
\left\{\begin{array}{ll}\Imm(u_{\sigma(i)})<\Imm(u_{\sigma(j)})&\textup{ or }\\
\Imm(u_{\sigma(i)})=\Imm(u_{\sigma(j)})&\textup{ and }
\Ree(u_{\sigma(i)})>\Ree(u_{\sigma(j)}). \end{array}\right.
\end{eqnarray}

(b) Fix a tuple $(u_1,...,u_\mu)$ as in (a) with good permutation $\sigma\in S_\mu$,
and fix additionally a $\tau\in \R_{>0}$ with $\tau>\max_i|u_i|$.
A {\it good} distinguished system of paths is a distinguished system
of paths $\gamma_1,...,\gamma_\mu$ such that $\gamma_j$ goes from $u_{\sigma(j)}$
to $\tau$.
\end{definition}

For a fixed tuple $(u_1,...,u_\mu,\tau)$ as above, 
all good distinguished systems of paths are homotopy equivalent
with respect to a natural notion of homotopy equivalence. 
And if $F_t$ is a morsification of a singularity $f$ with 
critical values $u_1,...,u_\mu$, all good distinguished systems 
of paths give the same distinguished basis up to the action 
of the sign group $G_{sign,\mu}$.
%The following picture gives an idea of a good 
%distinguished system of paths.
%Here the good permutation is $\sigma=\id$.
%PICTURE

%\bigskip

\begin{definition}\label{t3.10}
Fix one reference singularity $f^{(0)}$. 

\medskip
(a) The set of {\it Stokes walls}
within the space $M^{mar}(f^{(0)})$ in theorem \ref{t3.8} is the set
\begin{eqnarray}\label{3.21}
W_{Stokes}&:=& \{t\in M^{mar}(f^{(0)})\, |\, 
\textup{the eigenvalues }u_1,...,u_\mu\textup{ of }\\
&& E\circ:T_tM\to T_tM \textup{ satisfy }\Imm(u_i)=\Imm(u_j)\textup{ for some }
i\neq j\}\nonumber
\end{eqnarray}
The set $W_{Stokes}$ of Stokes walls is a real codimension 1 subvariety and
contains $\KK_3^{mar}\cup \KK_2^{mar}$. The components of its complement
$M^{mar}(f^{(0)})-W_{Stokes}$ are called {\it Stokes regions}.
Let $R_{Stokes}$ be the set of Stokes regions,
and let $R_{Stokes}^0$ be the subset of those Stokes regions 
which are in the component
$M^{mar}(f^{(0)})^0$ of $M^{mar}(f^{(0)})$ 
which contains $[(f^{(0)},\pm\id)]$. 

\medskip
(b) The set $\BB^{ext}(f^{(0)})$ is the orbit of $\BB(f^{(0)})$ under
the action of $G_\Z$. It contains (all?) $\Z$-bases $(\delta_1,...,\delta_\mu)$ 
of $Ml(f^{(0)})$ whose elements $\delta_j$ are vanishing cycles and such that
$s_{\delta_1}\circ ...\circ s_{\delta_\mu}=M_h$.
It consists of $G_{sign,\mu}\rtimes \textup{Br}_\mu$ orbits.
One of these orbits is the set $\BB(f^{(0)})$ of distinguished bases.

\medskip
(c) The {\it Looijenga-Deligne map} is the map
\begin{eqnarray}\label{3.22}
LD:R_{Stokes}\to \BB^{ext}(f^{(0)})/G_{sign,\mu}
\end{eqnarray}
which is defined as follows. 
For a Stokes region in $M^{mar}(f^{(0)})$, 
choose a point $t$ in it
and a point $[(f,\pm\rho)]\in M_\mu^{mar}(f^{(0)})$ with 
$t\in U_{[(f,\pm\rho)]}$. 
Let $(u_1,...,u_\mu)$ be the eigenvalues of $E\circ:T_tM\to T_tM$
in the good ordering (definition \ref{t3.9}). 
Consider a good distinguished system of paths from $(u_1,...,u_\mu)$ 
to a value $\tau>\max_i|u_i|$ (definition \ref{t3.9} (b)).
The usual construction of distinguished bases gives a distinguished
basis in $\BB(f)$ up to the action of the sign group $G_{sign,\mu}$.
Shift this basis with the isomorphism $\rho:Ml(f)\to Ml(f^{(0)})$ to an element of 
$\BB^{ext}(f^{(0)})/G_{sign,\mu}$.
\end{definition}

\begin{remarks}\label{t3.11}
(i) We claim that $LD$ restricts to a map
\begin{eqnarray}\label{3.23}
LD:R^0_{Stokes}\to \BB(f^{(0)})/G_{sign,\mu}.
\end{eqnarray}
We prove this by a different description of the image $LD(t)$ for 
$t\in R^0_{Stokes}$.
Choose $[(f,\pm\rho)]\in M^{mar}_\mu(f^{(0)})^0$ with 
$t\in U_{[(f,\pm\rho)]}$, 
choose $\tau>\max_i |u_i|$, and choose a good
distinguished system of paths from
$(u_1,...,u_\mu)$ to $\tau$. One can move $t$ within 
$M^{mar}(f^{(0)})^0-(\KK_3^{mar}\cup\KK_2^{mar})$ 
to a point in 
$U_{[(f^{(0)},\pm\id)]}\subset M^{mar}_\mu(f^{(0)})^0$. 
Then the good distinguished system of
paths moves to some new distinguished system of paths. 
Now the construction of distinguished bases for $f^{(0)}$ 
gives directly 
the class of bases $LD(t)\in\BB(f^{(0)})/G_{sign,\mu}$. 
This follows with \eqref{3.17}.

\medskip
(ii) The action of $G_\Z$ on $M^{mar}(f^{(0)})$ induces an action
on $R_{Stokes}$. And $R^0_{Stokes}=R_{Stokes}$ if and only if
$M^{mar}(f^{(0)})$ is connected.
The map \eqref{3.22} is $G_\Z$-equivariant.
Therefore, if $M^{mar}(f^{(0)})$ is connected, then \eqref{3.23} 
and the definition of $\BB^{ext}(f^{(0)})$ give also 
$\BB^{ext}(f^{(0)})=\BB(f^{(0)})$. 

\medskip
(iii) But if $M^{mar}(f^{(0)})$ is not connected, we do not know
whether $\BB^{ext}(f^{(0)})=\BB(f^{(0)})$ or
$\BB^{ext}(f^{(0)})\supsetneqq\BB(f^{(0)})$ holds. 
The first open cases are the subseries in remark \ref{t3.7} (ii) of the
eight bimodal series.

\medskip
(iv) Looijenga considered the map $LD$ for the simple singularities and
proved that it is an isomorphism for the $A_\mu$ singularities \cite{Lo74}.
Deligne \cite{De74} proved the same for the $D_\mu$ and $E_\mu$ singularities.
We will reprove their results and extend them to the simple elliptic 
singularities in section \ref{s7}. We will study $LD$ in the general
case in \cite{GH18}.
\end{remarks}

\section{Unfoldings of the simple and the
simple elliptic singularities}\label{s4}
\setcounter{equation}{0}

\noindent
The first singularities in Arnold's lists of isolated hypersurface singularities
\cite[ch 15.1]{AGV85} are the simple and the simple elliptic singularities.
They are distinguished by many properties. Especially, they possess
universal unfoldings such that all members are defined globally on $\C^{n+1}$. 
We start by giving well known normal forms. Then we choose universal unfoldings.

\subsection{Normal forms for the simple and the simple elliptic singularities}\label{s4.1}

The first table lists normal forms from \cite{AGV85} 
for the simple singularities $f:\C^{n+1}\to\C$,
\begin{eqnarray}\label{4.1}
\begin{array}{cccr}
\textup{name} & \mu & n & f(x_0,...,x_n)\\ \hline
A_\mu & \geq 1 & \geq 0 & x_0^{\mu+1}+\sum_{i=1}^nx_i^2 \\
D_\mu & \geq 4 & \geq 1 & x_0^{\mu-1}+x_0x_1^2+\sum_{i=2}^nx_i^2 \\
E_6 & 6 & \geq 1 & x_0^4+x_1^3+\sum_{i=2}^nx_i^2 \\
E_7 & 7 & \geq 1 & x_0^3x_1+x_1^3+\sum_{i=2}^nx_i^2 \\ 
E_8 & 8 & \geq 1 & x_0^5+x_1^3+\sum_{i=2}^nx_i^2 
\end{array}
\end{eqnarray}

The simple elliptic singularities can be represented as 1-parameter families
in different ways \cite[1.9 and 1.11]{Sa74}\cite[ch 15.1]{AGV85}.
We choose the Legendre normal forms $f=f_\lambda:\C^{n+1}\to\C$ 
from \cite[1.9]{Sa74} in the following table 
with $\lambda\in\C-\{0,1\}$,
\begin{eqnarray}\label{4.2}
\begin{array}{cccr}
\textup{name} & \mu & n & f_\lambda(x_0,...,x_n)\\ \hline
\www E_6 & 8 & \geq 2 & x_1(x_1-x_0)(x_1-\lambda x_0)-x_0x_2^2 +\sum_{i=3}^nx_i^2 \\
\www E_7 & 9 & \geq 1 & x_0x_1(x_1-x_0)(x_1-\lambda x_0) +\sum_{i=2}^nx_i^2 \\
\www E_8 & 10 & \geq 1 & x_1(x_1-x_0^2)(x_1-\lambda x_0^2) +\sum_{i=2}^nx_i^2 
\end{array}
\end{eqnarray}

\subsection{Universal unfoldings}\label{s4.2}

For the simple singularities, we reproduce the universal unfoldings which are
given in \cite{Lo74}. They are as follows, 
\begin{eqnarray}\label{4.3}
&&F^{alg}:\C^{n+1}\times M^{alg}\to\C\quad\textup{with}\quad M^{alg}=\C^\mu,\\
&&F^{alg}(x_0,...,x_n,t_1,...,t_\mu)=F^{alg}(x,t)=F^{alg}_t(x)
=f(x)+\sum_{j=1}^\mu t_jm_j \nonumber
\end{eqnarray}
with $f=F^{alg}_0$ and $m_1,...,m_\mu$ the monomials in the tables
\eqref{4.4} and \eqref{4.5},
\begin{eqnarray}\label{4.4}
\begin{array}{ccccccc}
\textup{name}  & m_1 & m_2 & m_3 & m_4 & ... & m_\mu \\ \hline 
A_\mu & 1 & x_0 & x_0^2 & x_0^3 & ... & x_0^{\mu-1} \\
D_\mu & 1 & x_1 & x_0 & x_0^2 & ... & x_0^{\mu-2} 
\end{array}
\end{eqnarray}
\begin{eqnarray}\label{4.5}
\begin{array}{ccccccccc}
\textup{name}  & m_1 & m_2 & m_3 & m_4 & m_5 
& m_6 & m_7 & m_8  \\ \hline
E_6 & 
1 & x_0 & x_1 & x_0^2 & x_0x_1 & x_0^2x_1 & &   \\
E_7 & 
1 & x_0 & x_1 & x_0^2 & x_0x_1 & x_0^3 & x_0^4 &  \\
E_8 & 
1 & x_0 & x_1 & x_0^2 & x_0x_1 & x_0^3 & x_0^2x_1 & x_0^3x_1 
\end{array}
\end{eqnarray}
One checks easily that the monomials form a basis 
of the Jacobi algebra $\OO_{\C^{n+1},0}/J_f$.
Therefore the unfolding $F^{alg}$ is indeed universal
(compare \eqref{2.9}).

For each of the three Legendre families of simple elliptic singularities, 
we give a global family of functions as follows,
\begin{eqnarray}
&&F^{alg}:\C^{n+1}\times M^{alg}\to\C \quad\textup{with}\quad
M^{alg}=\C^{\mu-1}\times(\C-\{0,1\}),\nonumber\\
&&F^{alg}(x_0,...,x_n,t_1,...,t_{\mu-1},\lambda)
=F^{alg}(x,t',\lambda)\\
&&=F^{alg}_{t',\lambda}(x) \nonumber 
=f_\lambda(x)+\sum_{j=1}^{\mu-1} t_jm_j \label{4.6}
\end{eqnarray}
with $f_\lambda=F^{alg}_{0,\lambda}$ and $m_1,...,m_{\mu-1}$ 
the monomials in the table \eqref{4.7},
\begin{eqnarray}\label{4.7}
\begin{array}{cccccccccc}
\textup{name}  & m_1 & m_2 & m_3 & m_4 & m_5 
& m_6 & m_7 & m_8 & m_9 \\ \hline
\www E_6 &
1 & x_0 & x_1 & x_2 & x_0^2 & x_0x_1 & x_1x_2 & & \\
\www E_7 &
1 & x_0 & x_1 & x_0^2 & x_0x_1 & x_1^2 & x_0^2x_1 & x_0x_1^2 &  \\
\www E_8 & 
1 & x_0 & x_0^2 & x_1 & x_0^3 & x_0x_1 & x_0^2x_1 & x_1^2 & x_0x_1^2 
\end{array}
\end{eqnarray}
Let $\lambda:\H\to\C-\{0,1\},t_\mu\mapsto \lambda(t_\mu)$, be the standard
universal covering. For each of the three Legendre families of 
simple elliptic singularities, we will also 
consider the global family of functions
\begin{eqnarray}\label{4.8}
&&F^{mar}:\C^{n+1}\times M^{mar}\to\C, (x,t)\mapsto F^{alg}(x,t',\lambda(t_\mu))\\
&&\textup{where}\quad M^{mar}=\C^{\mu-1}\times\H.\nonumber
\end{eqnarray}
For the simple singularities, we set
\begin{eqnarray}\label{4.9}
M^{mar}=M^{alg}=\C^\mu,\quad \lambda:=t_\mu,\quad F^{mar}=F^{alg}.
\end{eqnarray}

\begin{lemma}\label{t4.1}
Consider any of the three global families of functions in \eqref{4.6}.
At each point $(0,0,\lambda)\in\C^{n+1}\times\C^{\mu-1}\times(\C-\{0,1\})$, 
the germ of the family $F^{alg}$ is a universal unfolding of $f_\lambda$.

Also, at each point $(0,0,t_\mu)\in\C^{n+1}\times \C^{\mu-1}\times \H$,
the germ of the family $F^{mar}$ in \eqref{4.8} 
is a universal unfolding of $f_{\lambda(t_\mu)}$.
\end{lemma}

{\bf Proof:}
It suffices to prove the statement for $F^{alg}$. 
Because of \eqref{2.9}, it suffices to show
that for any $\lambda\in\C-\{0,1\}$ the monomials $m_1,...,m_{\mu-1}$
together with the weighted homogeneous polynomial
$\frac{\paa f_\lambda}{\paa\lambda}$ form a basis of the Jacobi algebra
$\OO_{\C^{n+1},0}/J_{f_\lambda}$. 
We carry out the least trivial case, which is the case $\www E_6$,
and leave the cases $\www E_7$ and $\www E_8$ to the reader.
In the case $\www E_6$, we work with the minimal number $n+1=3$ of 
variables. The normalized weight system is
${\bf w}=(w_0,w_1,w_2)=(\frac{1}{3},\frac{1}{3},\frac{1}{3})$,
and $\deg_{\bf w}f_\lambda=1$. As $f$ is quasihomogeneous, the Jacobi
algebra is isomorphic to $\C[x]/J_{f_\lambda}^{pol}$ where $J_{f_\lambda}^{pol}$
denotes the Jacobi ideal in $\C[x_0,x_1,x_2]=\C[x]$. For $q\in\Q_{\geq 0}$
denote by $\C[x]_q$ the sub vector space of $\C[x]$ 
generated by the monomials of weighted degree q.
\begin{eqnarray*}
\begin{array}{lllll}
\frac{\paa f_\lambda}{\paa x_0}&=& -(\lambda+1)x_1^2+2\lambda\cdot x_0x_1-x_2^2
&\in& J_{f_\lambda}^{pol}\cap\C[x]_{2/3},\\
\frac{\paa f_\lambda}{\paa x_1}&=& 3x_1^2-2(\lambda+1)\cdot x_0x_1+\lambda\cdot x_0^2
&\in& J_{f_\lambda}^{pol}\cap\C[x]_{2/3},\\
\frac{\paa f_\lambda}{\paa x_2}&=& -2x_0x_2
&\in& J_{f_\lambda}^{pol}\cap\C[x]_{2/3},\\
\frac{\paa f_\lambda}{\paa \lambda}&=& x_0^2x_1-x_0x_1^2
&\in& \C[x]_{1}.
\end{array}
\end{eqnarray*}
We have to show for $q\in \Q_{\geq 0}$
\begin{eqnarray}\label{4.10}
(J_{f_\lambda}^{pol})\cap \C[x]_q+\sum_{j:\, \deg_{\bf w}m_j=q}\C\cdot m_j
+\left\{\begin{array}{ll}0&\textup{if }q\neq 1\\ 
\C\cdot\frac{\paa f_\lambda}{\paa \lambda}&\textup{if }q=1\end{array}\right\} =\C[x]_q.
\end{eqnarray}
The only nontrivial cases are $q\in\{\frac{2}{3},1,\frac{4}{3}\}$. 

The case $q=\frac{2}{3}$:
$\C[x]_{2/3}$ is generated by the six monomials 
$x_0^2,x_0x_1,x_1^2,x_0x_2,x_1x_2,x_2^2$.
Here $m_5=x_0^2,m_6=x_0x_1,m_7=x_1x_2,$ and 
$\frac{\paa f_\lambda}{\paa x_2}= -2x_0x_2$. Modulo these four monomials,
$\frac{\paa f_\lambda}{\paa x_0}$ and $\frac{\paa f_\lambda}{\paa x_1}$
are congruent to $-(\lambda+1)x_1^2-x_2^2$ and $3x_1^2$. Thus the left hand
side of \eqref{4.10} contains also the monomials $x_1^2$ and $x_2^2$.

The case $q=1$: 
$\C[x]_1$ is generated by the 10 monomials
$x_0^3,x_0^2x_1,x_0x_1^2,x_1^3,x_0^2x_2,x_0x_1x_2,x_1^2x_2,x_0x_2^2,x_1x_2^2,x_2^3$.
The partial derivatives 
$x_0\frac{\paa f_\lambda}{\paa x_2},x_1\frac{\paa f_\lambda}{\paa x_2},
x_2\frac{\paa f_\lambda}{\paa x_2},x_2\frac{\paa f_\lambda}{\paa x_1}$ and
$x_2\frac{\paa f_\lambda}{\paa x_0}$ generate the monomials
$x_0^2x_2,x_0x_1x_2,x_0x_2^2,x_1^2x_2$ and $x_2^3$. For $\lambda\in\C-\{0,1\}$,
the polynomials $x_0^2x_1-x_0x_1^2$ and $x_0\frac{\paa f_\lambda}{\paa x_0}$
(and $x_0x_2^2$) generate the monomials $x_0^2x_1$ and $x_0x_1^2$.
Modulo these seven monomials, the three remaining partial derivatives  
$x_1\frac{\paa f_\lambda}{\paa x_0},x_0\frac{\paa f_\lambda}{\paa x_1},
x_1\frac{\paa f_\lambda}{\paa x_1}$ are congruent to
$-(\lambda+1)x_1^3-x_1x_2^2, \lambda x_0^3, 3x_1^3$.
Thus also the three remaining monomials 
$x_0^3,x_1^3,x_1x_2^2$ are in the 
left hand side of \eqref{4.10}. 

The case $q=\frac{4}{3}$: 
The ideal $J^{pol}_{f_\lambda}$ contains 
$\frac{\paa f_\lambda}{\paa x_2}$, so $x_0x_2$,
and $x_2\frac{\paa f_\lambda}{\paa x_1}$, so $x_1^2x_2$,
and $x_2\frac{\paa f_\lambda}{\paa x_0}$, so $x_2^3$.
Thus $J^{pol}_{f_\lambda}$ contains all monomials in 
$\C[x]_{4/3}$ which contain $x_2$.
For $g\in\{x_0^2,x_0x_1,x_1^2\}$, the intersection
$J^{pol}_{f_\lambda}\cap \C[x]_{4/3}$ contains
$g\cdot \frac{\paa f_\lambda}{\paa x_0}$, so 
$g(-(\lambda+1)x_1^2+2\lambda x_0x_1)$, and
$g\cdot\frac{\paa f_\lambda}{\paa x_1}$, and thus also for 
$i\in\{0,1\}$
\begin{eqnarray*}
&&x_ix_1\frac{\paa f_\lambda}{\paa x_1}+
(\frac{-1}{2}x_ix_0+\frac{3}{4}\frac{\lambda+1}{\lambda}x_ix_1)
(-(\lambda+1)x_1^2+2\lambda x_0x_1) \\
&&= \frac{-3}{4}\frac{(\lambda-1)^2}{\lambda} x_ix_1^3.
\end{eqnarray*}
Therefore $J^{pol}_{f_\lambda}$ contains $x_1^4$,
$x_0x_1^3$, and with 
$g\cdot \frac{\paa f_\lambda}{\paa x_1}$ also
$x_0^2x_1^2$, $x_0^3x_1$, $x_0^4$. This shows
$J^{pol}_{f_\lambda}\supset\C[x]_{4/3}$.
\hfill$\Box$

\begin{remarks}\label{t4.2}
(i) A priori, we do not see a reason why the monomials 
$m_1,...,m_{\mu-1}$
can be chosen such that they and 
$\frac{\paa f_\lambda}{\paa \lambda}$
generate $\C[x]/J^{pol}_{f_\lambda}$ 
for each $\lambda\in\C-\{0,1\}$ simultaneously.
It is a nice fact, but not crucial. 

(ii) Thus the global family $F^{alg}$ is nice, 
but it is not a unique global unfolding 
of the Legendre family in \eqref{4.2}. 
Though the global family $F^{mar}$ is a unique global unfolding of 
its restriction to $t'=0$, which is a family over $\H$ of functions
on $\C^{n+1}$ and which is the pull back by $\lambda:\H\to\C-\{0,1\}$
of the Legendre family in \eqref{4.2}. The family $F^{mar}$ is unique
because $\H$ is contractible.

(iii) For other 1-parameter families of simple elliptic singularities,
sets $m_1,...,m_{\mu-1}$ of monomials 
with the analogous property as in lemma \ref{t4.1} had been chosen 
in \cite[2.1]{Ja86} and in \cite[in 2.1 around formula (2)]{MS16}.
\end{remarks}

The next theorem is the special case of theorem \ref{t3.8}
for the simple and simple elliptic singularities.

\begin{theorem}\label{t4.3}
Consider for a simple or a simple elliptic singularity
the global family of functions $F^{mar}_t$ above the space $M^{mar}$
in \eqref{4.9} respectively \eqref{4.8}.
For any $t\in M^{mar}$, the global Milnor number
$\mu_{global}(F^{mar}_t):=\sum_{x\in \textup{Crit}(F^{mar}_t)}\mu(F^{mar}_t,x)$ is $\mu$. 
The manifold $M^{mar}$ is an F-manifold with Euler field $E$.
It is a thickening with the properties in theorem \ref{t3.8}
of the moduli space $M^{mar}_\mu(f^{(0)})$.
Here $f^{(0)}=f$ in the case of the simple singularities,
and we choose $f^{(0)}=f_{1/2}$ in the case of the simple elliptic singularities.
The bundle $H_\Z$ is simply the bundle
$H_\Z=\bigcup_{(\tau,t)\in \C\times M^{mar}-\DD^{mar}}H_n((F^{mar}_t)^{-1}(\tau),\Z)$.
\end{theorem}

{\bf Proof:}
In all cases, $f$ respectively $f_\lambda$ 
is a quasihomogeneous polynomial of weighted degree 1
with respect to a weight system ${\bf w}=(w_0,...,w_n)\in\Q\cap (0,\frac{1}{2}]$.
In the case of a simple singularity, all monomials $m_1,...,m_\mu$ have
weighted degree in $\Q\cap(0,1)$. In the case of a simple elliptic
singularity, all monomials $m_1,...,m_{\mu-1}$ have 
weighted degree in $\Q\cap (0,1)$. 
Therefore in all cases, the highest part (with respect to the weighted degree)
of $\frac{\paa F^{mar}_t}{\paa x_i}$ is equal to 
$\frac{\paa f}{\paa x_i}$ respectively 
$\frac{\paa f_{\lambda(t_\mu)}}{\paa x_i}$.
Denote $\C[x]_{\leq q}:=
\{f\in\C[x]\, |\, \deg_{\bf w}f\leq q\}$ for $q\in \Q_{\geq 0}$.
In the case of the simple elliptic singularities,
\eqref{4.10} implies for $q\in\Q_{\geq 0}$
\begin{eqnarray}\label{4.11}
(J_{F^{mar}_t}^{pol})\cap \C[x]_{\leq q}
&+&\sum_{j:\, \deg_{\bf w}m_j\leq q}\C\cdot m_j\\
&+&\left\{\begin{array}{ll}0&\textup{if }q< 1\\ 
\C\cdot\frac{\paa f_\lambda}{\paa \lambda} &
\textup{if }q\geq 1\end{array}\right\} =\C[x]_{\leq q}.\nonumber
\end{eqnarray}
And in the case of the simple singularities,
\begin{eqnarray}\label{4.12}
(J_{F^{mar}_t}^{pol})\cap \C[x]_{\leq q}
+\sum_{j:\, \deg_{\bf w}m_j\leq q}\C\cdot m_j
 =\C[x]_{\leq q}.
\end{eqnarray}
holds. This shows
\begin{eqnarray}\label{4.13}
\dim \C[x]/J_{F^{mar}_t}^{pol}=\mu.
\end{eqnarray}
The left hand side is the global Milnor number $\mu_{global}(F_t^{mar})$.

Observe 
\begin{eqnarray}\label{4.14}
\frac{\paa F^{mar}_t}{\paa t_j}=\left\{\begin{array}{ll}
m_j &\textup{for the ADE singularities},\\
m_j&\textup{for the }\www E_k\textup{ cases with }j\leq \mu-1,\\
\frac{\paa f_\lambda}{\paa\lambda}
\cdot\frac{\paa\lambda}{\paa t_\mu} &
\textup{for the }\www E_k\textup{ cases with }j=\mu. \end{array}\right.
\end{eqnarray}
This and \eqref{4.11} and \eqref{4.12} show that here the algebraic variant 
\begin{eqnarray}\label{4.15}
{\bf a}_C^{alg,t}:T_{M^{mar},t}\to \C[x]/J^{alg}_{F^{mar}_t},\quad 
\frac{\paa }{\paa t_j}\mapsto [\frac{\paa F^{mar}_t}{\paa t_j}],
\end{eqnarray}
(which is here for simplicity written pointwise) of the Kodaira-Spencer
map in \eqref{2.4} is an isomorphism and equips $T_{M^{mar},t}$
with a multiplication, a unit field vector and an Euler field vector.
This gives $M^{mar}(f^{(0)})$ the structure of an F-manifold with Euler field.
The proof of \cite[Theorem 5.3]{He02} works also here.
The unit field $e$ and the Euler field $E$ are here
\begin{eqnarray}\label{4.16}
e=\frac{\paa}{\paa t_1},\quad 
E=\sum_{j=1}^{\mu} (1-\deg_{\bf w}m_j)t_j\frac{\paa}{\paa t_j}.
\end{eqnarray}

In the cases of the simple singularities, the weights $1-w_j$ are all
positive. The $F$-manifold structure on $M^{mar}$ can be obtained from that
of the germ $(M^{mar},0)$ by using the flow of the Euler field
and $\Lie_E(\circ)=\circ$. 

In the cases of the simple elliptic singularities, $1-w_\mu=0$,
but all other weights $1-w_j$ are positive. The $F$-manifold structure
on $\C^{\mu-1}\times U$ for $U\subset\H$ a neighborhood of a point
$t_\mu$ can be obtained from that on the germ $(M^{mar},(0,t_\mu))$
by using the flow of the Euler field and $\Lie_E(\circ)=\circ$.

A polynomial $g\in\C[x_0,...,x_n]=\C[x]$ is {\it tame}
in the sense of Broughton \cite[Definition 3.1]{Br88} if
a compact neighborhood $U\subset\C^{n+1}$ exists such that
$||(\frac{\paa g}{\paa x_0},...,\frac{\paa g}{\paa x_n})||$
is bounded away from 0 on $\C^{n+1}-U$. 
He proved that $g$ is tame if and only of 
$\mu_{global}(g)=\mu_{global}(g+\sum_{i=0}^n x_is_i)$ for any
$(s_0,...,s_n)\in\C^{n+1}$ \cite[Proposition 3.1]{Br88}.
His main result is that for a tame polynomial $g$, the fiber 
$g^{-1}(\tau)$ for an arbitrary $\tau\in\C$ has the homotopy type of 
$\mu_{global}(g)
-\sum_{x\in \textup{Crit}(g^{-1}(\tau))}\mu(g,x)$ 
many $n$-spheres \cite[Theorem 1.2]{Br88}.

This applies to all the polynomials $F^{mar}_t$.
Because $\mu_{global}(F^{mar}_t)=\mu$ holds for all of them,
and because the unfolding $F^{mar}$ comprises the unfolding
by the terms $+\sum_{i=0}^n x_is_i$, they are all tame.
The eigenvalues of $E\circ:T_tM\to T_tM$ are the critical values of $F_t$. 
Therefore a fiber $(F^{mar}_t)^{-1}(\tau)$ is smooth
if and only if $(\tau,t)\in\C\times M^{mar}-\DD^{mar}$. 
By Broughton such a fiber has the homotopy type of a bouquet of
$\mu$ $n$-spheres. Therefore the middle homology groups
of these fibers glue to a flat $\Z$-lattice bundle of rank $\mu$,
\begin{eqnarray}\label{4.17}
H_\Z&:=& \bigcup_{(\tau,t)\in \C\times M^{mar}-\DD^{mar}}
H_n((F^{mar}_t)^{-1}(\tau),\Z)
\end{eqnarray}

Many of the properties of $M^{mar}(f^{(0)})$ and $H_\Z$ in theorem
\ref{t3.8} are now clear: (a)(i)+(ii) and (b)(i) are obvious.
(b)(ii) holds because $M^{mar}$ is simple connected.
In the case of the simple singularities, (b)(iii) is empty,
as $M^{mar}_\mu$ consists of a single point.
In the case of the simple elliptic singularities,
(b)(iii) holds by the proof in \cite[Theorem 6.1]{GH17-1}
that $M^{mar}_\mu$ is isomorphic to $\H$. There the markings on the 
points in $\H$ were defined essentially by the commutativity of
the diagram \eqref{3.17}.

It rests to prove (a)(iii) and (b)(iv).
First we treat the simple singularities, where this is easier.
As $M^{mar}_\mu(f)$ consists of only one point, 
the stabilizer of this point in $G_\Z(f)$ is the whole group $G_\Z(f)$, 
so \eqref{3.15} becomes
\begin{eqnarray}\label{4.18}
G_\Z(f)=G^{mar}_\RR(f).
\end{eqnarray} 
The homomorphism in theorem \ref{t3.1} (c)
becomes a natural surjective homomorphism $G_\Z(f)\to\Aut_M$ 
with kernel $\{\pm\id\}$. 
Because of the positive $\C^*$-action
by the flow of the Euler field on $M^{mar}$, 
\begin{eqnarray}\label{4.19}
\Aut_M=\Aut(M^{mar},\circ,e,E).
\end{eqnarray}
By lemma \ref{t2.3}, any such automorphism lifts to an up to 
$\pm 1$ unique automorphism of the canonical $\Z$-lattice
bundle $H_\Z$. The group of these automorphisms is $G_\Z(f)$.
This proves (a)(iii) and (b)(iv) in theorem \ref{t3.8}
for the simple singularities.

Now we treat the simple elliptic singularities.
Consider a group element $\chi\in G_\Z(f^{(0)})$, a point
$$(0,t_\mu)=[(f_{\lambda(t_\mu)},\pm\rho)]
\in M^{mar}_\mu(f^{(0)})\subset M^{mar}(f^{(0)}),$$
and its image 
$$\chi_{mar}(t)=[(f_{\lambda(t_\mu)},\pm\chi\circ\rho)]
=(0,\www t_\mu)=[(f_{\lambda(\www t_\mu)},\pm\www\rho)]
\in M^{mar}(f^{(0)})$$
under the action of $\chi_{mar}$ on $M^{mar}_\mu(f^{(0)})$.
Consider the isomorphism of F-manifolds with Euler fields
\begin{eqnarray}\label{4.20}
\psi_{(0,\www t_\mu)}\circ \psi_{(0,t_\mu)}^{-1}:
U_{(0,t_\mu)}\to U_{(0,\www t_\mu)}.
\end{eqnarray}
which (a)(i) in theorem \ref{t3.8} provides.
We claim that these isomorphisms for varying $t_\mu$ 
glue to an automorphism of $M^{mar}(f^{(0)})$
and that this lifts to an automorphism of $H_\Z$
which restricts on the trivial $\Z$-lattice bundle above 
$\bigcup_{t\in M^{mar}(f^{(0)})}\R_{>r(t)}\times\{t\}$
in theorem \ref{t3.8} (b)(ii) to $\chi$.

In several steps one sees that one local isomorphism in 
\eqref{4.20} extends to a global automorphism of $M^{mar}(f^{(0)})$.
First step: Its restriction to $M^{mar}_\mu(f^{(0)})$ is well defined
and given by $\chi_{mar}$.
Second step: The local isomorphism in \eqref{4.20}
extends to an automorphism of a suitable neighborhood of
$M^{mar}_\mu(f^{(0)})$ in $M^{mar}(f^{(0)})$ because 
$M^{mar}_\mu(f^{(0)})=\H$ is contractible.
Third step: 
With the $\C^*$-action by the flow of the Euler field
 on $M^{mar}(f^{(0)})$,
this extends to a global automorphism of $M^{mar}(f^{(0)})$.

Above the extension to $\C\times U_{(0,t_\mu)}\to \C\times U_{(0,\www t_\mu)}$
of the isomorphism in \eqref{4.20}, one has an isomorphism of the
corresponding restrictions of $H_\Z$, because they are isomorphic
to the canonical $\Z$-lattice bundles above 
$U_{(0,t_\mu)}$ and $U_{(0,\www t_\mu)}$ in definition/lemma \ref{t3.2}.
This isomorphism is unique up to $\pm 1$ by definition/lemma \ref{t3.2}.
The commuting diagram \eqref{3.17} tells that the restricted 
isomorphism from $H_\Z$ above 
$\bigcup_{s\in U_{(0,t_\mu)}}\R_{>r(s)}\times\{s\}$
to $H_\Z$ above 
$\bigcup_{s\in U_{(0,\www t_\mu)}}\R_{>r(s)}\times\{s\}$
is compatible with $\pm\chi$. 

Because of the uniqueness up to $\pm 1$, 
all the local isomorphisms of restrictions of $H_\Z$
to neighborhoods of points $(0,0,t_\mu)\in\C\times M^{mar}(f^{(0)})$
glue (possibly after changing some by $\pm\id$) to one automorphism of $H_\Z$.
Its restriction to 
$\bigcup_{s\in M^{mar}}\R_{>r(s)}\times\{s\}$ is $\chi$.

Only now it becomes clear that the automorphism of $M^{mar}(f^{(0)})$
restricts for {\it any} $s_\mu\in \H$ to the isomorphism in \eqref{4.20}:
Its restriction to $U_{(0,s_\mu)}$ and the automorphism in \eqref{4.20}
are in the same way compatible with $\pm\chi$, therefore they coincide
if $U_{(0,s_\mu)}$ and $U_{(0,t_\mu)}$ overlap.
Now (a)(iii) and (b)(iv) in theorem \ref{3.10} are proved for the
simple elliptic singularities. \hfill$\Box$

%The diagram \eqref{3.17} for $t=(0,t_\mu)$ and for $\www t=(0,\www t_\mu)$ gives the 
%commuting diagram 
%\begin{eqnarray}\label{4.21}
%\begin{xy}
%\xymatrix{ H_{\Z,(\tau,t)} \ar[r]^{\cong}_{(b)(i)}  \ar[d]_{\cong}^{(b)(ii)}
%& Ml(f_{\lambda(t_\mu)})  \ar[d]_{\cong}^{\pm\www\rho^{-1}\circ\rho}\\
%H_{\Z,(\tau,\www t)} \ar[r]^{\cong}_{(b)(i)} & Ml(f_{\lambda(\www t_\mu)})  }
%\end{xy}
%\end{eqnarray}

\bigskip
In the next section, we will be more concrete about the action of
$G_\Z(f^{(0)})$ on $M^{mar}(f^{(0)})$.

\section{Symmetries of the simple and the simple elliptic singularities}\label{s5}
\setcounter{equation}{0}

\noindent
In this section, we will write down concrete formulas for
the action on $M^{mar}(f^{(0)})$ of generating elements of $G_\Z(f^{(0)})$.
We need these formulas for an explicit calculation of certain numbers
in section \ref{s10}.
The formulas will also reprove a part of (a)(iii) and (b)(iv) in theorem \ref{t3.8}
for the simple and the simple elliptic singularities. 
But we prefer to keep the entity of the conceptual arguments in the last
part of the proof of theorem \ref{t4.3}, than to drop some of them
and mix the others with the concrete calculations below.

\subsection{Symmetries of the simple singularities}\label{s5.1}
We discussed the symmetries in the proof of theorem \ref{t4.3},
in the paragraph which contains the formulas \eqref{4.18}
and \eqref{4.19}.
In the case of a simple singularity $f$,
$M^{mar}_\mu(f)=\{\textup{pt}\}$, $G_\Z(f)=G^{mar}_{\RR}(f)$,
and
\begin{eqnarray}\label{5.1}
\Aut(M^{mar},\circ,e,E)=\Aut_M\cong G_\Z(f)/\{\pm\id\}.
\end{eqnarray}
By the theorems \ref{t8.3} and \ref{t8.4} in \cite{He11}
\begin{eqnarray}\label{5.2}
G_\Z(f)&=& \{\pm M_h^k\, |\, k\in\Z\}\times U\textup{ with}\\
U&\cong& \left\{
\begin{array}{ll}
S_1 & \textup{for }A_\mu,D_{2l+1},E_6,E_7,E_8,\\
S_2 & \textup{for }D_{2l}\textup{ with }l\geq 3,\\
S_3 & \textup{for }D_4.
\end{array}\right.  \nonumber
\end{eqnarray}
Remark \ref{t3.4} applies to $f$ and $F^{mar}=F^{alg}$
and gives 
\begin{eqnarray*}
()_M\circ ()_{hom}^{-1}(M_h) = 
(t\mapsto (e^{2\pi i \deg_{\bf w}t_0}t_0,...,
e^{2\pi i\deg_{\bf w}t_n}t_n))\in\Aut_M.
\end{eqnarray*}

In all cases with $U=S_1$, this automorphism of
$(M^{mar},\circ,e,E)$ generates $\Aut_M$.

In the cases $D_{2l}$, the coordinate change 
\begin{eqnarray}\label{5.3}
\varphi_2=(x\mapsto (x_0,-x_1,x_2,...,x_n))\in
\textup{Stab}_{G_{\bf w}}(f)\subset\RR^f
\end{eqnarray}
and $\Phi_2=(\id_X,(\varphi)_M)$ satisfy
\begin{eqnarray}\label{5.4}
(\varphi_2)_M&=& (t\mapsto (t_1,-t_2,t_3,...,t_\mu))\\
&\notin& \{()_M\circ ()_{hom}^{-1}(M_h^k)\, |\, k\in\Z\},
\nonumber\\
(\varphi_2)_{hom}&\notin& \{\pm M_h^k\, |\, k\in\Z\},\nonumber\\
\Phi_2\circ\varphi_2 &=& (\varphi_2,(\varphi_2)_M),\nonumber\\
F&=& F\circ (\Phi_2\circ\varphi_2),\quad
\pr_M\circ (\Phi_2\circ\varphi_2)
=(\varphi_2)_M\circ \pr_M.\nonumber
\end{eqnarray}

So, in the cases $D_{2l}$ with $l\geq 3$, 
$(\varphi_2)_{hom}$ can be chosen as a generator of $U$.

In the case $D_4$, $U$ is generated by $(\varphi_2)_{hom}$
and $(\varphi_3)_{hom}$ where
\begin{eqnarray}\label{5.5}
\varphi_3&:=& 
(x\mapsto (\frac{-1}{2}x_0+\frac{-i}{2}x_1,
\frac{3i}{2}x_0+\frac{1}{2}x_1,x_2,...,x_n))\\
&\in &\textup{Stab}_{G_{\bf w}}(f)\subset \RR^f.\nonumber
\end{eqnarray}
This follows from theorem \ref{t3.1}, theorem \ref{t3.3}
and the fact, which can be checked easily, that the group
$\textup{Stab}_{G_{\bf w}}\cong R_f$ is in the case 
$D_4$ with $n=1$ generated by $\varphi_1,\varphi_2$ 
and $\varphi_3$ where $\varphi_1$ is as in remark
\ref{t3.4}.

Only the unfolding morphism $(\Phi_3,(\varphi_3)_M)$
which induces $F\circ \varphi_3^{-1}$ by $F$,
i.e., which satisfies \eqref{3.8},
\begin{eqnarray}
F\circ (\Phi_3\circ\varphi_3)=F,\quad 
\pr_M\circ(\Phi_3\circ\varphi_3) =(\varphi_3)_M\circ \pr_M,
\label{5.6}
\end{eqnarray}
is much more complicated than $(\Phi_1,(\varphi_1)_M)$
and $(\Phi_2,(\varphi_2)_M)$. It is given by 
\begin{eqnarray}\label{5.7}
(\Phi_3(x_0),\Phi_3(x_1))&=& 
(x_0+\frac{-1}{4}t_4,x_1+\frac{i}{4}x_1),\\
(\varphi_3)_M^{-1}(t_1,t_2,t_3,t_4)&=&
(t_1+\frac{i}{4}t_2t_4+\frac{-1}{4}t_3t_4+\frac{1}{16}t_4^3,
\nonumber\\
&&\frac{1}{2}t_2+\frac{-i}{2}t_3+\frac{i}{8}t_4^2,\nonumber\\
&&\frac{3i}{2}t_2+\frac{-1}{2}t_3+\frac{3}{8}t_4^2,\ t_4).
\label{5.8}
\end{eqnarray}
Here one calculates \eqref{5.8} with the ansatz \eqref{5.7} and
\begin{eqnarray}\label{5.9}
F_t((\Phi_3\circ\varphi_3)(x)) &=& F_{(\varphi_3)_M^{-1}}(x).
\end{eqnarray}
For the simple elliptic singularities, we will encounter 
something similar, one coordinate change $\varphi_3$
for which $\Phi_3$ looks difficult.

\begin{remark}\label{t5.1}
For the simple singularities, it is rather obvious
(and it will be shown in the proof of theorem \ref{t7.1})
that $\Aut_M$ is the group of covering transformations
of the covering 
\begin{eqnarray*}
LL^{mar}:M^{mar}-(\KK^{mar}_2\cup\KK^{mar}_3)\to M_{LL}^{(\mu)}
-\DD_{LL}^{(\mu)}
\end{eqnarray*}
in theorem \ref{t6.1}.

This given, the results above (together with the shape of
$\{\pm M_h^k\, |\, k\in\Z\}$, see e.g. the theorems
8.3 and 8.4 in \cite{He11}) prove the main theorem in 
\cite{Li81} which describes this covering group.
This theorem and the isomorphism
$\Aut_M\cong G_\Z/\{\pm\id\}$ have also been 
(re)proved in \cite[Theorem 1 and Theorem 2]{Yu99}.
\end{remark}

\subsection{Symmetries of the simple elliptic singularities}
\label{s5.2}

\noindent
The group $G_\Z=G_\Z(f^{(0)})$ of the simple elliptic
reference singularity $f^{(0)}=f_{1/2}$ (see theorem \ref{t4.3})
sits by theorem 3.1 in \cite{GH17-1} in an exact sequence
\begin{eqnarray}\label{5.10}
1\to (U_1^0\rtimes U_2)\times \{\pm\id\}\to G_\Z\hspace*{2cm}\\
\to \Aut(Ml(f^{(0)})_{(-1)^n,\Z},L)/\{\pm\id\}\to 1 \nonumber
\end{eqnarray}
where 
\begin{eqnarray}\label{5.11}
\Aut(Ml(f^{(0)})_{(-1)^n,\Z},L)&\cong & \textup{SL}(2,\Z)
\end{eqnarray}
and
\begin{eqnarray}
&& U_1^0\cong  
\{(\alpha,\beta,\gamma)\in\Z/p\Z\times \Z/q\Z\times\Z/r\Z\, 
\nonumber \\
&& \hspace*{2cm}
\frac{\alpha}{p}+\frac{\beta}{q}+\frac{\gamma}{r}\equiv
0\mod\Z\}\nonumber \\
&&\begin{array}{c|c|c|c}
 & \www E_6 & \www E_7 & \www E_8 \\
(p,q,r) & (3,3,3) & (4,4,2) & (6,3,2) \\
U_2\cong & S_3 & S_2 & S_1 
\end{array}\label{5.12}
\end{eqnarray}

By theorem 6.1 in \cite{GH17-1}, the action of $G_\Z$ on
$M^{mar}_\mu$ pulls down to an action of the quotient
$\Aut(Ml(f^{(0)})_{(-1)^n,\Z},L)/\{\pm\id\}$ in the 
exact sequence \eqref{5.10} on $M^{mar}_\mu$,
and by the isomorphisms \eqref{5.11} and 
$M^{mar}_\mu\cong\H$ this becomes the standard action
of $\textup{PSL}(2,\Z)$ on $\H$. 

The action of $\textup{PSL}(2,\Z)$ on $\H$ descends
to the action of $S_3$ on $\C-\{0,1\}$ where
$S_3$ acts via
\begin{eqnarray}\label{5.13}
S_3\cong \{\lambda\mapsto \lambda, 1-\lambda,
\frac{1}{\lambda},\frac{\lambda-1}{\lambda},
\frac{\lambda}{\lambda-1},\frac{1}{1-\lambda}).
\end{eqnarray}
This and theorem \ref{t3.6} (c),
$M^{mar}_\mu/G_\Z\cong M_\mu$, reprove the well known
fact that the orbits of this action of $S_3$ on
$\C-\{0,1\}$ give the right equivalence classes of one
family of Legendre normals forms in table \eqref{4.2}.

The kernel $(U_1^0\rtimes U_2)\times\{\pm\id\}$
in the exact sequence \eqref{5.10} acts on the fibers
of the projection
\begin{eqnarray*}
\pr_\mu^{mar}:M^{mar}=\C^{\mu-1}\times\H\to\H,\ 
t\mapsto t_\mu.
\end{eqnarray*}
This action pulls down to an action on the fibers of the
projection
\begin{eqnarray*}
\pr_\mu^{alg}:M^{alg}=\C^{\mu-1}\times (\C-\{0,1\})\to
\C-\{0,1\},\ 
(t',\lambda)\mapsto \lambda.
\end{eqnarray*}

But the action of $G_\Z$ on $M^{mar}$ does not pull down
to an action of a quotient of $G_\Z$ on $M^{alg}$, 
because the covering group
($\cong \Gamma(2)/\{\pm{\bf 1}_2\}\subset\textup{PSL}(2,\Z)$)
of the coverings $M^{mar}\to M^{alg}$
and $\H\to \C-\{0,1\}$ is not a normal subgroup of 
$G_\Z$.

The action of $G_\Z$ on $M^{mar}$ pulls only down to an action
of a {\it groupoid} (see e.g. \cite{ALR07} for the definition
of a groupoid) on $M^{alg}$, whose quotient is
$M^{mar}/G_\Z$. We will not delve into this.
The groupoid structure comes from all isomorphisms
$(M^{alg},(t'^{(1)},\lambda^{(1)}))
\to (M^{alg},(t'^{(2)},\lambda^{(2)}))$ of germs of F-manifolds
(they all underlie isomorphisms of universal unfoldings).

As global maps, many of these isomorphisms become multivalued.
For the use in section \ref{s10}, we make some of them
explicit below. Before, we state a lemma on the relation
between $M^{alg}$ and $M^{mar}/G_\Z$, which will be
used in corollary \ref{t7.3}.

\begin{lemma}\label{t5.2}
The map $M^{alg}\to M^{mar}/G_\Z$ is finite and flat
and has the following degree,
\begin{eqnarray}
\www E_6:&& 6\cdot 2\cdot 3\cdot 3^2=326,\nonumber\\
\www E_7:&& 6\cdot 1\cdot 4\cdot 2^2=96,\label{5.14}\\
\www E_8:&& 6\cdot 1\cdot 6\cdot 1^2=36.\nonumber
\end{eqnarray}
\end{lemma}

{\bf Proof:}
Finiteness and flatness are clear.
The degree is $|S_3|\cdot |U_1^0|\cdot |U_2|$.
By \eqref{5.12}, this is the number in \eqref{5.14}.
\hfill$\Box$

\bigskip
In section \ref{s10}, we need to compare neighborhoods in
$M^{alg}=\C^{\mu-1}\times(\C-\{0,1\})$ of
$\C^{\mu-1}\times\{0\}$, $\C^{\mu-1}\times\{1\}$ and 
$\C^{\mu-1}\times\{\infty\}$.
For this, we give now multivalued maps 
$\psi_2,\psi_3:M^{alg}\dashrightarrow M^{alg}$ which
underlie locally isomorphisms of unfoldings and which lift
the automorphisms $\lambda\mapsto\frac{1}{\lambda}$
and $\lambda\mapsto 1-\lambda$ of $\C-\{0,1\}$.
In each case $\www E_k$, $k\in\{6,7,8\}$, we will give 
two coordinate changes $\varphi_2$ and $\varphi_3$
and multivalued maps
\begin{eqnarray*}
\Psi_2,\Psi_3&:&\C^{n+1}\times M^{alg}\dashrightarrow
\C^{n+1},\\
\psi_2,\psi_3&:& M^{alg}\dashrightarrow M^{alg}
\end{eqnarray*}
with 

\begin{eqnarray}\label{5.15}
F^{alg}_{t',\lambda}((\Psi_i\circ\varphi_i)(x))
&=& F^{alg}_{\psi_i(t',\lambda)}(x),\\
\pr_M^{alg}\circ\psi_2 &=& 
(\lambda\mapsto\frac{1}{\lambda})\circ \pr_M^{alg},\label{5.16}\\
\pr_M^{alg}\circ\psi_3 &=& 
(\lambda\mapsto 1-\lambda)\circ \pr_M^{alg}.
\label{5.17}
\end{eqnarray}
Then $\Phi_i :=(\Psi_i,\psi_i^{-1})$ and $\varphi_i$
satisfy \eqref{3.8}.

The choice of $\phi_2,\Psi_2$ and $\psi_2$ is rather
obvious, and there $f_\lambda\circ \varphi_2=f_{1/\lambda}$,
\eqref{5.15}, \eqref{5.16} and \eqref{3.8}
are easy to check. Also the property of $\varphi_3$,
$$f_\lambda\circ \varphi_3=f_{1-\lambda}$$
is easy to see. But $\Psi_3$ looks more difficult.
It is determined by the requirement that
$F_{t',\lambda}^{alg}(\Psi_3\circ \varphi_3(x)) 
[=F^{alg}(\Psi_3(\varphi_3(x),t',\lambda),t',\lambda)]$
is an unfolding of $f_{1-\lambda}$ only in the monomials
in table \eqref{4.7}.
That means that the coefficients of the following monomials
must vanish:
\begin{eqnarray}
\textup{For }\www E_6:&& x_0x_2,\ 
x_1^2\ (\textup{automatic}), \ x_2^2\ (\textup{automatic}),
\nonumber\\
\textup{For }\www E_7:&& x_0^3,\ 
x_1^3\ (\textup{automatic}),
\label{5.18}\\
\textup{For }\www E_8:&& x_0^5,\ 
x_0^3x_1\ (\textup{automatic}), \ x_0^4\ (\textup{automatic}).
\nonumber
\end{eqnarray}
Having $\Psi_3$ and $\varphi_3$, $\psi_3$ is determined
by \eqref{5.15}. In the case of $\www E_6$ it takes 
two lines, in the case of $\www E_7$ it takes 11
lines, but in the case of $\www E_8$ it would take 3 pages.
There we do not write down $\psi_3$ completely,
we write down only the part of it which is relevant
in section \ref{s10}.

\medskip
{\bf The case $\www E_6$:}
\begin{eqnarray}
\varphi_2(x_0,x_1,x_2)&=& (\lambda^{-1}x_0,x_1,\lambda^{1/2}x_2),
\label{5.19}\\
\Psi_2(x,t',\lambda)&=& x,\label{5.20}\\
\psi_2(t',\lambda)&=& (t_1,\lambda^{-1}t_2,t_3,\lambda^{1/2}t_4,
\nonumber\\
&&\lambda^{-2}t_5,\lambda^{-1}t_6,\lambda^{1/2}t_7,\lambda^{-1}).
\label{5.21}
\end{eqnarray}

\begin{eqnarray}\label{5.22}
\varphi_3(x)&=& (-x_0,x_1-x_0,ix_2),\\ 
\Psi_3(x,t',\lambda) &=& (x_0,x_1,x_2-\frac{i}{2}t_7),
\label{5.23}\\
\psi_3(t',\lambda)&=& (t_1+\frac{1}{2}t_4t_7, 
-t_2-t_3-\frac{1}{4}t_7^2,t_3+\frac{1}{2}t_7^2,\nonumber\\
&& it_4,t_5+t_6,-t_6,it_7,1-\lambda).\label{5.24}
\end{eqnarray}

\medskip
{\bf The case $\www E_7$:}
\begin{eqnarray}\label{5.25}
\varphi_2(x)&=& (\lambda^{-3/4}x_0,\lambda^{1/4}x_1),\\
\Psi_2(x,t',\lambda) &=& x,\label{5.26}\\
\psi_2(t',\lambda)&=& (t_1,\lambda^{-3/4}t_2,\lambda^{1/4}t_3,
\lambda^{-3/2}t_4,\lambda^{-1/2}t_5,\nonumber\\
&& \lambda^{1/2}t_6,\lambda^{-5/4}t_7,\lambda^{-1/4}t_8,
\lambda^{-1})\label{5.27}.
\end{eqnarray}

\begin{eqnarray}\label{5.28}
\varphi_3(x)&=& (-\xi x_0,\xi(x_1-x_0))\quad\textup{ with }
\xi = e^{2\pi i 1/8},\\
\Psi_3(x,t',\lambda) &=& (x_0,x_1-\frac{t_7+t_8}{1-\lambda}),
\label{5.29}\\ 
\psi_3(t',\lambda)&=& (\www t_1,...,\www t_8,1-\lambda)
\quad\textup{with }\label{5.30}\\
\www t_1&=& t_1+(-1)\frac{t_7+t_8}{1-\lambda}t_3+
\left(\frac{t_7+t_8}{1-\lambda}\right)^2 t_6,\nonumber\\
\www t_2&=& (-\xi)t_2+(-\xi)t_3 + \xi \frac{t_7+t_8}{1-\lambda}t_5 
+2\xi \frac{t_7+t_8}{1-\lambda}t_6 \nonumber\\
&&\hspace*{1cm}+(-\xi)\left(\frac{t_7+t_8}{1-\lambda}\right)^2t_8
+\xi\left(\frac{t_7+t_8}{1-\lambda}\right)^3,\nonumber\\
\www t_3&=&  \xi t_3+(-2\xi)\frac{t_7+t_8}{1-\lambda}t_6,\nonumber
\end{eqnarray}

\begin{eqnarray*}
\www t_4&=& \xi^2(t_4+t_5+t_6)
+(-\xi^2)\frac{t_7+t_8}{1-\lambda}t_7
+(-2\xi^2)\frac{t_7+t_8}{1-\lambda}t_8 \nonumber\\
&& \hspace*{1cm}+\xi^2(2-\lambda)\left(\frac{t_7+t_8}{1-\lambda}
    \right)^2, \nonumber\\
\www t_5&=& (-\xi^2)t_5+(-2\xi^2)t_6 
+2\xi^2 \frac{t_7+t_8}{1-\lambda}t_8 
+(-3\xi^2)\left(\frac{t_7+t_8}{1-\lambda}\right)^2,\nonumber\\
\www t_6&=& \xi^2t_6,\nonumber\\
\www t_7&=& \frac{\xi^3}{1-\lambda}((-3+\lambda)t_7
+(-2)t_8),\nonumber\\
\www t_8&=& \frac{\xi^3}{1-\lambda}(3t_7+(2+\lambda)t_8)
\nonumber.
\end{eqnarray*}

\medskip
{\bf The case $\www E_8$:}
\begin{eqnarray}\label{5.31}
\varphi_2(x)&=& (\lambda^{-1/2}x_0,x_1),\\
\Psi_2(x,t',\lambda) &=& x,\label{5.32}\\
\psi_2(t',\lambda)&=& (t_1,\lambda^{-1/2}t_2,\lambda^{-1}t_3,
t_4,\lambda^{-3/2}t_5,\nonumber\\
&&\lambda^{-1/2}t_6,\lambda^{-1}t_7,t_8,
\lambda^{-1/2}t_9,\lambda^{-1}).\label{5.33}
\end{eqnarray}

\begin{eqnarray}\label{5.34}
\varphi_3(x)&=& (ix_0,x_1-x_0^2),\\
\Psi_3(x,t',\lambda) &=& (x_0
+\frac{1}{2}\frac{t_9}{(1-\lambda)^2},
x_1+\frac{t_7+t_8}{1-\lambda}\nonumber\\
&&\hspace*{1cm}
+i\lambda \frac{t_9}{(1-\lambda)^2}x_0
+\frac{1}{4}\frac{t_9^2(4\lambda^2-2\lambda-1)}{(1-\lambda)^4}),
\label{5.35}\\
\psi_3(t',\lambda)&=& (\www t_1,...,\www t_9,1-\lambda)
\quad\textup{with}\label{5.36}\\
\www t_1&=& t_1 + (\textup{a term in }
\C[\lambda,t_2,...,t_9,\frac{t_7+t_8}{1-\lambda},
\frac{t_9}{(1-\lambda)^2}]),\nonumber\\
\www t_2&=& it_2 + (\textup{a term in }
\C[\lambda,t_3,...,t_9,\frac{t_7+t_8}{1-\lambda},
\frac{t_9}{(1-\lambda)^2}]),\nonumber\\
\www t_3&=& -t_3-t_4 + (\textup{a term in }
\C[\lambda,t_5,...,t_9,\frac{t_7+t_8}{1-\lambda},
\frac{t_9}{(1-\lambda)^2}]),\nonumber\\
\www t_4&=& t_4 + (\textup{a term in }
\C[\lambda,t_6,...,t_9,\frac{t_7+t_8}{1-\lambda},
\frac{t_9}{(1-\lambda)^2}]),\nonumber
\end{eqnarray}

\begin{eqnarray}
\www t_5&=& (-i)(t_5+t_6) + (\textup{a term in }
\C[\lambda,t_7,t_8,t_9,\frac{t_7+t_8}{1-\lambda},
\frac{t_9}{(1-\lambda)^2}]),\nonumber\\
\www t_6&=& it_6 + (\textup{a term in }
\C[\lambda,t_7,t_8,t_9,\frac{t_7+t_8}{1-\lambda},
\frac{t_9}{(1-\lambda)^2}]),\nonumber\\
\www t_7&=& \frac{\lambda-3}{\lambda-1}t_7 + 
\frac{-2}{\lambda-1}t_8 + 
\frac{6\lambda+1}{2(1-\lambda)^3}t_9^2,\nonumber\\
\www t_8&=& \frac{3}{\lambda-1}t_7
+\frac{\lambda+2}{\lambda-1}t_8 + 
\frac{14\lambda^2-11\lambda-2}{4(1-\lambda)^4}t_9^2,\nonumber\\
\www t_9&=& i\frac{\lambda^2}{(1-\lambda)^2}t_9.\nonumber
\end{eqnarray}

\begin{remarks}\label{t5.3}
(i) In the case of the minimal number of variables
($n=1$ for $\www E_7$ and $\www E_8$, and $n=2$ for $\www E_6$),
the subgroup $U_1^0\rtimes U_2$ of the kernel
$(U_1^0\rtimes U_2)\times \{\pm\id\}$ in the 
exact sequence in \eqref{5.10} comes via
$\textup{Stab}_{G_{\bf w}}(f_\lambda)\cong R_f
\stackrel{()_{hom}}{\longrightarrow} G_\Z$
from $\textup{Stab}_{G_{\bf w}}(f_\lambda)$ for generic
$\lambda$. This follows from the fact that the kernel of the
exact sequence in \eqref{5.10} is the subgroup of $G_\Z$
which acts trivially on $M^{mar}_\mu\cong\H$.
In the cases of $\www E_7$ and $\www E_8$, the elements
of $\textup{Stab}_{G_{\bf w}}(f_\lambda)$ for generic $\lambda$
can be determined easily explicitly.
In the case of $\www E_6$, this is more difficult.

\medskip
(ii) In any case, one can avoid at the beginning of 
this subsection \ref{s5.2} the use of theorem \ref{t3.1}
in \cite{GH17-1}, which gives the facts in
\eqref{5.10}--\eqref{5.12} on $G_\Z$. 
One can recover these by the following steps:

(1) Determine $\textup{Stab}_{G_{\bf w}}(f_\lambda)$ for 
generic $\lambda$.

(2) Use (i).

(3) Show that 
$\textup{Stab}_{G_{\bf w}}(f_\lambda)$ for generic $\lambda$
and $\varphi_2$ and $\varphi_3$ generate all quasihomogeneous
coordinate changes which map each $f_\lambda$
to some $f_{\www\lambda}$.
\end{remarks}

\section{Lyashko-Looijenga maps for the simple and the 
simple elliptic singularities}\label{s6}
\setcounter{equation}{0}

\subsection{Lyashko-Looijenga maps and their degrees}\label{s6.1}
Lyashko-Looijenga maps in general were discussed in subsection
\ref{s2.4}.
Here we consider the Lyashko-Looijenga maps for the families of functions
defined in section \ref{s4}, the maps 
\begin{eqnarray}
LL^{alg}:M^{alg}&\to& M_{LL}^{(\mu)}  \quad \textup{and}\label{6.1} \\
LL^{mar}:M^{mar}&\to& M_{LL}^{(\mu)},\quad \textup{with} 
\label{6.2}\\
t\in M^{mar}&\mapsto& \prod_{j=1}^\mu (y-u_j)
\quad \textup{ with }u_1,...,u_\mu
\textup{ the} \nonumber\\
&&\textup{critical values of }F^{mar}_t
\textup{ (with multiplicities).}\nonumber
\end{eqnarray}
The caustic $\KK_3^{mar}\subset M^{mar}$ 
and the Maxwell stratum $\KK_2^{mar}\subset M^{mar}$ 
had been defined in \eqref{3.18} and \eqref{3.19}.
They are analytic hypersurfaces.
The caustic $\KK_3^{alg}\subset M^{alg}$ 
and the Maxwell stratum $\KK_2^{alg}\subset M^{alg}$ 
are defined analogously.
They are algebraic hypersurfaces as $LL^{alg}$ is even an
algebraic map.

By Looijenga \cite{Lo74} and Lyashko \cite{Ly79}\cite{Ly84},
the map $LL^{alg}$ restricts to a locally biholomorphic map
$LL^{alg}:M^{alg}-(\KK_3^{alg}\cup \KK_2^{alg})
\to M_{LL}^{(\mu)}-\DD_{LL}^{(\mu)}$,  
it maps $\KK_3^{alg}\cup \KK_2^{alg}$ to $\DD_{LL}^{(\mu)}$,
and it is a branched covering of order 3 respectively 2 
at generic points of $\KK_3^{alg}$ respectively $\KK_2^{alg}$,
and analogous statements hold 
for $LL^{mar}$ and $\KK_3^{mar},\KK_2^{mar}\subset M^{mar}$.

In the case of the simple and the simple elliptic singularities,
we have the following more precise results.

Theorem \ref{t6.1} concerns the simple singularities and was
proved by Looijenga \cite{Lo74} 
and Lyashko \cite{Ly79}\cite{Ly84}.  

The covering result in theorem \ref{t6.2} 
for the simple elliptic singularities 
is an achievement of Jaworski \cite[Theorem 2]{Ja86}\cite[Proposition 1]{Ja88}.

The refinement in theorem \ref{t6.3} of theorem \ref{t6.2}
is a major result of this paper. 
It will be proved in section \ref{s10},
which builds on the sections \ref{s5}, \ref{s8} and \ref{s9}.
It reproves Jaworski's result. 
But the main point is the degree $\deg LL^{alg}$
for the simple elliptic singularities,
which was not calculated before.

Though for the bijections in the main result theorem \ref{t7.1},
we do not need the degree $\deg LL^{alg}$.
Theorem \ref{t6.2} and the analogous part of
theorem \ref{t6.1} are sufficient.

\begin{theorem}\label{t6.1} \cite{Lo74}\cite{Ly79}\cite{Ly84}
In the case of the simple singularities, 
$LL^{alg}$ is a branched covering of degree 
\begin{eqnarray}\label{6.3}
\deg LL^{alg} = \frac{\mu!}{\prod_{j=1}^\mu \deg_{\bf w} t_j}.
\end{eqnarray}
Here $\deg_{\bf w}t_j:=1-\deg_{\bf w}m_j$.
The degree $\deg LL^{alg}$ 
is given explicitly in table \eqref{6.4}.
\begin{eqnarray}\label{6.4}
\begin{array}{l|l|l|l|l|l}
\textup{name} & A_\mu & D_\mu & E_6 & E_7 & E_8 \\
\deg LL^{alg} & (\mu+1)^{\mu-1} & 2(\mu-1)^\mu & 2^9\cdot 3^4 
& 2\cdot 3^{12} & 2\cdot 3^5\cdot 5^7
\end{array}
\end{eqnarray}
And the restriction 
$LL^{alg}:M^{alg}-(\KK_3^{alg}\cup \KK_2^{alg})
\to M_{LL}^{(\mu)}-\DD_{LL}^{(\mu)}$
is a covering.
\end{theorem}

\begin{theorem}\label{t6.2}\cite{Ja86}\cite{Ja88}
 In the case of the simple elliptic singularities, 
the restriction $LL^{alg}:M^{alg}-(\KK_3^{alg}\cup \KK_2^{alg})
\to M_{LL}^{(\mu)}-\DD_{LL}^{(\mu)}$ is a covering.
\end{theorem}

\begin{theorem}\label{t6.3}
In the case of the simple elliptic singularities, an extension $M^{orb}$
\begin{eqnarray}\label{6.5}
\begin{array}{ccccc}(t',\lambda)& \in & M^{alg} & \subset & M^{orb} \\
\downarrow  & & \downarrow & & \downarrow \pi_{orb} \\
\lambda & \in & \C-\{0,1\} & \subset & \P^1
\end{array}
\end{eqnarray}
of $M^{alg}$ to an orbibundle above $\P^1\supset \C-\{0,1\}$ 
exists such that $LL^{alg}$ extends to a surjective 
holomorphic map $LL^{orb}:M^{orb}\to M_{LL}^{(\mu)}$ 
with the following properties.
The two-dimensional subspace $M^{orb}_0\subset M^{orb}$ with 
\begin{eqnarray*}
M^{orb}_0&=&(\textup{closure in }M^{orb}\textup{ of }
\{(t',\lambda)\in M^{alg}\, |\, t_2=...=t_{\mu-1}=0\})\\
&\cong&\C\times\P^1 
\end{eqnarray*}
(which is the $\mu$-constant stratum and its translates
under the unit field)
is mapped to the one-dimensional subspace $M_{LL,0}^{(\mu)}
\subset M_{LL}^{(\mu)}$ with 
\begin{eqnarray*}
M_{LL,0}^{(\mu)}:= \{p(y)\in M_{LL}^{(\mu)}\, |\, 
p(y)=(y-u_1)^\mu ,\, u_1\in\C\} \cong\C.
\end{eqnarray*} 
The restriction 
\begin{eqnarray}\label{6.6}
LL^{orb}&:& M^{orb}-M^{orb}_0 \to M_{LL}^{(\mu)}-M_{LL,0}^{(\mu)}
\end{eqnarray}
is a branched covering of degree 
\begin{eqnarray}\label{6.7}
\deg LL^{orb}=\deg LL^{alg} = 
\frac{\mu!\cdot \frac{1}{2}\cdot\sum_{j=2}^{\mu-1}\frac{1}{\deg_{\bf w}t_j}}
{\prod_{j=2}^{\mu-1} \deg_{\bf w} t_j}.
\end{eqnarray}
Here $\deg_{\bf w}t_j:=1-\deg_{\bf w}m_j$. 
The degree $\deg LL^{alg}$ 
is given explicitly in table \eqref{6.8}.  
\begin{eqnarray}\label{6.8}
\begin{array}{l|l|l|l}
\textup{name} & \www E_6 & \www E_7 & \www E_8 \\
\deg LL^{alg} & 2^2\cdot 3^{11}\cdot 5\cdot 7 
& 2^{18}\cdot 3\cdot 5^3\cdot 7 & 2^9\cdot 3^{10}\cdot 7\cdot 101 
\end{array}
\end{eqnarray}
And $LL^{orb}$ maps 
$\KK_3^{alg}\cup\KK_2^{alg}\cup\pi_{orb}^{-1}(\{0,1,\infty\})$
to $\DD_{LL}^{(\mu)}$.
\end{theorem}

\begin{remark}\label{t6.4}
(i) Let $N_{Coxeter}$ be the Coxeter number of an ADE root lattice,
and $W$ its Weyl group.
By \cite{Bo68} $|W|=N_{Coxeter}^\mu\cdot \prod_{j=1}^\mu\deg_{\bf w}t_j$. 
Therefore
\begin{eqnarray}\label{6.9}
\deg LL^{alg} = \frac{\mu!}{\prod_{j=1}^\mu \deg_{\bf w} t_j}
=\frac{\mu!\cdot N_{Coxeter}^\mu}{|W|}.
\end{eqnarray}
This was observed for example in \cite{Yu90}.

\medskip
(ii) In order to make the tables \eqref{6.4} and \eqref{6.8} transparent,
here we give the weights $\deg_{\bf w}x_i=w_i$, the weights 
$\deg_{\bf w}t_j$, 
in the ADE cases the Coxeter numbers $N_{Coxeter}$,
and in the simple elliptic cases the number
$\frac{1}{2}\sum_{j=2}^{\mu-1}\frac{1}{\deg_{\bf w}t_j}$,
\begin{eqnarray*}
\begin{array}{l|l|l|l|l|l|l|l|l|l|l|l}
 & N_{Coxeter} & x_0 & x_1 & t_1 & t_2 & t_3 & t_4 & ... & 
 t_\mu \\
A_\mu & \mu+1 & \frac{1}{\mu+1} & & 1 & \frac{\mu}{\mu+1} & \frac{\mu-1}{\mu+1} 
& \frac{\mu-2}{\mu+1} & ... & \frac{2}{\mu+1} \\
D_\mu & 2(\mu-1) & \frac{1}{\mu-1} & \frac{\mu-2}{2(\mu-1)} & 1 & \frac{\mu}{2(\mu-1)}
& \frac{\mu-2}{\mu-1} & \frac{\mu-3}{\mu-1} & ... & \frac{1}{\mu-1} 
\end{array}
\end{eqnarray*}
\begin{eqnarray*}
\begin{array}{l|l|l|l|l|l|l|l|l|l|l|l}
 & N_{Coxeter} & x_0 & x_1 & t_1 & t_2 & t_3 & t_4 & t_5 & t_6 & t_7 & t_8 \\
E_6 & 12 & \frac{1}{4} & \frac{1}{3} & 1 & \frac{3}{4} 
& \frac{2}{3} & \frac{1}{2} & \frac{5}{12} & \frac{1}{6} & & 
\\[1mm] 
E_7 & 18 & \frac{2}{9} & \frac{1}{3} & 1 & \frac{7}{9} 
& \frac{2}{3} & \frac{5}{9} & \frac{4}{9} & \frac{1}{3} 
& \frac{1}{9} & \\[1mm]
E_8 & 30 & \frac{1}{5} & \frac{1}{3} & 1 & \frac{4}{5} 
& \frac{2}{3} & \frac{3}{5} & \frac{7}{15} & \frac{2}{5} 
& \frac{4}{15} & \frac{1}{15} 
\end{array}
\end{eqnarray*}
\begin{eqnarray*}
\begin{array}{l|l|l|l|l|l|l|l|l|l|l|l|l|l}
 & x_0 & x_1 & x_2 & t_1 & t_2 & t_3 & t_4 & t_5 & t_6 & t_7 & t_8 & t_9 
 & \frac{1}{2}\sum_{j=2}^{\mu-1}\frac{1}{\deg_{\bf w}t_j}\\
\www E_6 & \frac{1}{3} & \frac{1}{3} & \frac{1}{3} & 1 & \frac{2}{3} 
 & \frac{2}{3} & \frac{2}{3} & \frac{1}{3} & \frac{1}{3} & \frac{1}{3} & & 
 & \frac{27}{4} \\ 
\www E_7 & \frac{1}{4} & \frac{1}{4} & & 1 & \frac{3}{4} & \frac{3}{4} 
 & \frac{1}{2} & \frac{1}{2} & \frac{1}{2} & \frac{1}{4} & \frac{1}{4} &
 & \frac{25}{3} \\ 
\www E_8 & \frac{1}{6} & \frac{1}{3} & & 1 & \frac{5}{6} & \frac{2}{3} 
 & \frac{2}{3} & \frac{1}{2} & \frac{1}{2} & \frac{1}{3} & \frac{1}{3} 
 & \frac{1}{6} & \frac{101}{10}
\end{array}
\end{eqnarray*}
\end{remark}

\bigskip

\subsection{Limit behaviour of $LL^{alg}$ after Jaworski}\label{s6.2}
Jaworski's proof of theorem \ref{t6.2} 
required an understanding of the 
limit behaviour of $LL^{alg}$ near $\lambda\in\{0,1,\infty\}$.
Here we will explain a result of him which concerns this limit behaviour.
It will be crucial for the proof of the main theorem \ref{t7.1}
in the case of the simple elliptic singularities.

The intersection form $I$ on the Milnor lattice $Ml(f)$ 
of a simple elliptic singularity is positive semidefinite if $n\equiv 0(4)$, 
see e.g. \cite{AGV88}. 
Consider the Stokes matrix $S$ of any distinguished basis as defined in \eqref{2.16}. 
Because of \eqref{2.18}, $S+S^t$ is positive semidefinite.
Therefore all entries of $S$ are in $\{0,\pm 1,\pm 2\}$. 
Therefore for any two vanishing cycles $\delta_i$ and $\delta_j$
and any $n$ (not necessarily $n\equiv 0(4)$) 
$I(\delta_i,\delta_j)\in\{0,\pm 1,\pm 2\}$.

$LL^{alg}:M^{alg}-(\KK_3^{alg}\cup \KK_2^{alg})\to M_{LL}^{(\mu)}-\DD_{LL}^{(\mu)}$ 
is a covering by theorem \ref{t6.2}.
Now consider a ($C^\infty$ or real analytic) path 
\begin{eqnarray}\label{6.10}
r:[0,\varepsilon)\to M_{LL}^{(\mu)}
&\textup{with}& r((0,\varepsilon))\subset M_{LL}^{(\mu)}-\DD_{LL}^{(\mu)},\\
&\textup{and}& r(0)\in\DD_{LL}^{(\mu),reg}.\nonumber 
\end{eqnarray}
Consider any lift 
$\rho:(0,\varepsilon)\to M^{alg}-(\KK_3^{alg}\cup \KK_2^{alg})$
of the restriction of the path $r$ to $(0,\varepsilon)$.
For $s\in (0,\varepsilon)$, denote by $u_1(s),...,u_\mu(s)$ the
critical values of $F_{\rho(s)}$.
They are pairwise different. 
Because of $r(0)\in \DD_{LL}^{(\mu),reg}$,
precisely two of them will tend to one another if $s\to 0$. 
We can suppose that they are numbered $u_i(s)$ and $u_{i+1}(s)$.
Write 
$\rho(s)=(t_1^{(\rho)}(s),...,t_{\mu-1}^{(\rho)}(s),
\lambda^{(\rho)}(s))$ 
for $s\in(0,\varepsilon)$. 

For fixed $s\in (0,\varepsilon)$, 
consider the $\Z$-lattice bundle 
$\bigcup_{\tau\in\C-\{u_1(s),...,u_\mu(s)\}}
H_n(F_{\rho(s)}^{-1}(\tau),\Z)$.
Move the vanishing cycles at $u_i(s)$ and $u_{i+1}(s)$ 
along straight
lines to the $\Z$-lattice 
$H_n(F_{\rho(s)}^{-1}(\frac{u_i(s)+u_{i+1}(s)}{2}),\Z)$
and call the images $\delta_i(s),\delta_{i+1}(s)$. 
They are unique up to the sign. Then the following holds.

\begin{theorem}\cite[Proposition 2]{Ja88}\label{t6.5}
\begin{eqnarray}
I(\delta_i(s),\delta_{i+1}(s))=0 &\iff& \rho\textup{ extends to 0 and }
\rho(0)\in\KK_2^{alg,reg},\label{6.11} \\
I(\delta_i(s),\delta_{i+1}(s))=\pm 1&\iff& \rho\textup{ extends to 0 and }
\rho(0)\in\KK_3^{alg,reg},\nonumber\\
I(\delta_i(s),\delta_{i+1}(s))=\pm 2&\iff& 
\lambda^{(\rho)}\textup{ extends to 0 and }\lambda^{(\rho)}(0)\in\{0,1,\infty\}.
\nonumber
\end{eqnarray}
\end{theorem}

{\bf Proof:}
The statement in \cite[Proposition 2]{Ja88} is slightly weaker.
Therefore here we provide the additional arguments.
Because of theorem \ref{t6.3}, $\lambda^{(\rho)}$ extends in any case to 0,
and if $\lambda^{(\rho)}(0)\notin\{0,1,\infty\}$, then
$\rho$ extends to $0$ and $\rho(0)\in\KK_2^{alg,reg}\cup \KK_3^{alg,reg}$.
Therefore exactly one of the three cases on the right hand side
of \eqref{6.11} holds. In the third case, Proposition 2 in \cite{Ja88}
applies and gives $I(\delta_i(s),\delta_{i+1}(s))=\pm 2$.

In the second case, $F_{\rho(0)}$ has $\mu-2$ $A_1$ singularities
and one $A_2$-singularity, with pairwise different critical values,
and then $\delta_i(s)$ and $\delta_{i+1}(s)$
are a distinguished basis of the Milnor lattice of the $A_2$ singularity.
Then $I(\delta_i(s),\delta_{i+1}(s))=\pm 1$.

In the first case, $F_{\rho(0)}$ has $\mu$ $A_1$ singularities,
and two of them have the same critical value, the others have pairwise different
critical values. Then $\delta_i(s)$ and $\delta_{i+1}(s)$ are 
vanishing cycles of the two $A_1$ singularities with the same critical value.
Then $I(\delta_i(s),\delta_{i+1}(s))=0$.
\hfill $\Box$ 

\bigskip

The situation for the simple singularities is analogous, but simpler. 
The Milnor lattice $Ml(f)$ with intersection form of an ADE singularity is
isomorphic to the ADE root lattice if $n\equiv 0(4)$, see e.g.
\cite{AGV88}. Consider the Stokes matrix $S$ of any distinguished basis
as defined in \eqref{2.16}. Because of \eqref{2.18}, $S+S^t$ is positive definite.
Therefore all entries of $S$ are in $\{0,\pm 1\}$. 
Therefore for any two vanishing cycles $\delta_i$ and $\delta_j$
and any $n$ (not necessarily $n\equiv 0(4)$) 
$I(\delta_i,\delta_j)\in\{0,\pm 1\}$.

Now consider a ($C^\infty$ or real analytic) path 
$r:[0,\varepsilon)\to M_{LL}^{(\mu)}$ as in \eqref{6.10},
and consider any lift 
$\rho:(0,\varepsilon)\to M^{alg}-(\KK_3^{alg}\cup \KK_2^{alg})$
of the restriction of the path $r$ to $(0,\varepsilon)$.
Because $LL^{alg}:M^{alg}\to M_{LL}^{alg}$ is a branched covering,
$\rho$ extends to $0$, and $\rho(0)\in \KK_2^{alg,reg}\cup \KK_3^{alg,reg}$.

For $s\in (0,\varepsilon)$, denote again by $u_1(s),...,u_\mu(s)$ the
critical values of $F_{\rho(s)}$.
They behave for $s\to 0$ as above, and we obtain vanishing
cycles $\delta_i(s)$ and $\delta_{i+1}(s)$ as above.
Then the following holds.

\begin{lemma}\label{t6.6}
\begin{eqnarray}
I(\delta_i(s),\delta_{i+1}(s))=0 &\iff& 
\rho(0)\in\KK_2^{alg,reg},\label{6.12} \\
I(\delta_i(s),\delta_{i+1}(s))=\pm 1&\iff& 
\rho(0)\in\KK_3^{alg,reg}.\nonumber
\end{eqnarray}
\end{lemma}

The proof is a subset of the proof of theorem \ref{t6.5}.

\section{The main theorem, its proof and consequences}\label{s7}
\setcounter{equation}{0}

\noindent
In subsection \ref{s3.4} we introduced for any reference singularity
$f^{(0)}$ a {\it Looijenga-Deligne map} 
\begin{eqnarray}\label{7.1}
LD:R_{Stokes}\to\BB^{ext}(f^{(0)})/G_{sign,\mu}.
\end{eqnarray}
Recall that $R_{Stokes}$ is the set of {\it Stokes regions},
which are the components of the complement of the Stokes walls
$W_{Stokes}$ in $M^{mar}(f^{(0)})$,
and $\BB^{ext}(f^{(0)})$ is the orbit under $G_\Z$ of the set 
$\BB(f^{(0)})$ of distinguished bases of $f^{(0)}$.
The map $LD$ is $G_\Z$ equivariant.

In the case of a simple singularity $f^{(0)}=f$ or
of a simple elliptic singularity $f^{(0)}=f_{1/2}$
(with $f_\lambda$ the Legendre normal form from \eqref{4.2}),
$M^{mar}(f^{(0)})$ had been constructed in section \ref{s4}.

The main theorem is as follows. 
For the simple singularities, the bijection \eqref{7.3} 
was proved in a different way in \cite{Lo74} and \cite{De74}, 
see remark \ref{t7.2} (iv) below. Yu \cite[6.3 Satz]{Yu90}
built on this and proved the bijection \eqref{7.4} for the
simple singularities.

\begin{theorem}\label{t7.1}
Consider a simple singularity $f^{(0)}=f$ or
a simple elliptic singularity $f^{(0)}=f_{1/2}$. Then
\begin{eqnarray}\label{7.2}
R_{Stokes}=R^{0}_{Stokes},\quad \BB^{ext}(f^{(0)})=\BB(f^{(0)}).
\end{eqnarray}
The Looijenga-Deligne map
\begin{eqnarray}\label{7.3}
LD:R^{0}_{Stokes}\to\BB(f^{(0)})/G_{sign,\mu}
\end{eqnarray}
and the induced quotient map
\begin{eqnarray}\label{7.4}
LD/G_\Z:R^{0}_{Stokes}/G_\Z\to 
\{\textup{Stokes matrices}\}/G_{sign,\mu}
\end{eqnarray}
are bijections.
\end{theorem}

{\bf Proof:} In \cite{GH17-1} it was proved that the moduli space
$M^{mar}(f^{(0)})$ is connected (see remark \ref{t3.7} (ii)).
Therefore $R_{Stokes}=R_{Stokes}^0$.
Recall the argument in remark \ref{t3.11} (ii) for $\BB^{ext}(f^{(0)})=\BB(f^{(0)})$:
The map $LD$ is $G_\Z$ equivariant, and remark \ref{t3.11} (i) shows
that $R_{Stokes}^0$ is mapped to $\BB(f^{(0)})$. Therefore
$\BB^{ext}(f^{(0)})=\BB(f^{(0)})$.

The Stokes matrix of a distinguished basis $\uuuu{\delta}$
of the Milnor lattice $Ml(f^{(0)})$ was defined in \eqref{2.16}
as $S:=(-1)^{(n+1)(n+2)/2}\cdot 
L(\uuuu{\delta}^t,\uuuu{\delta})^t$.
Obviously the set of Stokes matrices can be identified 
with the quotient $\BB(f^{(0)})/G_\Z$. 

It suffices to prove that the map $LD$ in \eqref{7.3} 
is a bijection. Then the quotient map $LD/G_\Z$ in 
\eqref{7.4} is a bijection, too.

Looijenga's argument \cite{Lo74} 
that $LD$ is surjective for the simple singularities,
works because of Jaworski's theorem \ref{t6.2} 
also for the simple elliptic singularities.
The argument is as follows.
Let $U\in R_{Stokes}^0$ be any Stokes region, 
let $t\in U$, and let $\uuuu{\delta}$ be the up to the action of $G_{sign,\mu}$
unique distinguished basis which is constructed from the morsification
$F_t$ and the good distinguished system of paths in definition \ref{t3.9} (b).
Then $\uuuu{\delta}$ is in $LD(U)$.
Let $\uuuu{\gamma}$ be any distinguished basis. 
It is the image of $\uuuu{\delta}$ under the action 
of a certain braid in $\textup{Br}_\mu$
and possibly a sign change in $G_{sign,\mu}$. 
The braid gives a (homotopy class of a) closed path in 
$M_{LL}^{(\mu)}-D_{LL}^{(\mu)}$.
The path has a unique lift to $M^{mar}-(\KK_3^{mar}\cup\KK_2^{mar})$
which starts at $t\in U$ because the Lyashko-Looijenga map
$LL^{mar}:M^{mar}-(\KK_3^{mar}\cup\KK_2^{mar})\to M_{LL}^{(\mu)}-\DD_{LL}^{(\mu)}$
is a covering by the theorems \ref{t6.1} and \ref{t6.2}.
Let $\www t$ be the endpoint of this lift and let $\www U$ be the Stokes
region which contains $\www t$. Then $\uuuu{\gamma}\in LD(\www U)$.
Therefore $LD$ is surjective.

It rests to prove that $LD$ is injective.
Let $U^{(1)}$ and $U^{(2)}$ be two Stokes regions with 
$LD(U^{(1)})=LD(U^{(2)})$. Because the Lyashko-Looijenga map
$LL^{mar}:M^{mar}-(\KK_3^{mar}\cup\KK_2^{mar})\to M_{LL}^{(\mu)}-\DD_{LL}^{(\mu)}$
is a covering, both Stokes regions are mapped by $LL^{mar}$ bijectively to
the open subset 
\begin{eqnarray}\label{7.5}
\{p(y)\in M_{LL}^{(\mu)}&|& p(y)=\prod_{j=1}^\mu(y-u_j)\\
&&\textup{ with }\Imm(u_i)\neq\Imm(u_j)\textup{ for }
i\neq j\}.\nonumber
\end{eqnarray}
of $M_{LL}^{(\mu)}$. There is a unique isomorphism
$\psi^U:U^{(1)}\to U^{(2)}$ which is compatible with $LL^{mar}$.
Obviously it is an isomorphism of semisimple F-manifolds.

We claim that it extends to an automorphism 
$\psi^{mar}:M^{mar}(f^{(0)})\to M^{mar}(f^{(0)})$.
If this is true then the rest of the proof is an elegant application
of theorem \ref{t4.3}: Then $\psi^{mar}$ comes from an element 
$\psi\in G_\Z(f^{(0)})$ (which is unique up to $\pm 1$).
The element $\psi$ must map $LD(U^{(1)})$ to $LD(U^{(2)})$.
As they coincide by assumption, $\psi=\pm\id$.
Thus $\psi^{mar}=\id$ on $M^{mar}$, thus $U^{(1)}=U^{(2)}$.

It rests to show that $\psi^U$ extends an automorphism 
$\psi^{mar}$ of $M^{mar}(f^{(0)})$. 
Roughly, the reason is that the covering 
$LL^{mar}:M^{mar}-(\KK_3^{mar}\cup\KK_2^{mar})
\to M_{LL}^{(\mu)}-\DD_{LL}^{(\mu)}$
with base point in $U^{(k)}$ is determined by the class 
of Stokes matrices
modulo $G_{sign,\mu}$ which are associated to the distinguished bases
in $LD(U^{(k)})$. As this class coincides for $k=1,2$, 
a deck transformation 
$\psi^{mar}:M^{mar}-(\KK_3^{mar}\cup\KK_2^{mar})
\to M^{mar}-(\KK_3^{mar}\cup\KK_2^{mar})$ exists, which 
extends $\psi^{U}$.
It extends to $\KK_3^{mar}\cup\KK_2^{mar}$ as there $LL^{mar}$
is generically branched of order 3 respectively 2.

More precisely, we can argue as follows.
Let $t^{(k)}\in U^{(k)}$ 
be points with $LL^{mar}(t^{(1)})=LL^{mar}(t^{(2)})$.
Then $\psi^U(t^{(1)})=t^{(2)}$.
Choose a path within 
$M^{mar}(f^{(0)})-(\KK_3^{mar,sing}\cup\KK_2^{mar,sing})$
from $t^{(1)}$ to any point in this space. 
We claim that $\psi^U$ extends from $U^{(1)}$ to a well defined map 
\begin{eqnarray}\label{7.6}
\psi^{U\cup\textup{path}}&:&U^{(1)}\cup(\textup{a neighborhood of this path})\\
&&\to M^{mar}-(\KK_3^{mar,sing}\cup \KK_2^{mar,sing})\nonumber 
\end{eqnarray}
and that this is locally an isomorphism of F-manifolds.
If the path does not meet $\KK_3^{mar}\cup \KK_2^{mar}$, this is obvious.
Now suppose that it meets $\KK_3^{mar,reg}\cup \KK_2^{mar,reg}$.
Let $\rho^{(1)}$ be the restriction of the path to a path 
from $t^{(1)}$
to a point $\www t^{(1)}$ just before the first meeting point 
with $\KK_3^{mar,reg}\cup \KK_2^{mar,reg}$. Then 
\begin{eqnarray*}
&&\psi^{U\cup\rho^{(1)}}: U^{(1)}\cup 
\{\textup{a neighborhood of }\rho^{(1)}\} \\
&\to&  M^{mar}-(\KK_3^{mar,sing}\cup \KK_2^{mar,sing})
\end{eqnarray*}
is well defined. Let
$\rho^{(2)}:=\psi^{U\cup\rho^{(1)}}\circ \rho^{(1)}$ 
be the image of $\rho^{(1)}$ under 
$\psi^{U\cup\rho^{(1)}}$.
Then $\rho^{(2)}$ starts at $t^{(2)}$ and ends at
$\www t^{(2)}:=\psi^{U\cup\rho^{(1)}}(\www t^{(1)})$.
Then $LL(\www t^{(1)})=LL(\www t^{(2)})$. 

Let $\www U^{(1)}$ and $\www U^{(2)}$ be the Stokes regions which contain
$\www t^{(1)}$ and $\www t^{(2)}$. 
Then still $LD(\www U^{(1)})=LD(\www U^{(2)})$,
and also the associated Stokes matrices are equal up to the action of 
$G_{sign,\mu}$. 

Write $LL(\www t^{(1)})=LL(\www t^{(2)})
=\prod_{i=1}^\mu (y-u_i)$ with $(u_1,...,u_\mu)$ in good
ordering (definition \ref{t3.9} (a)).
By lemma \ref{t2.3}, the $\Z$-lattice bundles
$\bigcup_{\tau\in\C-\{u_1,...,u_\mu\}}H_n(F_{\www t^{(k)}})^{-1}(\tau),\Z)$
for $k=1$ and $k=2$ are isomorphic. 

Near $\www t^{(k)}$ the path $\rho^{(k)}$ is a lift of a path
$r$ as in \eqref{6.10}. 
By the construction before theorem \ref{t6.5} and lemma \ref{t6.6},
we obtain vanishing cycles $\delta_i^{(k)}$ and $\delta_{i+1}^{(k)}$
in $H_n(F_{\www t^{(k)}}^{-1}(\frac{u_i+u_{i+1}}{2}),\Z)$.
Because the $\Z$-lattice bundles are isomorphic, 
\begin{eqnarray*}
I(\delta_i^{(1)},\delta_{i+1}^{(1)})=I(\delta_i^{(2)},\delta_{i+1}^{(2)}).
\end{eqnarray*}
By theorem \ref{t6.5} and lemma \ref{t6.6}, 
this is either $0$ or $\pm 1$, 
and the first meeting point of the extension of $\rho^{(1)}$ 
is in $\KK_2^{mar}$ in the first case and in $\KK_3^{mar}$ in the second case,
and $\rho^{(2)}$ extends to $\KK_2^{mar}$ in the first case and to
$\KK_3^{mar}$ in the second case. 
Therefore the isomorphism $\psi^{U\cup\rho^{(1)}}$ extends to a local
isomorphism of F-manifolds beyond this first meeting point.

Therefore \eqref{7.6} holds. As $\KK_3^{mar,sing}\cup\KK_2^{mar,sing}$
has codimension two in $M^{mar}$, the extensions of $\psi^U$ to 
local isomorphisms of F-manifolds along all paths in 
$M^{mar}-\KK_3^{mar,sing}\cup\KK_2^{mar,sing}$
glue to one global automorphism $\psi^{mar}$ of $M^{mar}$.
\hfill$\Box$

\begin{remarks}\label{t7.2}
(i) In the case of a simple singularity, the sets $R^0_{Stokes}$,
$\BB(f)/G_{sign,\mu}$ and $\{\textup{Stokes matrices}\}/G_{sign,\mu}$
are all finite. $|R_{Stokes}^0|=\deg LL^{alg}$ is finite 
because the Lyashko-Looijenga map is algebraic. 
$\BB(f)$ and the quotient sets are finite 
because there are only finitely many vanishing cycles, 
they are the roots of the ADE lattice.

In the case of a simple elliptic singularity, the set $R_{Stokes}^0$ 
is not finite, because the universal covering $\lambda:\H\to\C-\{0,1,\infty\}$
has infinite degree. 
The set $\BB(f_{1/2})$ of distinguished bases is not finite 
because the group $G_\Z$ acts on it because of \eqref{7.2}, and
the group $G_\Z$ is not finite \cite{GH17-1}.
The set of Stokes matrices is finite because
the entries of each Stokes matrix are in $\{0,\pm 1,\pm 2\}$, 
see the beginning of subsection \ref{s6.2}.

Ebeling \cite{Eb18} showed that for all other singularities the set $\BB(f)$ and
the set of Stokes matrices are infinite. Together this gives the
picture in table \eqref{1.1}.
%\begin{eqnarray}\label{7.7}
%\begin{array}{lll}
%f & |\BB(f)| & |\{\textup{Stokes matrices}\}| \\ \hline 
%\textup{simple singularity} & \textup{finite} & \textup{finite} \%\
%\textup{simple elliptic singularity} & \textup{infinite} & %\textup{finite} \\
%\textup{any other singularity} & \textup{infinite} & %\textup{infinite}
%\end{array}
%\end{eqnarray}

\medskip
(ii) Theorem \ref{t7.1} together with the degrees of $LL^{alg}$ 
in the theorems \ref{t6.1} and \ref{t6.3} and in the case of the
simple singularities the number $|G_\Z|$ allows now to calculate
all finite numbers in \eqref{1.1}.
Corollary \ref{t7.3} gives the result.

\medskip
(iii) All numbers in corollary \ref{t7.3} 
except the number of Stokes matrices for $\www E_8$ 
had already been determined. Deligne \cite{De74} 
determined the number $|\BB(f)|$ for the simple singularities,
Yu \cite{Yu90}\cite{Yu96}\cite{Yu99} determined
the number $|\{\textup{Stokes matrices}\}|$ for the simple
singularities. Kluitmann determined the number 
$|\{\textup{Stokes matrices}\}|$ for the simple elliptic singularities
$\www E_6$ \cite{Kl83} and $\www E_7$ \cite{Kl87}.
Deligne and Kluitmann worked directly with the braid group orbits
$\BB(f)$. Their calculations are hard, especially those of Kluitmann.
It is satisfying that theorem \ref{t7.1} together with 
$\deg LL^{alg}$ gives the same numbers
$|\{\textup{Stokes matrices}\}|$ for $\www E_6$ and $\www E_7$, 
and that they allow to find the missing number, 
the number $|\{\textup{Stokes matrices}\}|$ for $\www E_8$.

\medskip
(iv) For the simple singularities, Deligne \cite{De74}
and Looijenga \cite{Lo74} proved the bijection \eqref{7.3}
by comparison of numbers. Looijenga proved that \eqref{7.3} is surjective
(see the proof of theorem \ref{t7.1} for his argument)
and calculated $|R_{Stokes}^0|=\deg LL^{alg}$. 
Deligne calculated $|\BB(f)|/G_{sign,\mu}$ and observed
that it coincides with $|R_{Stokes}^0|$. 
Therefore \eqref{7.3} is a bijection.
But for the simple elliptic singularities, both sides of \eqref{7.3}
are infinite, and this argument does not work.

\medskip
(v) For the simple singularities observe 
\begin{eqnarray}\label{7.8}
|\BB(f)|&=&2^\mu\cdot |\BB(f)/G_{sign,\mu}|,
\end{eqnarray}
For the simple and simple elliptic singularities observe
\begin{eqnarray}
|\{\textup{Stokes matrices}\}|&=&2^{\mu-1}\cdot 
|\{\textup{Stokes matrices}\}/G_{sign,\mu}| .\label{7.9}
\end{eqnarray}
The last equality holds because any Coxeter-Dynkin diagram is connected.
\end{remarks}

\begin{corollary}\label{t7.3}
For any simple singularity $|\BB(f)/G_{sign,\mu}|=\deg LL^{alg}$, and
this number is given in table \eqref{6.4}. The other numbers are as follows.
\begin{eqnarray}\label{7.10}
\begin{array}{lll}
  & |G_\Z| & |\{\textup{Stokes matrices}\}/G_{sign,\mu}| \\ 
\hline
A_\mu & 2(\mu+1) & (\mu+1)^{\mu-2}\\
D_4 & 36 & 9\\
D_\mu,\mu\geq 5 & 4(\mu-1) & (\mu-1)^{\mu-1}\\
E_6 & 24  & 2^7\cdot 3^3 = 3456\\
E_7 & 18  & 2\cdot 3^{10} = 118098\\
E_8 & 30  & 2\cdot 3^4\cdot 5^6 = 2531250
\end{array}
\end{eqnarray}
In the case of the simple elliptic singularities
\begin{eqnarray}\label{7.11}
\begin{array}{lll}
 & \deg(M^{alg}\to M^{mar}/G_\Z) & |\{\textup{Stokes matrices}\}/G_{sign,\mu}| \\
\hline 
\www E_6 & 6\cdot 2\cdot 3\cdot 3^2 =326 & 3^7\cdot 5\cdot 7=76545\\
\www E_7 & 6\cdot 1\cdot 4\cdot 2^2 = 96 & 2^{13}\cdot 5^3\cdot 7 =7168000\\
\www E_8 & 6\cdot 1\cdot 6\cdot 1^2 = 36 & 2^7\cdot 3^8\cdot 7\cdot 101
=593744256
\end{array}
\end{eqnarray}
Here $\deg(M^{alg}\to M^{mar}/G_\Z)$ means the generic degree.
\end{corollary}

{\bf Proof:} 
First we consider the simple singularities.
The bijection \eqref{7.3} gives
$$\deg LL^{alg}=|R_{Stokes}^0|=|\BB(f)/G_{sign,\mu}|.$$
%Observe also $|\BB(f)|=2^\mu\cdot |\BB(f)/G_{sign,\mu}|$.
The group $G_\Z$ acts on $M^{alg}=M^{mar}$ 
with kernel $\{\pm\id\}$.
This and the bijection \eqref{7.4} give 
\begin{eqnarray*}
|\{\textup{Stokes matrices}\}/G_{sign,\mu}|&=& 
|R_{Stokes}^0/G_\Z|=2\cdot |R_{Stokes}^0|/|G_\Z|\\
&=&2\cdot\deg LL^{alg}/|G_\Z|.
\end{eqnarray*}
The values $|G_\Z|$ can be found in \cite[Theorem 8.3 and Theorem 8.4]{He11}.
Together with \eqref{6.4}, this gives \eqref{7.10}.

Now we consider the simple elliptic singularities. Obviously
\begin{eqnarray*}
|R^0_{Stokes}/G_\Z|&=&\deg(LL:M^{mar}/G_\Z\to M_{LL}^{(\mu)})\\
&=& \frac{\deg LL^{alg}}{\deg(M^{alg}\to M^{mar}/G_\Z)}.
\end{eqnarray*}
Therefore the degree of the map $M^{alg}\to M^{mar}/G_\Z$
in the second column of \eqref{7.11}, 
the bijection \eqref{7.4} and the table \eqref{6.8}
give the third column of \eqref{7.11}.
The degree $\deg(M^{alg}\to M^{mar}/G_\Z)$ is calculated
in lemma \ref{t5.2}.
\hfill$\Box$

\section{Segre classes of smooth cone bundles}\label{s8}
\setcounter{equation}{0}

\noindent
The calculation of the degrees $\deg LL^{alg}$ for the simple
elliptic singularities in section \ref{s10} will use corollary \ref{t8.6} below.
For this corollary, we have to extend some notions and results
in \cite{Fu84}. We do not need new ideas, just some new details.
We follow closely this book.

In \cite[B.5]{Fu84} {\it cones} are defined in the category of algebraic
schemes as follows. $X$ is an algebraic scheme.
$S^\bullet=\sum_{d\geq 0}S^d$ is a graded sheaf of $\OO_X$-algebras
such that the canonical map $\OO_X\to S^0$ is an isomorphism,
$S^1$ is a coherent $\OO_X$-module, and $S^\bullet$ is 
(locally) generated by $S^1$ as an $\OO_X$-algebra.

Then $C:=\textup{Spec}(S^\bullet)$ with the projection $C\to X$ is a {\it cone}.
Its fibers are affine and come equipped with a $\C^*$-action.
The bundle of $\C^*$-orbits is $P(C):=\textup{Proj}(S^\bullet)$.
The projection $p_C:P(C)\to X$ is proper. 
The rational functions on $C$ which are homogeneous of degree $d$ induce
the line bundle $\OO_{P(C)}(d)$. 

If $f:Y\to X$ is a morphism, then the pull-back $f^*C=C\times_XY$ is the 
cone on $Y$ defined by the sheaf $f^*S^\bullet$ of $\OO_Y$-algebras.
If $C_1$ and $C_2$ are two cones on $X$ defined by $S_1^\bullet$ and
$S_2^\bullet$, their direct sum $C_1\oplus C_2$ is the cone on $X$
defined by the graded sheaf $S_1^\bullet\otimes S_2^\bullet$.

Now suppose that the cone $C$ is pure dimensional. Then its {\it Segre class}
is by \cite[Example 4.1.2]{Fu84}
\begin{eqnarray}\label{8.1}
s(C)=(p_C)_*\left(\sum_{i\geq 0}
c_1(\OO(1))^i\cap[P(C)]\right)\in A_*X.
\end{eqnarray}
Here $A_*P(C)$ and $A_*X$ are the spaces of cycles modulo rational 
equivalence \cite[1.3]{Fu84},
$[P(C)]\in A_*P(C)$, $(p_C)_*:A_*P(C) \to A_*X$ is the 
push-forward
\cite[1.4]{Fu84}, $\OO(1)$ is the canonical line bundle on $P(C)$,
and the Chern class $c_1(\OO(1))$ is understood in the operational sense,
as a map
\begin{eqnarray*}
c_1(\OO(1))\cap :A_kP(C)\to A_{k-1}P(C)
\end{eqnarray*}
\cite[3.2]{Fu84}. 

We are interested in the (more special and more general) situation where
the base $X$ is pure dimensional, the fibers are smooth, and the
fibration $C\to X$ is locally trivial, but where the condition that $S^1$
generates $S^\bullet$ is not necessarily satisfied.  The following 
definition fixes this situation.

\begin{definition}\label{t8.1}
For some $n\in\Z_{\geq 1}$, let ${\bf a}=(a_1,...,a_n)\in\Z_{\geq 1}^n$
with $a_1\leq a_2\leq ...\leq a_n$.
Let $R=\C[x_1,...,x_n]$ be the $\C$-algebra with the grading 
$R=R^\bullet =\sum_{d\geq 0}R^d$ such that $x_i\in R^{a_i}$.
Let $X$ be an algebraic pure dimensional scheme.
Let $S^\bullet=\sum_{d\geq 0}S^d$ be a graded sheaf of $\OO_X$-algebras
such that there is a covering of $X$ by open affine charts 
$U_i,i\in I,$ and there are isomorphisms 
$S|_{U_i}\cong \uuuu{R}_{U_i}:=\OO_{U_i}\otimes R$ of
graded $\OO_{U_i}$-algebras.

Then $C:=\textup{Spec}(S^\bullet)$ is a {\it smooth cone bundle}
({\it smooth} because the fibers of $C\to X$ are smooth,
{\it bundle} because $C\to X$ is locally trivial).
\end{definition}

By the next lemma, a smooth cone bundle comes equipped with a chain of
smooth cone subbundles and quotients, which are vector bundles.

\begin{lemma}\label{t8.2}
The situation in definition \ref{t8.1} is kept.
For $k\in\Z$ with $a_1\leq k\leq a_n+1$, 
define $I_k\subset S^\bullet$ as the sheaf of
homogeneous ideals generated by $S^1+...+S^{k-1}$ (so $I_{a_1}=\{0\}$),
define $S^\bullet_k:= S^\bullet/I_k$ as the quotient sheaf of
$\OO_X$-algebras with the induced grading, and define
$S_{k,sub}^\bullet\subset S^\bullet_k$ as the subring sheaf
of $S^\bullet_k$ generated by $S^k_k$
(obviously $S_{k,sub}^d=0$ if $k\not|d$).
Define $S^\bullet_{(k)}$ essentially as $S^\bullet_{k,sub}$, but
with the new grading $S^d_{(k)}:=S^{d\cdot k}_{k,sub}$. 

Then the $C_k:=\textup{Spec}(S^\bullet_k)$ are smooth cone bundles
on $X$ and form a chain
\begin{eqnarray}\label{8.2}
(\textup{zero section})=C_{a_n+1}\subset C_{a_n}\subset 
C_{a_n-1}\subset ...\subset C_{a_1+1}\subset C_{a_1}=C.
\end{eqnarray}
The cones $C_{(k)}:=\textup{Spec}(S^\bullet_{(k)})$ 
are vector bundles
on $X$ with $\rank C_{(k)}=|\{i\, |\, k=a_i\}|$ 
(so many of them may be 0,
and $\sum_{k=1}^{a_n}\rank C_{(k)}=\rank C$).
The smooth cone bundle $C_{k+1}$ is the kernel of 
the projection $\pr_k: C_k\to C_{(k)}$
(from the inclusion $S^\bullet_{k,sub}\hookrightarrow S^\bullet_k$).
\end{lemma}

The proof is clear. 

Given a smooth cone bundle $C\to X$, we want to define a {\it Segre class}.
The rational functions of degree $d$ on $C$ induce a sheaf
$\OO_{P(C)}(d)$. But $S^\bullet$ is in general not generated by $S^1$,
therefore $\OO_{P(C)}(1)$ is not necessarily invertible.
For example if $\gcd(a_1,...,a_n)>1$ then $S^d=0$ and 
$\OO_{P(C)}(d)=0$ if $\gcd(a_1,...,a_n)\not| d$. 
But for certain larger $d$, the sheaf $\OO_{P(C)}(d)$ is 
good enough.

\begin{definition}\label{t8.3}
The situation in definition \ref{t8.1} is kept.
Choose $d\in \lcm(a_1,...,a_n)\cdot\Z_{\geq 1}$, and define
a Segre class
\begin{eqnarray}\label{8.3}
s^{(d)}(C):= (p_C)_*\left( \sum_{i\geq 0}
\left(\frac{c_1(\OO(d))}{d}\right)^i
\cap \frac{[P(C)]}{\gcd(a_1,...,a_n)}\right) \in A_*^\Q X,
\end{eqnarray}
here $A_*^\Q X=A_*X\otimes_\Z\Q$, $A_*^\Q P(C)
=A_*P(C)\otimes_\Z\Q$,
and $c_1(\OO(d))\cap:A_k^\Q X P(C)\to A_{k-1}^\Q P(C)$, 
$(p_C)_*:A_k^\Q  P(C)\to A_k^\Q X$.
\end{definition}

Part (b) in the following proposition generalizes
\cite[Proposition 4.1 (a)]{Fu84}.

\begin{proposition}\label{t8.4}
The situation in definition \ref{t8.1} is kept.

(a) $s^{(d)}(C)$ is independent of the choice of $d$
and is called $s^{(scb)}(C)$. 
If $a_1=...=a_n=1$, then $C$ is a vector bundle
and $s^{(scb)}(C)$ is the classical Segre class in 
\cite[Example 4.1.2]{Fu84}.

(b) 
\begin{eqnarray}\label{8.4}
s^{(scb)}(C)&=& \frac{1}{a_1...a_n}\cdot \prod_{k=a_1}^{a_n}
c^{pol}_{1/k}(C_{(k)})^{-1}\cap [X],
\end{eqnarray}
here $C_{(k)}$ is the vector bundle associated to $C$ in lemma \ref{t8.2},
and $c^{pol}_t(C_{(k)})=c_0+tc_1+t^2c_2+...$ is its Chern polynomial (in \eqref{8.4}
the variable $t$ is replaced by the number $\frac{1}{k})$.
\end{proposition}

{\bf Proof:}
(a) The independence will follow from the formula \eqref{8.4}
in (b). If $a_1=...=a_n=1$, then $C$ is a vector bundle,
$\OO(d)=\OO(1)^{\otimes d}$, $c_1(\OO(d))=d\cdot c_1(\OO(1))$,
and the definition \eqref{8.3} agrees with 
\cite[Example 4.1.2]{Fu84}.

(b) This will be proved by induction on the dimension $n$ of the
fibers of the smooth cone bundle. We carry out the first step of
the induction. Part of it is close to \cite[Example 4.1.5]{Fu84}.

By the splitting construction \cite[Proof of Theorem 3.2]{Fu84},
there is a flat morphism $f:Y\to X$ such that 
$f^*:A_*X\to A_*Y$ is injective and 
$f^*C_{(a_1)}$ has a filtration by subbundles
\begin{eqnarray}\label{8.5}
f^*C_{(a_1)} = E_r\supset E_{r-1}\supset ...\supset 
E_1\supset E_0=0
\end{eqnarray}
with line bundle quotients $L_i=E_i/E_{i-1}$ 
(and $r=|\{i\, |\, a_i=a_1\}|=\rk C_{(a_1)}$).

The smooth cone bundle $F:=f^*(C)$ in $Y$
contains the smooth cone bundle 
$G:=\ker(F\to L_r)$ in codimension one,
where $F\to L_r$ is the composition of the projections
$f^*\pr_{a_1}:F\to E_r$ and $E_r\to L_r$. 
Denote by $\LL^\bullet$ the graded sheaf of $\OO_X$-algebras with
$L_r=\textup{Spec}(\LL^\bullet)$.
The projection $F\to L_r$ with kernel $G$ corresponds to an embedding
$\LL^\bullet\hookrightarrow f^*S_{k,sub}^\bullet\subset f^*S^\bullet$.
Denote by $p_F:P(F)\to Y$ the projection.

The line bundle $(p_F)^* L_r^{\otimes d/a_1}\otimes 
\OO_{P(F)}(d)$ has a global
section $\sigma$: If $U\subset Y$ is an open affine chart and
$f^*S^\bullet|_U =\OO_U\otimes R\cdot e_1$ is a trivialization and 
$\LL^\bullet|_U\cong \OO_U\otimes \C[x_1]$, then
$\sigma|_{(p_F)^{-1}(U)} \sim (p_F)^*e_1^{\otimes d/a_1}\otimes x_1^{d/a_1}$.
Its zero-scheme in $P(F)$ may be identified with $P(G)$
with multiplicity $\frac{d}{a_1}\cdot \frac{\gcd(a_1,...,a_n)}{\gcd(a_2,...,a_n)}$,
equivalently, 
\begin{eqnarray}\label{8.6}
\frac{d}{a_1}\cdot \frac{\gcd(a_1,...,a_n)}{\gcd(a_2,...,a_n)}
\cdot [P(G)] 
= c_1((p_F)^* L_r^{\otimes d/a_1}\otimes \OO_{P(F)}(d))
\cap [P(F)].
\end{eqnarray}
Here $\gcd(a_1,...,a_n)$ and $\gcd(a_2,...,a_n)$ are the
sizes of the kernels of the $\C^*$-actions on $F$ and 
on $G$.

Now we want to calculate $s^{(d)}(G)$ in terms of $s^{(d)}(F)$ and
the value at $t=\frac{1}{a_1}$ of the Chern polynomial $c^{pol}_t(L_r)$.
For this observe that the closed embedding 
$i:P(G)\hookrightarrow P(F)$
is proper. The formula $i^*\OO_{P(F)}(d)=\OO_{P(G)}(d)$ and the
projection formula for Chern classes \cite[Theorem 3.2 (c)]{Fu84} will be used.
\begin{eqnarray*}
&&\frac{1}{a_1}\cdot s^{(d)}(G)\\
&=& \frac{1}{a_1}\cdot (p_G)_* \left(\sum_{i\geq 0}
\left(\frac{c_1(\OO_{P(G)}(d))}{d}\right)^i
\cap \frac{[P(G)]}{\gcd(a_2,...,a_n)}\right) \\
&=& \frac{1}{d} \cdot (p_F)_* \left(\sum_{i\geq 0} 
\left(\frac{c_1(\OO_{P(F)}(d))}{d}\right)^i 
\cap c_1\Bigl((p_F)^* L_r^{\otimes d/a_1}
\otimes \OO_{P(F)}(d)\Bigr) \right.\\
&& \left. \hspace*{7cm}\cap
\frac{[P(F)]}{\gcd(a_1,...,a_n)}\right) \\
&=& (p_F)_* \left(\sum_{i\geq 0} \left( \frac{c_1(\OO_{P(F)}(d))}{d}\right)^i 
\cap \left(\frac{1}{a_1}c_1((p_F)^* L_r)+\frac{1}{d}c_1(\OO_{P(F)}(d)\right)\right. \\
&&  \left. \hspace*{7cm}\cap
\frac{[P(F)]}{\gcd(a_1,...,a_n)}\right) \\
&=& (p_F)_*\left(\sum_{i\geq 0}\left(\frac{c_1(\OO_{P(F)}(d))}{d}\right)^i
\cap c^{pol}_{1/a_1}((p_F)^* L_r)\cap
\frac{[P(F)]}{\gcd(a_1,...,a_n)}\right) \\
&=& c_{1/a_1}^{pol}(L_r)\cap s^{(d)}(F).
\end{eqnarray*}
In the second to last equality, the term 
$(p_F)_*\Bigl(c_0(p_F)^*(L_r)\cap 
\frac{[P(F)]}{\gcd(a_w,...,a_n)}\Bigr)$ was added. 
This term vanishes,
as $[P(F)]\in A_{\dim P(F)}$ is mapped by $(p_F)_*$ to 
$A_{\dim P(F)}Y$, which is zero because 
$\dim P(F)=\dim P(G)+1\geq \dim Y+1$. Therefore 
\begin{eqnarray}\label{8.7}
s^{(d)}(F) &=& \frac{1}{a_1}\cdot c^{pol}_{1/a_1}(L_r)^{-1}\cap s^{(d)}(G).
\end{eqnarray}
By induction and the product formula 
$c^{pol}_t(L_r)\cdot c^{pol}_t(E_{r-1}) = c^{pol}_t(f^* C_{(a_1)})$ we obtain
\begin{eqnarray}
s^{(d)}(F) &=& \frac{1}{a_1^r}\cdot \prod_{j=1}^r 
c_{1/a_1}^{pol}(L_j)^{-1}\cap s^{(d)}(f^*C_{(a_1+1)})\nonumber\\
&=& \frac{1}{a_1^r} \cdot c_{1/a_1}^{pol}(f^*C_{(a_1)})^{-1}
\cap s^{(d)}(f^*C_{(a_1+1)})\nonumber\\
&=& \frac{1}{a_1...a_n}\cdot\prod_{k=a_1}^{a_n}
c_{1/k}^{pol}(f^*C_{(k)})^{-1})\cap [Y]  \nonumber \\
&=& \frac{1}{a_1...a_n}\cdot f^*\left( \prod_{k=a_1}^{a_n}
c_{1/k}^{pol}(C_{(k)})^{-1})\cap [X]\right) . \label{8.8}
\end{eqnarray}
Now the injectivity of $f^*:A_*X\to A_*Y$ gives
\begin{eqnarray}\label{8.9}
s^{(d)}(C) &=& \frac{1}{a_1...a_n}\cdot\prod_{k=a_1}^{a_n}
c_{1/k}^{pol}(C_{(k)})^{-1})\cap [X].
\end{eqnarray}

\medskip
It rests to settle the beginning of the induction.
Consider the case $n=1$ and choose $d\in a_1\cdot \Z_{\geq 1}$.
Then $P(C)=X$, $p_C=\id$, $C=C_{(a_1)}\cong \OO_{P(C)}(-a_1)$, and 
\begin{eqnarray*}
\sum_{i\geq 0}\left(\frac{c_1(\OO_X(d))}{d}\right)^i &=&
\left(1-\frac{1}{d}c_1(\OO_X(d))\right)^{-1} 
= \left(1+\frac{1}{d}c_1(\OO_X(-d))\right)^{-1} \\
&=& \left(1+\frac{1}{a_1}c_1(\OO_X(-a_1))\right)^{-1}
=c^{pol}_{1/a_1}(C_{(a_1)})^{-1}.
\end{eqnarray*}
Therefore
\begin{eqnarray*}
s^{(d)}(C)&=& (p_C)_*\left( \sum_{i\geq 0}\left(\frac{c_1(\OO_{P(C)}}{d}\right)^i
\cap \frac{[P(C)]}{a_1}\right) \\
&=& \frac{1}{a_1} (c_{1/a_1}^{pol}(C_{(a_1)}))^{-1}\cap [X]. \hspace*{2cm}\Box
\end{eqnarray*}

\bigskip

As in \cite{Fu84}, a variety means a reduced and irreducible scheme.
A smooth cone bundle $C\to X$ on a variety is also a variety.
In the following $X\subset C$ means the embedding as zero section.

Suppose that $C_1\to X_1$ and $C_2\to X_2$ are smooth cone bundles
on complete varieties $X_1$ and $X_2$ with $\dim C_1=\dim C_2$
and $\dim X_1\geq \dim X_2$, 
and suppose that $f:C_1\to C_2$ is a $\C^*$-equivariant proper morphism
such that $f^{-1}(X_2)=X_1$ and the restriction
$f^{C-X}:C_1-X_1\to C_2-X_2$ is finite. Then $f$ and the restrictions
$f^{C-X}$ and $f^X:X_1\to X_2$ are surjective, and $f$ has a finite
degree $\deg f=[K(C_1):K(C_2)]$, which is the number of preimages
of a generic point in $C_2-X_2$.  
In \cite[Definition 1.4]{Fu84} a degree $\int_{X_i}:A_0X_i\to \Z$
is defined and extended by $\int_{X_i}:A_kX_i\to 0$ for $k>0$ to 
$\int_{X_i}:A_*X_i\to\Z$. It also extends to $\int_{X_i}A_*^\Q X_i\to\Q$.

\begin{proposition}\label{t8.5}
In the situation just described
\begin{eqnarray}\label{8.10}
f^X_* s^{(scb)}(C_1) &=& \deg f\cdot s^{(scb)}(C_2),\\
\int_{X_1} s^{(scb)}(C_1) &=& \deg f\cdot \int_{X_2}s^{(scb)}(C_2).\label{8.11}
\end{eqnarray}
\end{proposition}

{\bf Proof:}
Denote by ${\bf a}=(a_1,...,a_{n_1})$ respectively by
${\bf v}=(v_1,...,v_{n_2})$ the weights of $C_1$ respectively $C_2$.
Denote $\www w:=\gcd(a_i),\www v:=\gcd(v_i), d_1:=\lcm(a_i), d_2:=\lcm(v_i)$.
Because $f$ is $\C^*$-equivariant and does not map $C_1$ to $X_2$, 
$\www w$ divides $\www v$. The map $f^{C-X}$ induces a finite 
(and surjective) morphism $f^{PC}:P(C_1)\to P(C_2)$ with
\begin{eqnarray*}
\deg f^{PC}=\left(\frac{\www v}{\www w}\right)^{-1}\cdot \deg f.
\end{eqnarray*}
Furthermore
\begin{eqnarray*}
(f^{PC})^* \OO_{P(C_2)}(d) &=& \OO_{P(C_1)}(d),\\
(f^{PC})_* [P(C_1)] &=& \deg f^{PC}\cdot [P(C_2)].
\end{eqnarray*}
Choose $d\in\lcm(d_1,d_2)\cdot\Z_{\geq 1}$. Then
\begin{eqnarray*}
\OO_{P(C_i)}(d)=\OO_{P(C_i)}(d_i)^{\otimes d/d_i},\quad
\frac{1}{d_i} c_1(\OO_{P(C_i)}(d_i)) = \frac{1}{d} c_1(\OO_{p(C_i)}(d)).
\end{eqnarray*}
Therefore 
\begin{eqnarray*}
&&f^X_* s^{(d)}(C_1)\\
&=& f^X_*\circ (p_{C_1})_*\left( \sum_{i\geq 0}
\left( \frac{c_1(\OO_{P(C_1)}(d))}{d}\right)^i
\cap \frac{[P(C_1)]}{\www w}\right) \\
&=& (p_{C_2})_*\circ (f^{PC})_* \left( \sum_{i\geq 0} 
\left( \frac{c_1((f^{PC})^*\OO_{P(C_2)}(d))}{d}\right)^i
\cap \frac{[P(C_1)]}{\www w}  \right) \\
&=& (p_{C_2})_*\left( \sum_{i\geq 0}
\left( \frac{c_1(\OO_{P(C_2)}(d)}{d}\right)^i 
\cap \frac{(f^{PC})_* [P(C_1)]}{\www w} \right) \\
&=& \deg f\cdot (p_{C_2})_*\left( \sum_{i\geq 0}
\left( \frac{c_1(\OO_{P(C_2)}(d)}{d}\right)^i 
\cap \frac{[P(C_2)]}{\www v}\right) \\
&=& \deg f\cdot s^{(d)}(C_2).
\end{eqnarray*}
With the functoriality $\int_{X_1}\alpha = \int_{X_2} f^X_*\alpha$ 
\cite[Definition 1.4]{Fu84}, we obtain
\begin{eqnarray*}
\int_{X_1} s^{(scb)} (C_1) = \deg f\cdot \int_{X_2} s^{(scb)}(C_2).
\hspace*{1cm}\Box
\end{eqnarray*}

If $\int_{X_2}s^{(scb)}(C_2)\neq 0$, then \eqref{8.11} can be used
to calculate $\deg f$. We will use it in the following case.

\begin{corollary}\label{t8.6}
Keep the situation in and before proposition \ref{t8.5}.
Suppose additionally that $X_1$ is a smooth complete curve 
and $X_2$ is a point.
The weights of $C_1$ are denoted 
$(a_1,..,a_{n_1})$, the weights of $C_2$ are denoted
$(b_1,...,b_{n_2})$. The vector bundles on $X_1$ 
associated to $C_1$ in lemma \ref{t8.2} are 
denoted by $C_{1,(k)}$, $a_1\leq k\leq a_n$
Then $n_2=n_1+1$,
\begin{eqnarray}\label{8.12}
\int_{X_2} s^{(scb)}(C_2) &=& \frac{1}{b_1...b_{n_2}}>0,\\
\int_{X_1} s^{(scb)}(C_1) &=& \frac{1}{a_1...a_{n_1}}
\cdot \left(-\sum_{k=a_1}^{a_{n_1}}
\frac{1}{k}\deg C_{1,(k)}\right),\label{8.13}\\
\deg f &=& \frac{b_1...b_{n_2}}{a_1...a_{n_1}}\cdot \left(-\sum_{k=a_1}^{a_{n_1}}
\frac{1}{k}\deg C_{1,(k)}\right).\label{8.14}
\end{eqnarray}
\end{corollary}

{\bf Proof:}
This follows with proposition \ref{t8.4} (b) and proposition \ref{t8.5}.
\hfill$\Box$

\section{Extension to $\lambda=0$ 
of the Lyashko-Looijenga map for the simple
elliptic singularities}\label{s9}
\setcounter{equation}{0}

\noindent
Here we will do the first and biggest step in the proof
of theorem \ref{t6.3}. The $\lambda$-parameter space
$\C-\{0,1\}$ contains the punctured disk
$\Delta^*=\Delta-\{0\}$, where $\Delta=\{z\in\C\, |\, 
|z|<1\}$. Define $c:=3$ for $\www E_6$ and $\www E_8$
and $c:=2$ for $\www E_7$, and define the $c$-fold
coverings
\begin{eqnarray*}
\rho_{naive}:\C^{\mu-1}\times\Delta^* &\to&
\C^{\mu-1}\times\Delta^*,\\
(t',\kappa)&\mapsto& (t',\kappa^c)=(t',\lambda),\\
R_{naive}:\C^{n+1}\times\C^{\mu-1}\times\Delta^*&\to&
\C^{n+1}\times\C^{\mu-1}\times\Delta^*,\\
(x,t',\kappa)&\mapsto& (x,t',\kappa^c)=(x,t',\lambda).
\end{eqnarray*}
We will glue fibers above $\kappa=0$ into
$\C^{n+1}\times\C^{\mu-1}\times\Delta^*$
and $\C^{\mu-1}\times\Delta^*$ such that
$F^{alg}\circ R_{naive}$, its critical space
and its Lyashko-Looijenga map extend well to $\kappa=0$.
This is the content of the following theorem \ref{t9.1}
and its long proof.

\begin{theorem}\label{t9.1}
Consider for each of the three families of simple
elliptic singularities and their unfoldings the following
spaces and maps.

For $\www E_6$: 
\begin{eqnarray}
(x,y,s,\kappa) &=& (x_0,...,x_n,y_0,y_1,y_2,s_1,...,s_7,\kappa)
\nonumber\\
&\in& \C^{n+1}\times \C^3\times \C^7\times\Delta
=\C^{n+11}\times\Delta^*,\nonumber\\
Y&:=&\{(x,y,s,\kappa)\in\C^{n+11}\times\Delta\, |\, 
x_0(x_1+s_5)=\kappa y_0,\nonumber\\ 
&&(x_1+s_5)y_0=\kappa y_1,x_0x_2^2=\kappa^2y_2\},\label{9.1}\\
\pr_\mu:Y&\to&\C^7\times\Delta,\, (x,y,s,\kappa)\mapsto(s,\kappa),
\nonumber
\end{eqnarray}
\begin{eqnarray}
\rho:\C^7\times\Delta^*&\to&\C^7\times\Delta^*,\nonumber\\
(s,\kappa)&\mapsto& (s_1,\kappa^2s_2+\kappa s_5s_6-s_5^2,s_3,
s_4,\kappa^3s_5,\nonumber\\
&& \kappa s_6-2s_5,s_7,\kappa^3),\label{9.2}\\
R:Y\cap \C^{n+11}\times\Delta^*&\to&\C^{n+8}\times\Delta^*,
\nonumber\\
(x,y,s,\kappa)&\mapsto& (\kappa^{-2}x_0,x_1,...,x_n,
\rho(s,\kappa)).\label{9.3}
\end{eqnarray}

For $\www E_7$: 
\begin{eqnarray}
(x,y,s,\kappa) &=& (x_0,...,x_n,y,s_1,...,s_8,\kappa)
\nonumber\\
&\in& \C^{n+1}\times \C\times \C^8\times\Delta
=\C^{n+10}\times\Delta^*,\nonumber\\
Y&:=&\{(x,y,s,\kappa)\in\C^{n+10}\times\Delta\, |\, 
x_0x_1=\kappa y\},\label{9.4}\\
\pr_\mu:Y&\to&\C^8\times\Delta,\, (x,y,s,\kappa)\mapsto(s,\kappa),
\nonumber
\end{eqnarray}
\begin{eqnarray}
\rho:\C^8\times\Delta^*&\to&\C^8\times\Delta^*,\nonumber\\
(s,\kappa)&\mapsto& (s_1,\kappa s_2,s_3,
\kappa^2s_4,s_5,s_6,\kappa s_7,s_8,\kappa^2),\label{9.5}\\
R:Y\cap \C^{n+10}\times\Delta^*&\to&\C^{n+9}\times\Delta^*,
\nonumber\\
(x,y,s,\kappa)&\mapsto& (\kappa^{-1}x_0,x_1,...,x_n,
\rho(s,\kappa)).\label{9.6}
\end{eqnarray}

For $\www E_8$: 
\begin{eqnarray}
(x,y,s,\kappa) &=& (x_0,...,x_n,y,s_1,...,s_9,\kappa)
\nonumber\\
&\in& \C^{n+1}\times \C\times \C^9\times\Delta
=\C^{n+11}\times\Delta^*,\nonumber\\
Y&:=&\{(x,y,s,\kappa)\in\C^{n+11}\times\Delta\, |\, 
(x_0-\frac{1}{2}s_9)x_1=\kappa y\},\label{9.7}\\
\pr_\mu:Y&\to&\C^9\times\Delta,\, (x,y,s,\kappa)\mapsto(s,\kappa),
\nonumber
\end{eqnarray}
\begin{eqnarray}
\rho:\C^9\times\Delta^*&\to&\C^9\times\Delta^*,\nonumber\\
(s,\kappa)&\mapsto& (s_1,\kappa s_2,\kappa^2 s_3,\nonumber\\
&&s_4-\frac{1}{2}\kappa^{-1}s_6s_9-\frac{1}{4}\kappa^{-1}s_7s_9^2
-\frac{1}{16}\kappa^{-1}s_9^4,\kappa^3s_5,\nonumber\\
&&s_6,\kappa s_7,s_8-\frac{1}{4}\kappa^{-2}s_9^2,\kappa^{-1}s_9,
\kappa^3),\label{9.8}\\
R:Y\cap \C^{n+11}\times\Delta^*&\to&\C^{n+10}\times\Delta^*,
\nonumber\\
(x,y,s,\kappa)&\mapsto& (\kappa^{-1}x_0,x_1,...,x_n,
\rho(s,\kappa)).\label{9.9}
\end{eqnarray}

(a) The $\C^*$-action on $\C^{n+11}\times\Delta$ for
$\www E_6$ and $\www E_8$ and on $\C^{n+10}\times\Delta$
for $\www E_7$ with the following weights restricts to
a $\C^*$-action on $Y$,
\begin{eqnarray}
&&\deg_w x_i=w_i,\ \deg_w s_i=\deg_w t_i,\ 
\deg_w\kappa=\deg_w \lambda=0,\nonumber\\
&&\textup{for }\www E_7\textup{ and }\www E_8:\quad 
\deg_w y=w_0+w_1, \label{9.10}\\
&&\textup{ for }\www E_6:\quad\left\{\begin{array}{l}
\deg_w y_0=w_0w_1,\  \deg_w y_1=w_0+2w_1,\\
\deg_w y_2=w_0+2w_2.\end{array}\right. \nonumber
\end{eqnarray}
The map $R$ is $\C^*$-equivariant with respect to this
$\C^*$-action and the natural $\C^*$-action on the image
space with coordinates $(x,t',\lambda)$.

\medskip
(b) The maps $\rho$ and $R$ are coverings of degree $c$.
Especially, for each fixed 
$(s,\kappa)\in\C^{\mu-1}\times\Delta^*$,
\begin{eqnarray}
R:\pr_\mu^{-1}((s,\kappa))&\stackrel{\cong}{\to}&
\C^{n+1}\times\{\rho(s,\kappa)\}.\label{9.11}
\end{eqnarray}

\medskip
(c) The pull back $F^{alg}\circ R$ extends from
$Y\cap \C^{n+(10\textup{ or }11)}\times\Delta^*$
holomorphically to $\kappa=0$, that means, to a function
$(F^{alg}\circ R)^{ext}:Y\to \C$.

\medskip
(d) Let $C^{alg}:=\{(x,t',\lambda)\in\C^{n+1}\times M^{alg}\, |\, 
\frac{\paa F^{alg}}{\paa t_i}=\frac{\paa F^{alg}}{\paa \lambda}
=0\}$
be the critical space of the unfolding $F^{alg}$. 
Consider the closure 
$\oooo{R^{-1}(C^{alg})}$ in $Y$ of the pull back by $R$ of
$C^{alg}\cap \C^{n+1}\times \C^{\mu-1}\times\Delta^*$.
The restriction
\begin{eqnarray}\label{9.12}
\pr_\mu:\oooo{R^{-1}(C^{alg})}\to\C^{\mu-1}\times\Delta,\ 
(x,y,s,\kappa)\mapsto (s,\kappa)
\end{eqnarray}
is finite and flat of degree $\mu$.

\medskip
(e) The composition $LL^{alg}\circ\rho:\C^{\mu-1}\times\Delta^*
\to M_{LL}^{(\mu)}$ of the Lyashko-Looijenga map
$LL^{alg}$ with $\rho$ extends holomorphically to
$\C^{\mu-1}\times\Delta$. The restriction
\begin{eqnarray}\label{9.13}
(LL^{alg}\circ\rho)^{ext}:\C^{\mu-1}\times\Delta -
\{(s,\kappa)\, |\, s_2=...=s_{\mu-1}=0\}\\
\to M_{LL}^{(\mu)}-M_{LL,0}^{(\mu)}\hspace*{4cm}\nonumber
\end{eqnarray}
is finite and flat onto its image. And 
$\C^{\mu-1}\times\{0\}$ is mapped by
$(LL^{alg}\circ\rho)^{ext}$ to $D_{LL}^{(\mu)}$.
\end{theorem}

The rest of this section is devoted to the proof
of this theorem.

\bigskip

{\bf Proof:}
(a) This follows from comparison of the formulas
\eqref{9.1}--\eqref{9.9} with the weights in remark
\ref{t6.4} (ii) and in \eqref{9.10}.

\medskip
(b) The definition of $Y$ shows that for $\kappa\in\Delta^*$
$(x_0,...,x_n)$ serve as coordinates on
$\pr_\mu^{-1}((s,\kappa))$ and that this is isomorphic
to $\C^{n+1}$.

The following three statements show that $\rho$ and $R$
are coverings of degree $c$.
The last component of $\rho$ is $\rho_\mu(s,\kappa)=\kappa^c$.
Each other component $\rho_i$ has a nonvanishing linear term
in $s_i$ (and in the case of $\www E_6$ the linear terms
of $\rho_5$ and $\rho_6$ are $\kappa^3s_5$
and $\kappa s_6-2s_5$).
The map $R$ restricts to a linear isomorphism
$R:\pr_\mu^{-1}((s,\kappa))\to 
\C^{n+1}\times\{\rho((s,\kappa))\}$
for $(s,\kappa)\in\C^{\mu-1}\times \Delta^*$.

\medskip
(c) The pull back $F^{alg}\circ R
=\bigl(f_\lambda(x)+\sum_{i=1}^{\mu-1}m_it_i\bigr)\circ R$
can be written as follows.

For $\www E_6$: 
\begin{eqnarray}
F^{alg}\circ R&=& 
\kappa^3 \kappa^{-4}x_0^2x_1 
-(\kappa^3+1)\kappa^{-2}x_0x_1^2 +x_1^3
-\kappa^{-2}x_0x_2^2 +\sum_{i=3}^n x_i^2 \nonumber\\
&&+s_1+\kappa^{-2}x_0(\kappa^2s_2+\kappa s_5s_6-s_5^2)
+x_1s_3+x_2s_4 \nonumber\\
&&+ \kappa^{-4}x_0^2 \kappa^3s_5
+\kappa^{-2}x_0x_1 (\kappa s_6-2s_5)+x_1x_2s_7  \nonumber\\
&=& \kappa^{-1}x_0(x_0+s_6)(x_1+s_5) 
-\kappa^{-2}x_0(x_1+s_5)^2
-\kappa x_0x_1^2 +x_1^3 \nonumber\\
&& -\kappa^{-2}x_0x_2^2  +\sum_{i=3}^n x_i^2 +s_1
+x_0s_2 + x_1s_3 +x_2s_4 +x_1x_2s_7 \nonumber\\
&=& (x_0+s_6)y_0 -y_1 -\kappa x_0x_1^2 +x_1^3 -y_2
+\sum_{i=3}^n x_i^2\\
&&+ s_1+x_0s_2+x_1s_3+x_2s_4+x_1x_2s_7.\nonumber
\end{eqnarray}

For $\www E_7$: 
\begin{eqnarray}
F^{alg}\circ R&=& 
\kappa^2 \kappa^{-3}x_0^3x_1 
-(\kappa^2+1)\kappa^{-2}x_0^2x_1^2
+\kappa^{-1}x_0x_1^3 +\sum_{i=2}^n x_i^2 \nonumber\\
&&+s_1+\kappa^{-1}x_0\kappa s_2
+x_1s_3+ \kappa^{-2}x_0^2\kappa^2 s_4 \nonumber\\
&&+ \kappa^{-1}x_0x_1 s_5
+ x_1^2s_6 +\kappa^{-2}x_0^2x_1\kappa s_7
+\kappa^{-1}x_0x_1^2 s_8 \nonumber\\
&=& x_0^2y -(\kappa^2+1)y^2+ x_1^2y
+\sum_{i=2}^n x_i^2 \\
&&+ s_1+x_0s_2+x_1s_3+x_0^2s_4+ ys_5+x_1^2s_6
+x_0y s_7+x_1ys_8.\nonumber
\end{eqnarray}

For $\www E_8$: 
\begin{eqnarray}
F^{alg}\circ R&=& 
\kappa^3 \kappa^{-4}x_0^4x_1 
-(\kappa^3+1)\kappa^{-2}x_0^2x_1^2
+x_1^3 +\sum_{i=2}^n x_i^2 \nonumber\\
&&+s_1+\kappa^{-1}x_0\kappa s_2
+\kappa^{-2}x_0^2\kappa^2 s_3\nonumber\\
&&+x_1(s_4-\frac{1}{2}\kappa^{-1}s_6s_9 
-\frac{1}{4}\kappa^{-1} s_7s_9^2 
-\frac{1}{16} \kappa^{-1}s_9^4) \nonumber\\
&&+ \kappa^{-3}x_0^3\kappa^3 s_5
+ \kappa^{-1}x_0x_1s_6 +\kappa^{-2}x_0^2x_1\kappa s_7\nonumber\\
&& +x_1^2(s_8-\frac{1}{4}\kappa^{-2}s_9^2) 
+ \kappa^{-1}x_0x_1^2\kappa^{-1}s_9 \nonumber\\
&=& (x_0^3+x_0^2\frac{1}{2}s_9 +x_0\frac{1}{4}s_9^2 
+\frac{1}{8}s_9^3 + (x_0+\frac{1}{2}s_9)s_7 +s_6)y  \nonumber\\
&&-\kappa x_0^2x_1^2-y^2+ x_1^3
+\sum_{i=2}^n x_i^2\\
&&+ s_1+x_0s_2+x_0^2s_3+x_1s_4+ x_0^3s_5+x_1^2s_8.\nonumber
\end{eqnarray}

In all three cases, the terms after the last equality 
sign are in $\C[x,y,s,\kappa]$.

\medskip
(d) We call {\it $y$-relations} the elements
\begin{eqnarray}
\left. \begin{array}{l}
x_0(x_1+s_5)-\kappa y_0,\ (x_1+s_5)y_0-\kappa y_1,\\ 
x_0(x_1+s_5)^2-\kappa^2 y_1, x_0x_2^2-\kappa^2 y_2
\end{array}\right\} && \textup{for }\www E_6,\nonumber \\
x_0x_1-\kappa y && \textup{for }\www E_7,\label{9.17} \\
(x_0-\frac{1}{2}s_9)x_1-\kappa y && \textup{for }\www E_8
\nonumber
\end{eqnarray}
in $\C[x,y,s,\kappa,\kappa^{-1}]$ and in 
$\C[x,y,s,\kappa]$.
The compositions $\frac{\paa F^{alg}}{\paa x_i}\circ R$
of the partial derivatives of $F^{alg}$ with $R$ are
in $\C[x,y,s,\kappa,\kappa^{-1}]$.
We consider the following ideals,
\begin{eqnarray}
I_0&:=& \left( \frac{\paa F^{alg}}{\paa x_i}\circ R,
y\textup{-relations}\right)\subset \C[x,y,s,\kappa,\kappa^{-1}],
\nonumber\\
I_1&:=& I_0\cap \C[x,y,s,\kappa],\nonumber\\
I_2&:=& \{g(x,y,s,0)\, |\, g(x,y,s,\kappa)\in I_1\}
\subset \C[x,y,s], \label{9.18}\\
I_3&:=& \{g(x,y,0)\, |\, g(x,y,s)\in I_2\}\subset \C[x,y].
\nonumber
\end{eqnarray}

We will calculate generating elements of these ideals.
Then we will show $\dim \C[x,y]/I_3=\mu$.
This is sufficient for (d) because of the following.
$C^{alg}\subset \C^{n+1}\times M^{alg}$ and
$\oooo{R^{-1}(C^{alg})}\subset Y$ are invariant under the
$\C^*$-actions. Therefore it is sufficient to show that
the restriction of \eqref{9.12} to $s=0$
is finite and flat of degree $\mu$. This holds above 
$\Delta^*$. For $\kappa=0$ is is equivalent to
$\dim\C[x,y]/I_3=\mu$.

$I_1$ determines $\oooo{R^{-1}(C^{alg})}\subset Y$,
and $I_2$ determines $\oooo{R^{-1}(C^{alg})}\subset Y
\cap \C^{n+(10\textup{ or }11)}\times\{0\}$.
The information below on $I_2$ will also be useful in the 
proof of part (e).

\medskip
{\bf The case $\www E_6$:} 
\begin{eqnarray*}
\frac{\paa F^{alg}}{\paa x_0}\circ R &=& 
2\kappa x_0x_1 -(\kappa^3+1)x_1^2-x_2^2 
+(\kappa^2s_2+\kappa s_5s_6-s_5^2)\\
&&+ 2\kappa x_0s_5 +x_1(\kappa s_6-2s_5)\\
&\stackrel{y\textup{-relations}}{\equiv} & 
2\kappa^2 y_0-\kappa^3 x_1^2 -
(x_1+s_5-\kappa s_6)(x_1+s_5)\\
&& -x_2^2 +\kappa^2s_2.
\end{eqnarray*}
\begin{eqnarray*}
\frac{\paa F^{alg}}{\paa x_0}\circ R \ \&\  y\textup{-relations}
&\Rightarrow& (x_1+s_5)^2 +x_2^2\in I_2,\\
x_0\cdot \frac{\paa F^{alg}}{\paa x_0}\circ R \ \&\  
y\textup{-relations}
&\Rightarrow& 2x_0y_0-y_1+y_0s_6-y_2+x_0s_2\in I_2.
\end{eqnarray*}

\begin{eqnarray*}
\frac{\paa F^{alg}}{\paa x_1}\circ R &=& 
\kappa^{-1} x_0^2 -2(\kappa^3+1)\kappa^{-2}x_0x_1+3x_1^2 \\
&& + s_3 +\kappa^{-2}x_0(\kappa s_6-2s_5)+x_2s_7\\
&\stackrel{y\textup{-relations}}{\equiv} & 
\kappa^{-1}x_0^2 -2\kappa x_0x_1 -2\kappa^{-1}y_0 +3x_1^2
+s_3\\
&&+\kappa^{-1} x_0s_6 +x_2s_7.
\end{eqnarray*}
\begin{eqnarray}
\frac{\paa F^{alg}}{\paa x_1}\circ R \ \&\  y\textup{-relations}
&\Rightarrow& x_0^2-2y_0+x_0s_6\in I_2,\nonumber\\
(x_1+s_5)\cdot \frac{\paa F^{alg}}{\paa x_1}\circ R \ \&\  
y\textup{-relations}
&\Rightarrow& x_0y_0-2y_1+3(x_1+s_5)x_1^2\nonumber\\
+(x_1+s_5)s_3+y_0s_6
&&+(x_1+s_5)x_2s_7\in I_2.\label{9.19}
\end{eqnarray}

\begin{eqnarray*}
\frac{\paa F^{alg}}{\paa x_2}\circ R &=& 
-2\kappa^{-2}x_0x_2+s_4+x_1s_7.
\end{eqnarray*}
\begin{eqnarray*}
\frac{\paa F^{alg}}{\paa x_2}\circ R 
&\Rightarrow& x_0x_2\in I_2,\\
x_2\cdot \frac{\paa F^{alg}}{\paa x_2}\circ R \ \&\  
y\textup{-relations}
&\Rightarrow& -2y_2+x_2s_4+x_1x_2s_7\in I_2.
\end{eqnarray*}

\begin{eqnarray*}
&&x_2\cdot \frac{\paa F^{alg}}{\paa x_1}\circ R
\ \&\  \frac{\paa F^{alg}}{\paa x_2}\circ R
\ \&\  x_0\cdot \frac{\paa F^{alg}}{\paa x_2}\circ R\\
&\Rightarrow& 
-(x_1+s_5)(s_4+x_1s_7) + 3x_1^2x_2 + x_2s_3 +x_2^2s_7\in I_2.
\end{eqnarray*}
\begin{eqnarray*}
y\textup{-relation } x_0(x_1+s_5)-\kappa y_0 
&\Rightarrow& x_0(x_1+s_5)\in I_2.
\end{eqnarray*}

This gives the following $n+2$ elements of $I_2$.
The first three elements express $y_0,y_1$ and $y_2$
in terms of $(x,s)$, the last element is calculated from
these three elements and from \eqref{9.19}.
\begin{eqnarray}
&&-2y_0+x_0(x_0+s_6),
\ -2y_2+x_2(s_4+x_1s_7),\nonumber\\
&&-y_1-y_2+(2x_0+s_6)y_0+x_0s_2,\ x_3,...,\ x_n,
\ x_0(x_1+s_5), \ x_0x_2,\nonumber\\ 
&&(x_1+s_5)^2+x_2^2, \ 3x_1^2x_2-(x_1+s_5)(s_4+x_1s_7) +x_2s_3+x_2^2s_7,\nonumber\\
&&-\frac{1}{2}x_0(x_0+s_6)(3x_0+s_6) -2x_0s_2 
+3x_1^2(x_1+s_5)\nonumber\\
&&\hspace*{1cm} +x_2(s_4+x_1s_7)
+(x_1+s_5)s_3 +(x_1+s_5)x_2s_7.\label{9.20}
\end{eqnarray}
Restriction to $s=0$ gives the following $n+2$ elements of $I_3$.
\begin{eqnarray}
&& -2y_0+x_0^2,\ y_2,\ -y_1+2x_0y_0, \ x_3,...,\ x_n,
\ x_0x_1,\ x_0x_2,
\nonumber \\
&& x_1^2+x_2^2,\ x_1^2x_2,\ -x_0^3+2x_1^3.\label{9.21}
\end{eqnarray}
Therefore the monomials 
\begin{eqnarray*}
1,\ x_0,\ x_1,\ x_2,\ x_0^2,\ x_1^2,\ x_1x_2,\ x_0^3
\end{eqnarray*}
generate the quotient $\C[x,y]/I_3$. 
As this quotient cannot have dimension less than 8,
it has dimension 8, the elements in \eqref{9.21}
generate $I_3$, and the elements in \eqref{9.20}
generate $I_2$. 

\medskip
{\bf The case $\www E_7$:} 
\begin{eqnarray*}
\frac{\paa F^{alg}}{\paa x_0}\circ R &=& 
3 x_0^2x_1 -2(\kappa^2+1)\kappa^{-1}x_0x_1^2+x_1^3 \\
&& +\kappa s_2+ 2\kappa x_0s_4 +x_1s_5 +2x_0x_1s_7 +x_1^2s_8\\
&\stackrel{y\textup{-relation}}{\equiv} & 
3\kappa x_0y -2(\kappa^2+1)x_1y +x_1^3\\
&& + \kappa s_2 +2\kappa x_0s_4
+x_1s_5 +2\kappa ys_7 +x_1^2s_8.
\end{eqnarray*}
\begin{eqnarray*}
\frac{\paa F^{alg}}{\paa x_0}\circ R \ \&\  y\textup{-relation}
&\Rightarrow& -2x_1y+x_1^3+x_1s_5+x_1^2s_8\in I_2,\\
x_0\cdot \frac{\paa F^{alg}}{\paa x_0}\circ R \ \&\  
y\textup{-relation}
&\Rightarrow& 3x_0^2y -2y^2 +x_1^2y +x_0s_2 
+2x_0^2s_4 
\\&& +ys_5 +2x_0ys_7 +x_1ys_8\in I_2.
\end{eqnarray*}

\begin{eqnarray*}
\frac{\paa F^{alg}}{\paa x_1}\circ R &=& 
\kappa^{-1} x_0^3 -2(\kappa^2+1)\kappa^{-2}x_0^2x_1 
+3\kappa^{-1}x_0x_1^2 \\
&& +s_3 +\kappa^{-1}x_0s_5 +2x_1s_6 +\kappa^{-1}x_0^2s_7
 +2\kappa^{-1}x_0x_1s_8\\
&\stackrel{y\textup{-relation}}{\equiv} & 
\kappa^{-1}x_0^3 -2(\kappa^2+1)\kappa^{-1}x_0y +3x_1y \\
&& +s_3 +\kappa^{-1}x_0s_5 +2x_1s_6 +\kappa^{-1}x_0^2s_7 +2ys_8.
\end{eqnarray*}
\begin{eqnarray*}
\frac{\paa F^{alg}}{\paa x_1}\circ R \ \&\  y\textup{-relation}
&\Rightarrow& x_0^3-2x_0y+x_0s_5+x_0^2s_7\in I_2,\\
x_1\cdot \frac{\paa F^{alg}}{\paa x_1}\circ R \ \&\  
y\textup{-relation}
&\Rightarrow& x_0^2y-2y^2+3x_1^2y+x_1s_3 +ys_5\\
&& +2x_1^2s_6 +x_0ys_7 +2x_1ys_8\in I_2.
\end{eqnarray*}

\begin{eqnarray*}
y\textup{-relation } x_0x_1-\kappa y_0 
&\Rightarrow& x_0x_1\in I_2.
\end{eqnarray*}

This gives the following $n+4$ elements of $I_2$.
\begin{eqnarray}
&& x_2,...,\ x_n,\ x_0x_1,\ -2x_0y+x_0^3+x_0s_5+x_0^2s_7,
\nonumber\\
&& -2x_1y+x_1^3+x_1s_5+x_1^2s_8, \nonumber\\
&& -4y^2+(4x_0^2+4x_1^2+3x_0s_7+3x_1s_8+2s_5)y\nonumber\\
&& \hspace*{2cm}+x_0s_2+x_1s_3+2x_0^2s_4+2x_1^2s_6,\label{9.22}\\
&& (2x_0^2-2x_1^2+x_0s_7-x_1s_8)y +x_0s_2-x_1s_3+2x_0^2s_4
-2x_1^2s_6.\nonumber
\end{eqnarray}
Restriction to $s=0$ gives the following $n+4$ elements of $I_3$.
\begin{eqnarray}
x_2,...,\ x_n,\ x_0x_1,\ -2x_0y+x_0^3,\ -2x_1y+x_1^3,
\nonumber\\ 
-y^2+(x_0^2+x_1^2)y,\ (x_0^2-x_1^2)y.\label{9.23}
\end{eqnarray}
Therefore the monomials 
\begin{eqnarray*}
1,\ x_0,\ x_1,\ x_0^2,\ x_1^2,\ y, \ x_0^3,\ x_1^3,\ x_0^4
\end{eqnarray*}
generate the quotient $\C[x,y]/I_3$. 
As this quotient cannot have dimension less than 9,
it has dimension 9, the elements in \eqref{9.23}
generate $I_3$, and the elements in \eqref{9.22}
generate $I_2$.

\medskip
{\bf The case $\www E_8$:} 
\begin{eqnarray*}
\frac{\paa F^{alg}}{\paa x_0}\circ R &=& 
4 x_0^3x_1 -2(\kappa^3+1)\kappa^{-1}x_0x_1^2
+\kappa s_2+ 2\kappa x_0s_3\\
&& +3\kappa x_0^2s_5 
+ x_1s_6 +2x_0x_1s_7 +\kappa^{-1}x_1^2s_9\\
&\stackrel{y\textup{-relation}}{\equiv} & 
4\kappa x_0^2y +2\kappa x_0ys_9 + \kappa ys_9^2 
+\frac{1}{2}x_1s_9^3 -2\kappa^2 x_0x_1^2\\
&& -2x_1y +\kappa s_2 +2\kappa x_0s_3
+3\kappa x_0^2 s_5 +x_1s_6\\
&&  +2\kappa ys_7 +x_1s_7s_9.
\end{eqnarray*}
\begin{eqnarray*}
\frac{\paa F^{alg}}{\paa x_0}\circ R &\& &  
y\textup{-relation}\\
&\Rightarrow& \frac{1}{2}x_1s_9^3-2x_1y +x_1s_6+x_1s_7s_9
\in I_2,\\
(x_0-\frac{1}{2}s_9) 
\cdot \frac{\paa F^{alg}}{\paa x_0}\circ R &\& &  
y\textup{-relation}\\
&\Rightarrow& 4x_0^3y -2y^2  +(x_0-\frac{1}{2}s_9)
(s_2+2x_0s_3+3x_0^2s_5)\\
&& +ys_6 +2x_0ys_7\in I_2.
\end{eqnarray*}

\begin{eqnarray*}
\frac{\paa F^{alg}}{\paa x_1}\circ R &=& 
\kappa^{-1} x_0^4 -2(\kappa^3+1)\kappa^{-2}x_0^2x_1 
+3x_1^2 \\
&& +(s_4-\frac{1}{2}\kappa^{-1}s_6s_9 
-\frac{1}{4}\kappa^{-1}s_7s_9^2-\frac{1}{16}\kappa^{-1}s_9^4) \\
&& +\kappa^{-1}x_0s_6 +\kappa^{-1}x_0^2s_7 
+2x_1(s_8-\frac{1}{4}\kappa^{-2}s_9^2) \\
&& +2\kappa^{-2}x_0x_1s_9\\
&\stackrel{y\textup{-relation}}{\equiv} & 
\kappa^{-1}x_0^4 -2\kappa x_0^2x_1 -2\kappa^{-1}x_0y +3x_1^2 \\
&& (s_4-\frac{1}{2}\kappa^{-1}s_6s_9 
-\frac{1}{4}\kappa^{-1}s_7s_9^2-\frac{1}{16}\kappa^{-1}s_9^4)\\
&& +\kappa^{-1}x_0s_6 +\kappa^{-1}x_0^2s_7 
+2x_1s_8 +\kappa^{-1}ys_9.
\end{eqnarray*}
\begin{eqnarray*}
\frac{\paa F^{alg}}{\paa x_1}\circ R \ \&\  y\textup{-relation}
&\Rightarrow& x_0^4-2x_0y+ys_9-\frac{1}{2}s_6s_9 
-\frac{1}{4}s_7s_9^2\\
&&  -\frac{1}{16}s_9^4 +x_0s_6+x_0^2s_7\in I_2,\\
x_1\cdot \frac{\paa F^{alg}}{\paa x_1}\circ R \ \&\  
y\textup{-relation}
&\Rightarrow& (x_0^3 +x_0^2\frac{1}{2}s_9 +x_0\frac{1}{4}s_9^2 +\frac{1}{8}s_9^3)y \\
-2y^2+3x_1^3+x_1s_4 +ys_6
&& +(x_0+\frac{1}{2}s_9) ys_7+2x_1^2s_8\in I_2.
\end{eqnarray*}

\begin{eqnarray*}
y\textup{-relation } (x_0-\frac{1}{2}s_9)x_1-\kappa y_0 
&\Rightarrow& (x_0-\frac{1}{2}s_9)x_1\in I_2.
\end{eqnarray*}

This gives the following $n+4$ elements of $I_2$.
\begin{eqnarray}
&& x_2,...,\ x_n,\ (x_0-\frac{1}{2}s_9)x_1,\ 
x_1(-2y+s_6+s_7s_9+\frac{1}{2}s_9^3),\nonumber\\
&& y(-2y+s_6+2x_0s_7+4x_0^3) 
+(x_0-\frac{1}{2}s_9)(s_2+2x_0s_3+3x_0^2s_5),\nonumber\\
&& (x_0-\frac{1}{2}s_9)(-2y+s_6+(x_0+\frac{1}{2}s_9)s_7
+(x_0^3 +x_0^2\frac{1}{2}s_9 +x_0\frac{1}{4}s_9^2 
+\frac{1}{8}s_9^3)),\nonumber\\
&& y(-2y+s_6+(x_0+\frac{1}{2}s_9)s_7
+(x_0^3 +x_0^2\frac{1}{2}s_9 +x_0\frac{1}{4}s_9^2 
+\frac{1}{8}s_9^3)) \nonumber \\
&&\hspace*{2cm} + x_1(3x_1^2+s_4+2x_1s_8).\label{9.24}
\end{eqnarray}
Restriction to $s=0$ gives the following $n+4$ elements of $I_3$.
\begin{eqnarray}
x_2,...,\ x_n,\ x_0x_1,\ x_1y, \ -y^2+2x_0^3y, \nonumber\\ 
-2x_0y+x_0^4,\ 
-2y^2+ x_0^3y+3x_1^3.\label{9.25}
\end{eqnarray}
Therefore the monomials 
\begin{eqnarray*}
1,\ x_0,\ x_0^2,\ x_1,\ x_0^3,\ y, \ x_0^4,\ x_1^2,\ x_0^5,
\ x_0^6
\end{eqnarray*}
generate the quotient $\C[x,y]/I_3$. 
As this quotient cannot have dimension less than 10,
it has dimension 10, the elements in \eqref{9.25}
generate $I_3$, and the elements in \eqref{9.24}
generate $I_2$.

\medskip
(e) The critical space $C^{alg}$ of $F^{alg}$ is 
everywhere smooth as $F^{alg}$ is everywhere locally a
universal unfolding.
But the closure $\oooo{R^{-1}(C^{alg})}\subset Y$
is not smooth above $\C^{\mu-1}\times\{0\}$.
Denote by 
$$C^0:= \oooo{R^{-1}(C^{alg})}\cap 
\pr_\mu^{-1}(\C^{\mu-1}\times\{0\})$$
the restriction of it which lies above
$\C^{\mu-1}\times\{0\}$. 
It will turn out below that
$$C^{0,red}:= 
(C^0\textup{ with the reduced complex structure})$$
is in each of the three cases a union of four smooth
components,
$$C^{0,red}= C^{0,red}_1\cup C^{0,red}_2\cup C^{0,red}_3
\cup C^{0,red}_4.$$
The first three components appear in $C^0$ with their reduced
structure, $C^{0,red}_4$ appears in $C^0$ with a nonreduced
structure with multiplicity two.
The following table \eqref{9.27} 
collects facts, which will be proved
below in a case-by-case discussion.
{\it A-to-\eqref{9.26}} 
means the answer to the following question
\eqref{9.26}.
\begin{eqnarray}\label{9.26}
&& \textup{Is }C^{0,red}_i
\textup{ isomorphic to the critical space of a suitable}\\
&& \textup{unfolding which is obtained by a 
restriction of }F^{alg}\circ R\textup{?}\nonumber \\
&& \textup{And if yes, what is the type of the function
which is unfolded?}\nonumber
\end{eqnarray} 
\begin{eqnarray}\label{9.27}
\begin{array}{l|l|l|l|l|l}
 & & C^{0,red}_1 & C^{0,red}_2 & C^{0,red}_3 & C^{0,red}_4 \\
 \hline 
\www E_6 & \deg\pr_\mu|_{C^{0,red}_i} & 2 & 2 & 2 & 1 \\
 & \textup{A-to-\eqref{9.26}} & \textup{no} &
 \textup{yes},A_2 & \textup{yes},A_2 & \textup{no} \\ \hline 
\www E_7 & \deg\pr_\mu|_{C^{0,red}_i} & 3 & 3 & 1 & 1 \\
 & \textup{A-to-\eqref{9.26}} & \textup{yes},A_3 
 & \textup{yes},A_3 & \textup{no} & \textup{no} \\ \hline 
\www E_8 & \deg\pr_\mu|_{C^{0,red}_i} & 5 & 2 & 1 & 1 \\
 & \textup{A-to-\eqref{9.26}} & \textup{yes},A_5 
 & \textup{yes},A_2 & \textup{no} & \textup{no} 
\end{array}
\end{eqnarray}

The Lyashko-Looijenga map $LL^{alg}\circ\rho$ maps
$(s,\kappa)\in\C^{\mu-1}\times\Delta^*$ to the tuple
of the symmetric polynomials (with suitable signs)
in the values of $F^{alg}\circ R$ on
$R^{-1}(C^{alg})\cap\C^{n+1}\times\{(s,\kappa)\}$.

Because of (c), $F^{alg}\circ R$ extends to $\kappa=0$.
Because of (d), $R^{-1}(C^{alg})$ extends to
$\kappa=0$. Therefore the map $LL^{alg}\circ\rho$
extends to a holomorphic map
$(LL^{alg}\circ\rho)^{ext}:\C^{\mu-1}\times\Delta
\to M_{LL}^{(\mu)}$.
The table \eqref{9.27} shows that in each of the three cases
\begin{eqnarray}\label{9.28}
\sum_{i=1}^4 \deg\pr_\mu|_{C^{0,red}_i}=\mu-1.
\end{eqnarray}
Therefore $(LL^{alg}\circ\rho)^{ext}$ maps 
$\C^{\mu-1}\times\{0\}$ to $D_{LL}^{(\mu)}$.

It rests to show that the map in \eqref{9.13} is finite
and flat onto its image. Above $\C^{\mu-1}\times\Delta^*$
this holds. Therefore in order to prove it, it
rests to show 
\begin{eqnarray}\label{9.29}
(LL^{alg})^{ext}(s,0)=0&\Rightarrow & s=0.
\end{eqnarray}
This will be shown below in the case-by-case discussion.

\medskip
{\bf The case $\www E_6$:} 
\eqref{9.20} shows that $C^{0,red}$ has the following four
components $C^{0,red}_i$, $i\in\{1,2,3,4\}$.
The components are given in terms of functions which 
vanish on them.
In each case, they contain the first three functions
in \eqref{9.20} which express $y_0,y_1$ and $y_2$
in terms of $(x,s)$.
Of course, $x_3,...,x_n$ vanish on all four components.
\begin{eqnarray}
C^{0,red}_1:&& 
\textup{ the first three functions in \eqref{9.20}},
\ x_1+s_5,\ x_2,\nonumber\\
&& (x_0+s_6)(3x_0+s_6)+4s_2\  
(\textup{and generically }x_0\neq 0).\label{9.30} \\
C^{0,red}_2&&\textup{and }C^{0,red}_3: 
\textup{ the first three functions in \eqref{9.20}},\nonumber\\
&&x_0,\ x_2-\varepsilon\cdot i\cdot (x_1+s_5)\textup{ with}
\nonumber\\ 
&& \varepsilon=1\textup{ for }C^{0,red}_2\textup{ and }
\varepsilon=-1\textup{ for }C^{0,red}_3,\nonumber\\
&& 3x_1^2+\varepsilon i (2x_1s_7+s_4+s_5s_7)+s_3\nonumber \\
&& (\textup{and generically }x_1+s_5\neq 0,x_2\neq 0).
\label{9.31} \\
C^{0,red}_4:&& 
\textup{ the first three functions in \eqref{9.20}},\nonumber\\
&& x_0,\ x_1+s_5,\ x_2.\label{9.32}
\end{eqnarray}

Obviously, each $C^{0,red}_i$ is smooth, and 
$\deg\pr_\mu|_{C^{0,red}_i}$ is as claimed in table \eqref{9.27}.

It rests to prove \eqref{9.29}.
The restriction of $(F^{alg}\circ R)^{ext}$
(which was calculated in the proof of (c))
to $C^{0,red}_i$ is as follows:
\begin{eqnarray}\label{9.33}
(F^{alg}\circ R)^{ext}|_{C^{0,red}_1}&=& 
-\frac{1}{2}(x_0+s_6)x_0^2 + s_1 -s_3s_5-s_5^3,\\
(F^{alg}\circ R)^{ext}|_{C^{0,red}_j}&=&
x_1^3 +s_1+x_1s_3+\varepsilon i (x_1+s_5)s_4\nonumber \\
&& +\varepsilon i x_1(x_1+s_5)s_7\textup{ for }j\in\{2,3\},
\label{9.34} \\ 
(F^{alg}\circ R)^{ext}|_{C^{0,red}_4}&=& s_1-s_3s_5-s_5^3.
\label{9.35}
\end{eqnarray}

Consider a parameter $s\in \C^7$ with 
$(LL^{alg}\circ R)^{ext}(s,0)=0$.
We want to show $s=0$.
\eqref{9.35} gives $s_1-s_3s_5-s_5^3=0$. 
\eqref{9.33} and \eqref{9.30} give the existence of 
a number $x_0\in\C$ with 
$(x_0+s_6)(3x_0+s_6)+4s_2=0$ and $(x_0+s_6)x_0=0$.
The first quadratic polynomial has a double zero
if and only if $12s_2-s_6^2=0$,
and then the double zero is $x_0=-\frac{2}{3}s_6$.
It is a zero of $(x_0+s_6)x_0$ only if $s_6=0$,
and then $s_2=0$. 
In the case $12s_2-s_6^2\neq 0$, we must have
$$(x_0+s_6)(3x_0+s_6)+4s_2=3(x_0+s_6)x_0,
\quad\textup{thus again } s_6=0,s_2=0.$$
So, \eqref{9.33} gives in any case $s_6=0$ and $s_2=0$.

Now consider $j\in\{2,3\}$ and \eqref{9.34}.
It motivates the definition of the unfolding
\begin{eqnarray*}
G_j(x_1,s_1,s_3,s_4,s_5,s_7):=
x_1^3+s_1+x_1s_3+\varepsilon i (x_1+s_5)s_4 + 
\varepsilon i x_1(x_1+s_5)s_7
\end{eqnarray*}
in the variable $x_1$ and with parameters 
$s_1,s_3,s_4,s_5,s_7$ of the $A_2$-singularity $x_1^3$.
The derivative $\frac{\paa G_j}{\paa x_1}$
is in the ideal which defines $C^{0,red}_j$,
so 
$$\textup{Crit}(G_j)\cong C^{0,red}_j\textup{ and }
G_j|_{\textup{Crit}(G_j)} \cong 
(F^{alg}\circ R)^{ext}|_{C^{0,red}_j}.$$
Denote the Lyashko-Looijenga map of $G_j$ by $LL_{G_j}$.
Then $LL_{G_j}(s_1,s_3,s_4,s_5,s_7)=0$.

The unfolding $G_j$ is induced by the universal unfolding
$$G_{A_2}(z,t_1,t_2)=z^3+t_1+zt_2$$ 
via the morphism $(\Phi^{(j)},\varphi^{(j)})$ with
$G_{A_2}\circ\Phi^{(j)}=G_j$ and 
\begin{eqnarray*}
z=\Phi^{(j)}_1(x_1,s)&=& 
x_1+\frac{1}{3}\varepsilon i s_7,\\
t_1=\varphi^{(j)}_1(s)&=&
s_1+\varepsilon i s_4s_5+\frac{1}{27}\varepsilon i s_7^3 \\
&&-\frac{1}{3}\varepsilon i (s_3+\varepsilon i s_4+
\varepsilon i s_5s_7 +\frac{1}{3}s_7^2)s_7,\\
t_2=\varphi^{(j)}_2(s)&=& 
s_3+\varepsilon i s_4 +\varepsilon i s_5s_7 +\frac{1}{3}s_7^2.
\end{eqnarray*}
Then $LL_{G_j} = LL_{A_2}\circ \varphi^{(j)}$
where $LL_{A_2}$ is the Lyashko-Looijenga map of the 
universal unfolding  $G_{A_2}$.
The map $LL_{A_2}$ is a finite branched covering and has
value 0 only at 0. Therefore
$\varphi^{(j)}_2(s)=\varphi^{(j)}_1(s)=0$.
Now $\varphi^{(1)}_2\pm \varphi^{(2)}_2=0$ and 
$\varphi^{(1)}_1\pm \varphi^{(2)}_1=0$ give
$$ s_3+\frac{1}{3}s_7^2=0,\ s_4+s_5s_7=0,\ s_1=0,\ 
s_4s_5+\frac{1}{27}s_7^3=0.$$
Together with 
$$s_1-s_3s_5-s_5^3=0 \textup{ and }s_6=0,\ s_2=0,$$
this gives $s=0$. Now \eqref{9.29} is proved in the
case $\www E_6$.

\medskip
{\bf The case $\www E_7$:}
\eqref{9.22} shows that $C^{0,red}$ has the following four
components $C^{0,red}_i$, $i\in\{1,2,3,4\}$.
The components are given in terms of functions which 
vanish on them.
Of course, $x_2,...,x_n$ vanish on all four components.
\begin{eqnarray}
C^{0,red}_1:&& 
-2y+x_0^2+s_5+x_0s_7,\ x_1,\nonumber\\ 
&&(2x_0+s_7)y+s_2+2x_0s_4\nonumber \\
&& (\textup{and generically }x_0\neq 0).\label{9.36} \\
C^{0,red}_2:&&
-2y+x_1^2+s_5+x_1s_8,\ x_0,\nonumber\\ 
&&(2x_1+s_8)y+s_3+2x_1s_6\nonumber \\
&& (\textup{and generically }x_1\neq 0).
\label{9.37} \\
C^{0,red}_3:&& y,\ x_0,\ x_1.\label{9.38}\\
C^{0,red}_4:&& y-\frac{1}{2}s_5,\ x_0,\ x_1.\label{9.39}
\end{eqnarray}

Obviously, each $C^{0,red}_i$ is smooth, and 
$\deg\pr_\mu|_{C^{0,red}_i}$ is as claimed in table \eqref{9.27}.

It rests to prove \eqref{9.29}.
The restriction of $(F^{alg}\circ R)^{ext}$
to $C^{0,red}_i$ is as follows:
\begin{eqnarray}
(F^{alg}\circ R)^{ext}|_{C^{0,red}_1}&=& 
x_0^2y-y^2+s_1+x_0s_2+x_0^2s_4\nonumber\\
&&+ys_5+x_0ys_7,\label{9.40}\\
(F^{alg}\circ R)^{ext}|_{C^{0,red}_2}&=&
x_1^2y-y^2+s_1+x_1s_3+x_1^2s_6\nonumber\\
&&+ys_5+x_1ys_8,\label{9.41} \\ 
(F^{alg}\circ R)^{ext}|_{C^{0,red}_3}&=& s_1,
\label{9.42}\\
(F^{alg}\circ R)^{ext}|_{C^{0,red}_4}&=& s_1+\frac{1}{4}s_5^2.
\label{9.43}
\end{eqnarray}

Consider a parameter $s\in \C^8$ with 
$(LL^{alg}\circ R)^{ext}(s,0)=0$.
We want to show $s=0$.
\eqref{9.42} and \eqref{9.43} give $s_1=s_5=0$.
This and \eqref{9.40} motivate the definition of the unfolding
\begin{eqnarray*}
G(x_0,y,s_2,s_4,s_7):=
x_0^2y-y^2+x_0s_2+ x_0^2s_4+x_0ys_7
\end{eqnarray*}
of the $A_3$-singularity in the parameters 
$s_2,s_4,s_7$.
The derivatives $\frac{\paa G}{\paa x_0}$
and $\frac{\paa G}{\paa y}$
are in the ideal which defines $C^{0,red}_1|_{s_1=s_5=0}$,
so 
$$\textup{Crit}(G)\cong C^{0,red}|_{s_1=s_5=0}\textup{ and }
G|_{\textup{Crit}(G)} \cong 
(F^{alg}\circ R)^{ext}|_{C^{0,red}_1|_{s_1=s_5=0}}.$$
Denote the Lyashko-Looijenga map of $G$ by $LL_{G}$.
Then $LL_{G}(s_2,s_4,s_7)=0$.

The unfolding $G$ is induced by the universal unfolding
$$G_{A_3}(z,y_1,t_1,t_2,t_3)=z^2y_1-y_1^2+t_1+zt_2+xz^2t_3$$
via the morphism $(\Phi,\varphi)$ with
$G_{A_3}\circ\Phi=G$ and 
\begin{eqnarray*}
z=\Phi_1(x_0,y,s)&=& 
x_0+\frac{1}{2}s_7,\\
y_1=\Phi_1(x_0,y,s)&=& 
y+\frac{1}{8}s_7^2 ,\\
t_1=\varphi_1(s)&=&
-\frac{1}{2}s_2s_7 +\frac{1}{4}s_4s_7^2+\frac{1}{64}s_7^4, \\
t_2=\varphi_2(s)&=& 
s_2-s_4s_7,\\
t_3=\varphi_3(s)&=& 
s_4-\frac{1}{8}s_7^2.
\end{eqnarray*}
Then $LL_{G} = LL_{A_3}\circ \varphi$
where $LL_{A_3}$ is the Lyashko-Looijenga map of the 
universal unfolding  $G_{A_3}$.
The map $LL_{A_3}$ is a finite branched covering and has
value 0 only at 0. Therefore
$0=\varphi_1(s)=\varphi_2(s)=\varphi_3(s)$.
This gives $s_2=s_4=s_7=0$. 

\eqref{9.36}\ \& \eqref{9.37} and \eqref{9.40}\ \&\ \eqref{9.41}
are symmetric with respect to
$$(x_0,y,s_1,s_2,s_4,s_5,s_7)\longleftrightarrow
(x_1,y,s_1,s_3,s_6,s_5,s_8).$$
Therefore also $s_3=s_6=s_8=0$. 
This gives $s=0$. Now \eqref{9.29} is proved in the
case $\www E_7$.

\medskip
{\bf The case $\www E_8$:}
\eqref{9.24} shows that $C^{0,red}$ has the following four
components $C^{0,red}_i$, $i\in\{1,2,3,4\}$.
The components are given in terms of functions which 
vanish on them.
Of course, $x_2,...,x_n$ vanish on all four components.
\begin{eqnarray}
C^{0,red}_1:&& 
-2y+(x_0^3+x_0^2\frac{1}{2}s_9+x_0\frac{1}{4}s_9^2 
+\frac{1}{8}s_9^3)+s_6+(x_0+\frac{1}{2}s_9)s_7,\nonumber\\ 
&& x_1,\ y(s_7+3x_0^2+x_0s_9+\frac{1}{4}s_9^2)+s_2+2x_0s_3
+3x_0^2s_5\nonumber \\
&& (\textup{and generically }x_0-\frac{1}{2}s_9\neq 0).
\label{9.44} \\
C^{0,red}_2:&&
x_0-\frac{1}{2}s_9,\ -2y+s_6+s_7s_9+\frac{1}{2}s_9^3,\nonumber\\ 
&& 3x_1^2+s_4+2x_1s_8\  (\textup{and generically }x_1\neq 0).
\label{9.45} \\
C^{0,red}_3:&& x_0-\frac{1}{2}s_9,\ x_1,\ y.\label{9.46}\\
C^{0,red}_4:&& x_0-\frac{1}{2}s_9,\ x_1,\ 
-2y+s_6+s_7s_9+\frac{1}{2}s_9^3.\label{9.47}
\end{eqnarray}

Obviously, each $C^{0,red}_i$ is smooth, and 
$\deg\pr_\mu|_{C^{0,red}_i}$ is as claimed in table \eqref{9.27}.

It rests to prove \eqref{9.29}.
The restriction of $(F^{alg}\circ R)^{ext}$
to $C^{0,red}_i$ is as follows:
\begin{eqnarray}\label{9.48}
(F^{alg}\circ R)^{ext}|_{C^{0,red}_1}&=& 
y\bigl[(x_0^3+x_0^2\frac{1}{2}s_9+x_0\frac{1}{4}s_9^2 
+\frac{1}{8}s_9^3)+s_6\\
&&+(x_0+\frac{1}{2}s_9)s_7\bigr]
-y^2+s_1+x_0s_2+x_0^2s_3+x_0^3s_5,\nonumber\\
(F^{alg}\circ R)^{ext}|_{C^{0,red}_2}&=&
y(s_6+s_7s_9+\frac{1}{2}s_9^3)-y^2+x_1^3\label{9.49}\\
&&+s_1+\frac{1}{2}s_2s_9+\frac{1}{4}s_3s_9^2+x_1s_4
+\frac{1}{8}s_5s_9^3+x_1^2s_8,\nonumber \\ 
(F^{alg}\circ R)^{ext}|_{C^{0,red}_3}&=& 
s_1+\frac{1}{2}s_2s_9+\frac{1}{4}s_3s_9^2+\frac{1}{8}s_5s_9^3,
\label{9.50}\\
(F^{alg}\circ R)^{ext}|_{C^{0,red}_4}&=&
\frac{1}{4}(s_6+s_7s_9+\frac{1}{2}s_9^3)^2 \label{9.51}\\
&&+(s_1+\frac{1}{2}s_2s_9+\frac{1}{4}s_3s_9^2
+\frac{1}{8}s_5s_9^3).\nonumber
\end{eqnarray}

Consider a parameter $s\in \C^9$ with 
$(LL^{alg}\circ R)^{ext}(s,0)=0$.
We want to show $s=0$.
\eqref{9.50} and \eqref{9.51} give
\begin{eqnarray}\label{9.52}
0=s_1+\frac{1}{2}s_2s_9+\frac{1}{4}s_3s_9^2+\frac{1}{8}s_5s_9^3,
\quad 0= s_6+s_7s_9+\frac{1}{2}s_9^3.
\end{eqnarray}

This and \eqref{9.48} motivate the definition of the unfolding
\begin{eqnarray*}
&&G_1(x_0,y,s_2,s_3,s_5,s_7,s_9)\\
&:=& y\bigl[(x_0^3+x_0^2\frac{1}{2}s_9+x_0\frac{1}{4}s_9^2  
-\frac{3}{8}s_9^3)+(x_0-\frac{1}{2}s_9)s_7\bigr]\\
&&-y^2 -(\frac{1}{2}s_2s_9+\frac{1}{4}s_3s_9^2
+\frac{1}{8}s_5s_9^3)+x_0s_2+x_0^2s_3+x_0^3s_5
\end{eqnarray*}
of the $A_5$-singularity $yx_0^3-y^2$ in the parameters 
$s_2,s_3,s_5,s_7,s_9$.
The derivatives $\frac{\paa G_1}{\paa x_0}$
and $\frac{\paa G_1}{\paa y}$
are in the ideal which defines 
$C^{0,red}_1|_{s\textup{ with }\eqref{9.52}}$,
so 
$$\textup{Crit}(G_1)\cong C^{0,red}|_{s\textup{ with }
\eqref{9.52}}\textup{ and }
G_1|_{\textup{Crit}(G_1)} \cong 
(F^{alg}\circ R)^{ext}|_{C^{0,red}_1|_{s\textup{ with }
\eqref{9.52}}}.$$
Denote the Lyashko-Looijenga map of $G_1$ by $LL_{G_1}$.
Then $LL_{G_1}(s_2,s_3,s_5,s_7,s_9)=0$.

The unfolding $G_1$ is induced by the universal unfolding
$$G_{A_5}=y_1z^3-y_1^2+t_1+zt_2+z^2t_3+z^3t_4+zy_1t_5$$ 
via the morphism $(\Phi^{(1)},\varphi^{(1)})$ with
$G_{A_5}\circ\Phi^{(1)}=G_1$ and 
\begin{eqnarray*}
z=\Phi_1^{(1)}(x_0,y,s)&=& 
x_0+\frac{1}{6}s_9,\\
y_1=\Phi_2^{(1)}(x_0,y,s)&=& 
y+\frac{1}{2}(\frac{2}{3}s_7s_9+\frac{11}{27}s_9^3),\\
t_1=\varphi_1^{(1)}(s)&=&
-\frac{1}{4}(\frac{2}{3}s_7s_9+\frac{11}{27}s_9^3)^2\\
&&+(-\frac{2}{3}s_2s_9-\frac{2}{9}s_3s_9^2
-\frac{7}{54}s_5s_9^3),\\
t_2=\varphi_2^{(1)}(s)&=& s_2-\frac{1}{3}s_3s_9+
\frac{1}{12}s_5s_9^2,\\
t_3=\varphi_3^{(1)}(s)&=& s_3-\frac{1}{2}s_5s_9, \\
t_4=\varphi_4^{(1)}(s)&=& s_5 ,\\
t_5=\varphi_5^{(1)}(s)&=& s_7+\frac{1}{6}s_9^2.
\end{eqnarray*}
Then $LL_{G} = LL_{A_5}\circ \varphi^{(1)}$
where $LL_{A_5}$ is the Lyashko-Looijenga map of the 
universal unfolding  $G_{A_5}$.
The map $LL_{A_5}$ is a finite branched covering and has
value 0 only at 0. Therefore
$$0=\varphi_1^{(1)}(s)=\varphi_2^{(1)}(s)=\varphi_3^{(1)}(s)
=\varphi_4^{(1)}(s)=\varphi_5^{(1)}(s).$$
This gives 
\begin{eqnarray}\label{9.53}
s_2=s_3=s_5=s_7=s_9=0,
\textup{ and with \eqref{9.52}}\ s_1=s_6=0.
\end{eqnarray}

This and \eqref{9.49} motivate the definition of the unfolding
\begin{eqnarray*}
G_2(x_1,y,s_4,s_8)&:=&
-y^2+x_1^3+x_1s_4+x_1^2s_8
\end{eqnarray*}
of the $A_2$-singularity $-y^2+x_1^3$ in the parameters 
$s_4$ and $s_8$.
The derivatives $\frac{\paa G_2}{\paa x_1}$
and $\frac{\paa G_2}{\paa y}$
are in the ideal which defines 
$C^{0,red}_2|_{s\textup{ with }\eqref{9.53}}$,
so 
$$\textup{Crit}(G_2)\cong C^{0,red}_2|_{s\textup{ with }
\eqref{9.53}}\textup{ and }
G_2|_{\textup{Crit}(G_2)} \cong 
(F^{alg}\circ R)^{ext}|_{C^{0,red}_2|_{s\textup{ with }
\eqref{9.53}}}.$$
Denote the Lyashko-Looijenga map of $G_2$ by $LL_{G_2}$.
Then $LL_{G_2}(s_4,s_8)=0$.

The unfolding $G_2$ is induced by the universal unfolding
$$G_{A_2}(z,y_1,t_1,t_2)=-y_1^2+z^3+t_1+zt_2$$ 
via the morphism $(\Phi^{(2)},\varphi^{(2)})$ with
$G_{A_2}\circ\Phi^{(2)}=G_2$ and 
\begin{eqnarray*}
z=\Phi_1^{(2)}(x_1y,s)&=& 
x_1+\frac{1}{3}s_8,\\
y_1=\Phi_2^{(2)}(x_0,y,s)&=& y,\\
t_1=\varphi_1^{(2)}(s)&=&
-\frac{1}{3}s_4s_8+\frac{2}{27}s_8^3,\\
t_2=\varphi_2^{(2)}(s)&=& s_4-\frac{1}{3}s_8^2.
\end{eqnarray*}
Then $LL_{G_2} = LL_{A_2}\circ \varphi^{(2)}$
where $LL_{A_2}$ is the Lyashko-Looijenga map of the 
universal unfolding  $G_{A_2}$.
The map $LL_{A_2}$ is a finite branched covering and has
value 0 only at 0. Therefore
$$0=\varphi_1^{(2)}(s)=\varphi_2^{(2)}(s),\quad\textup{so}\quad
0=s_4=s_8.$$

This gives $s=0$. Now \eqref{9.29} is proved in the
case $\www E_8$.
This finishes the proof of theorem \ref{t9.1}.
\hfill$\Box$

\section{Degree of the Lyashko-Looijenga map $LL^{alg}$
for the simple elliptic singularities}\label{s10}
\setcounter{equation}{0}

\noindent
This section is devoted to the proof of theorem \ref{t6.3}.
The main work has already been done in the sections
\ref{s9}, \ref{s5} and \ref{s8}.
The maps $\rho$ in theorem \ref{t9.1} tell how to glue
into $M^{alg}=\C^{\mu-1}\times (\C-\{0,1\})$ a fiber
above $\lambda=0$. This and the maps $\psi_2$ and $\psi_3$
in subsection \ref{s5.2} tell how to glue into $M^{alg}$
fibers above $\lambda=1$ and $\lambda=\infty$.
Corollary \ref{t8.6} together with the maps
$\psi_2,\psi_3$ and $\rho$ allows to calculate the
degree of $LL^{alg}$. 

The maps $\psi_2$ in \eqref{5.21},\eqref{5.27} and \eqref{5.33}
and the maps $\rho$ in \eqref{9.2}, \eqref{9.5} and \eqref{9.8}
contain the following fractional powers of $\lambda$,
\begin{eqnarray}\label{10.1}
\begin{array}{llll}
 & \www E_6 & \www E_7 & \www E_8 \\
\psi_2 & \lambda^{1/2} & \lambda^{1/4} & \lambda^{1/2} \\
\rho & \kappa=\lambda^{1/c}=\lambda^{1/3} & 
\kappa=\lambda^{1/c}=\lambda^{1/2} &
\kappa=\lambda^{1/c}=\lambda^{1/3}
\end{array}
\end{eqnarray}
Therefore we consider coverings of $\C-\{0,1\}$
and of $M^{alg}$ which are  of order $2c$ at respectively
above $\lambda\in\{0,1,\infty\}$. 
Denote by $\P^1(2c,2c,2c)$ the orbifold $\P^1$ with 
orbifold points $0,1$ and $\infty$ which all have multiplicity
$2c$. Because of $2-3(1-\frac{1}{2c})=-1+\frac{3}{2c}<0$
it is a hyperbolic orbifold, so a good orbifold.
By a classical theorem of Fox \cite[Theorem 2.5]{Sc83},
a finite orbifold covering $p_X:X\to\P^1(2c,2c,2c)$
with $X$ a manifold exists. 
It is a branched covering of order $2c$ at each preimage
of $0$, $1$ and $\infty$ and a covering everywhere else.
Denote $N^{alg}:= \C^{\mu-1}\times (X-p_X^{-1}(\{0,1,\infty\})$,
and denote by $p_{alg}:=\id|_{\C^{\mu-1}}\times 
p_X:N^{alg}\to M^{alg}$ 
the lift to $N^{alg}$ of the restricted map 
$p_X:X-p_X^{-1}(\{0,1,\infty\})\to \C-\{0,1,\infty\}$. 

The bundles $M^{alg}\to \C-\{0,1,\infty\}$ and
$N^{alg}\to X-p_X^{-1}(\{0,1,\infty\})$ are smooth cone bundles
with the weights 
$(a_1,a_2,..,a_{\mu-1})=
(\deg_{\bf w}t_{\mu-1},\deg_{\bf w}t_{\mu-2},...,
\deg_{\bf w}t_1)\cdot d$. Here 
 \begin{eqnarray}\label{10.2}
d:=3\textup{ for }\www E_6, \quad 
d:=4\textup{ for }\www E_7, \quad
d:=6\textup{ for }\www E_8,
\end{eqnarray}
is chosen so that all weights $d\cdot \deg_{\bf w}t_i$ are
integers.
Now $N^{alg}$ will be extended to a smooth cone bundle
$N^{orb}\to X$, i.e. fibers above the points in
$p_X^{-1}(\{0,1,\infty\})$ will be glued into $N^{alg}$.

Let $\delta_0:\Delta\to X$ be an isomorphism from
the unit disk $\Delta$ to a neighborhood of any point
in $p_X^{-1}(0)$ with $p_X\circ \delta_0(z)=z^{2c}$.
Glue $\C^{\mu-1}\times \Delta$ into $N^{alg}$ with the map
\begin{eqnarray}\label{10.3}
\C^{\mu-1}\times\Delta^*&\hookrightarrow& N^{alg},\\
(t',z)&\mapsto& ((\rho_1,...,\rho_{\mu-1})(t',z^2),
\delta_0(z)).\nonumber
\end{eqnarray}
Let $\delta_1:\Delta\to X$ be an isomorphism from
the unit disk $\Delta$ to a neighborhood of any point
in $p_X^{-1}(1)$ with $p_X\circ \delta_1(z)=1-z^{2c}$.
Glue $\C^{\mu-1}\times \Delta$ into $N^{alg}$ with the map
\begin{eqnarray}\label{10.4}
\C^{\mu-1}\times\Delta^*&\hookrightarrow& N^{alg},\\
(t',z)&\mapsto& (((\psi_3)_1,...(\psi_3)_{\mu-1})
((\rho_1,...,\rho_{\mu-1})(t',z^2),z^{2c}),\delta_1(z)).
\nonumber
\end{eqnarray}
Let $\delta_\infty:\Delta\to X$ be an isomorphism from
the unit disk $\Delta$ to a neighborhood of any point
in $p_X^{-1}(\infty)$ with $p_X\circ \delta_\infty(z)=z^{-2c}$.
Glue $\C^{\mu-1}\times \Delta$ into $N^{alg}$ with the map
\begin{eqnarray}\label{10.5}
\C^{\mu-1}\times\Delta^*&\hookrightarrow& N^{alg},\\
(t',z)&\mapsto& (((\psi_2)_1,...,(\psi_2)_{\mu-1})
((\rho_1,...,\rho_{\mu-1})(t',z^2),z^{2c}),\delta_\infty(z)).
\nonumber
\end{eqnarray}
This is a univalued map although $\psi_2$ contains
$\lambda^{1/2}$ (in the cases $\www E_6$ and $\www E_8$)
and $\lambda^{1/4}$ (in the case $\www E_7$),
by setting $\lambda^{1/2}\circ z^{2c}=z^c$ and 
$\lambda^{1/4}\circ z^4=z$.

The resulting manifold $N^{orb}$ is a smooth cone bundle above 
$X$ with weights 
$(a_1,a_2,..,a_{\mu-1})=
(\deg_{\bf w}t_{\mu-1},\deg_{\bf w}t_{\mu-2},...,
\deg_{\bf w}t_1)\cdot d$
because $M^{alg}$ and $N^{alg}$ are smooth cone bundles with
these weights and all involved maps are $\C^*$-equivariant
with respect to the natural $\C^*$-actions.
The covering group of the covering
$p_{alg}:N^{alg}\to M^{alg}$ extends to an automorphism group
of $N^{orb}$. The quotient of $N^{orb}$ by this group
is an orbibundle $M^{orb}$ above $\P^1$ which extends
$M^{alg}\to\C-\{0,1\}$. 
Let $p_{orb}:N^{orb}\to M^{orb}$ be the quotient map.

Recall the definition of $M^{orb}_0\subset M^{orb}$
in theorem \ref{t6.3}, and define 
$N^{orb}_0:=p_{orb}^{-1}(M^{orb}_0)$.
We claim that $LL^{alg}\circ p_{alg}:N^{alg}\to M_{LL}^{(\mu)}$
extends to a holomorphic map 
$LL^{orb}_N:N^{orb}\to M_{LL}^{(\mu)}$, 
that the restriction
\begin{eqnarray}\label{10.6}
LL^{orb}_N: N^{orb}-N^{orb}_0\to M_{LL}^{(\mu)}-M_{LL,0}^{(\mu)}
\end{eqnarray}
is a branched covering of a finite degree,
and that $LL^{orb}_N$ maps the fibers of $N^{orb}$
above the points in $p_X^{-1}(\{0,1,\infty\})$ to 
$D_{LL}^{(\mu)}$.

Near the fibers of $N^{orb}$ above the points in 
$p_X^{-1}(0)$, this follows from theorem \ref{t9.1} (e).
Near the fibers of $N^{orb}$ above the points in
$p_X^{-1}(\{1,\infty\})$, this follows again from 
theorem \ref{t9.1} (e) and from the fact that 
$\psi_3$ and $\psi_2$ are locally isomorphisms of 
F-manifolds with Euler fields and thus
\begin{eqnarray}\label{10.7}
LL^{alg}(t',\lambda) 
&=& LL^{alg}(\psi_2(t',\lambda))
=LL^{alg}(\psi_3(t',\lambda)) \\
\textup{for}&&
(t',\lambda)\in \C^{\mu-1}\times\Delta^*\subset M^{alg}
=\C^{\mu-1}\times(\C-\{0,1\}).\nonumber
\end{eqnarray}

$LL^{alg}$ inherits the good properties from 
$LL^{alg}\circ p_{alg}$.
It extends to a holomorphic map 
$LL^{orb}:M^{orb}\to M_{LL}^{(\mu)}$,
the restriction in \eqref{6.6} is a branched covering,
and $\pi_{orb}^{-1}(\{0,1,\infty\})$ is mapped to
$D_{LL}^{(\mu)}$. 

It rests to determine the degree
of $LL^{alg}$. Of course,
$\deg LL^{orb}_N=\deg LL^{alg}\cdot \deg p_{alg}$.

The tuple $(LL^{orb}_N,N^{orb},M_{LL}^{(\mu)})$ satisfies
almost the properties of the tuple
$(f,C_1,C_2)$ in corollary \ref{t8.6}, but not completely.

The affine group $G_{\A_1}=(\C,+)$ acts freely on $N^{orb}$
and $M_{LL}^{(\mu)}$ as follows, 
and $LL^{orb}_N$ is equivariant with respect to these actions.
We have to divide out these actions.
The action of $G_{\A_1}$ on $N^{orb}$ comes from the lift
to $N^{alg}$ and extension to $N^{orb}$ of the action 
on $M^{alg}$,
\begin{eqnarray}\label{10.8}
G_{\A_1}\times M^{alg}\to M^{alg},\quad 
(s,t',\lambda)\mapsto (t_1+s,t_2,...,t_{\mu-1},\lambda).
\end{eqnarray}
The action of $G_{\A_1}$ on $M_{LL}^{(\mu)}$ is given by
\begin{eqnarray}\label{10.9}
G_{\A_1}\times M_{LL}^{(\mu)}\to M_{LL}^{(\mu)},\quad 
(s,p(y))\mapsto p(y-s).
\end{eqnarray}
The quotient triple $(LL^{orb}_N,N^{orb},M_{LL}^{(\mu)})/G_{\A_1}$
satisfies the properties of the triple 
$(f,C_1,C_2)$ in corollary \ref{t8.6}.

$C_1:=N^{orb}/G_{\A_1}$ is a smooth cone bundle with
weights $(a_1,...,a_{\mu-2})=(\deg_{\bf w}t_{\mu-1},...,
\deg_{\bf w}t_2)\cdot d$ 
and basis $X_1:=X$ of dimension 1. 
$C_2:=M_{LL}^{(\mu)}/G_{\A_1}$ is a smooth cone bundle with
weights $(b_1,b_2,...,b_{\mu-1})=(2,3,...,\mu-1)\cdot d$
with basis $X_2=(\textup{a point})$.
And $f:=LL^{orb}_N/G_{\A_1}$ satisfies the properties in the
situation before proposition \ref{t8.5}. Therefore
by corollary \ref{t8.6} 
\begin{eqnarray}
\deg LL^{orb}_N &=& \deg f = 
\frac{b_1...b_{\mu-1}}{a_1...a_{\mu-2}}\cdot\left( 
-\sum_{k=a_1}^{a_{\mu-2}}\frac{1}{k}\cdot\deg C_{1,(k)}\right)
\nonumber \\
&=& \frac{2\cdot 3\cdot ...\cdot \mu}
{\prod_{i=2}^{\mu-1}\deg_{\bf w}t_i}\cdot \left( 
-\sum_{k=a_1}^{a_{\mu-2}} \frac{d}{k}\cdot \deg C_{1,(k)}\right),
\label{10.10} \\
\deg LL^{alg}&=& 
\frac{\mu!}{\prod_{i=2}^{\mu-1}\deg_{\bf w}t_i}\cdot\left( 
\sum_{k=a_1}^{a_{\mu-2}} \frac{d}{k}\cdot\left(
-\frac{\deg C_{1,(k)}}{\deg p_{alg}}\right)\right).
\label{10.11}
\end{eqnarray}
For the proof of formula \eqref{6.7} it rests to show
\begin{eqnarray}\label{10.12}
-\frac{\deg C_{1,(k)}}{\deg p_{alg}} = 
\frac{1}{2}\cdot |\{j\, |\, a_j=k\}|.
\end{eqnarray}

A basis of trivial global sections of the trivial smooth cone bundle
$\C^{\mu-1}\times X\supset N^{alg}$ and the glueing maps 
\eqref{10.3}, \eqref{10.4} and \eqref{10.5} give a global
meromorphic section in the determinant bundle 
$\deg C_{1,(k)}$ of each vector bundle $C_{1,(k)}$. 
The sum of the orders of zeros and poles of this section
is $\deg C_{1,(k)}$. In fact, we can read of
$-\deg C_{1,(k)}/\deg p_{alg}$ directly from the sum of the 
orders of $\lambda$ in those parts of 
$\rho,\psi_3$ and $\psi_2$ which correspond to $C_{1,(k)}$.
Here $\rho$ is used three times, $\psi_3$ and $\psi_2$ are
each used one times. The following tables collect
the relevant data from the formulas 
\eqref{9.2}, \eqref{9.5}, \eqref{9.8}, \eqref{5.24}, 
\eqref{5.30}, \eqref{5.36}, \eqref{5.21}, \eqref{5.27}
and \eqref{5.33}.

\bigskip
The case $\www E_6$:
\begin{eqnarray*}
\begin{array}{l|l|l|l}
k & \textup{involved }t_i& \rho:\textup{order of} 
& \psi_3:\textup{order of} \\ 
 & & \lambda\textup{ in \eqref{9.2}} & 
 \lambda\textup{ in \eqref{5.24}} \\ \hline
1=a_1=a_2=a_3 & t_7,t_6,t_5 & 0+\frac{1}{3}+1 & 0+0+0 
\\[1mm]
2=a_4=a_5=a_6 & t_4,t_3,t_2 & 0+0+\frac{2}{3} & 0+0+0
\end{array}
\end{eqnarray*}

\begin{eqnarray*}
\begin{array}{l|l|l}
k & \psi_2:\textup{order of }\lambda\textup{ in }\eqref{5.21} &
 -\deg C_{1,(k)}/\deg p_{alg} \\[1mm] \hline
1 & \frac{1}{2}-1-2 & 
3\cdot \frac{4}{3}+0-\frac{5}{2}=\frac{3}{2} \\[1mm]
2 & \frac{1}{2}+0-1 & 
3\cdot \frac{2}{3}+0-\frac{1}{2}=\frac{3}{2} 
\end{array}
\end{eqnarray*}

\bigskip
The case $\www E_7$:
\begin{eqnarray*}
\begin{array}{l|l|l|l}
k & \textup{involved }t_i& \rho:\textup{order of} 
& \psi_3:\textup{order of} \\ 
 & & \lambda\textup{ in \eqref{9.5}} & 
 \lambda\textup{ in \eqref{5.30}} \\ \hline
1=a_1=a_2 & t_8,t_7 & 0+\frac{1}{2} & 
1\textup{ (see \eqref{10.13}})
\\[1mm]
2=a_3=a_4=a_5 & t_6,t_5,t_4 & 0+0+1 & 0+0+0 \\
3=a_6=a_7 & t_3,t_2 & 0+\frac{1}{2} & 0+0 
\end{array}
\end{eqnarray*}

\begin{eqnarray*}
\begin{array}{l|l|l}
k & \psi_2:\textup{order of }\lambda\textup{ in }\eqref{5.27} &
 -\deg C_{1,(k)}/\deg p_{alg} \\[1mm] \hline
1 & -\frac{1}{4}-\frac{5}{4} & 
3\cdot \frac{1}{2}+1-\frac{3}{2}=1 \\[1mm]
2 & \frac{1}{2}-\frac{1}{2}-\frac{3}{2} & 
3\cdot 1+0-\frac{3}{2}=\frac{3}{2} \\[1mm]
3 & \frac{1}{4}-\frac{3}{4} & 
3\cdot \frac{1}{2}+0-\frac{1}{2} = 1
\end{array}
\end{eqnarray*}

\begin{eqnarray}\label{10.13}
\det\begin{pmatrix}-3+\lambda & -2 \\ 3 & 2+\lambda\end{pmatrix}
=-\lambda(1-\lambda).
\end{eqnarray}

\bigskip
The case $\www E_8$:
\begin{eqnarray*}
\begin{array}{l|l|l|l}
k & \textup{involved }t_i& \rho:\textup{order of} 
& \psi_3:\textup{order of} \\ 
 & & \lambda\textup{ in \eqref{9.8}} & 
 \lambda\textup{ in \eqref{5.36}} \\ \hline
1=a_1 & t_9 & -\frac{1}{3} & 2
\\[1mm]
2=a_2=a_3 & t_8,t_7 & 0+\frac{1}{3} &  
1\textup{ (see \eqref{10.13})}\\
3=a_4=a_5 & t_6,t_5 & 0+1 & 0+0  \\
4=a_6=a_7 & t_4,t_3 & 0+\frac{2}{3} & 0+0 \\
5=a_8 & t_2 & \frac{1}{3} & 0 
\end{array}
\end{eqnarray*}

\begin{eqnarray*}
\begin{array}{l|l|l}
k & \psi_2:\textup{order of }\lambda\textup{ in }\eqref{5.33} &
 -\deg C_{1,(k)}/\deg p_{alg} \\[1mm] \hline
1 & -\frac{1}{2} & 
3\cdot \frac{-1}{3}+2-\frac{1}{2}=\frac{1}{2} \\[1mm]
2 & 0-1 & 
3\cdot \frac{1}{3}+1-1=1 \\[1mm]
3 & -\frac{1}{2}-\frac{3}{2} & 3\cdot 1+0-2 = 1 \\
4 & 0-1 & 3\cdot \frac{2}{3}+0-1=1 \\
5 & -\frac{1}{2} & 3\cdot\frac{1}{3}+0-\frac{1}{2}=\frac{1}{2}
\end{array}
\end{eqnarray*}

In all cases \eqref{10.12} holds. This and \eqref{10.11}
show \eqref{6.7}. This completes the proof of theorem
\ref{t6.3}.

\end{document}